\documentclass[toc=bib,toc=idx,abstract=false,parskip=half]{scrartcl}
\usepackage[utf8]{inputenc}
\usepackage[T1]{fontenc}
\usepackage[english]{babel}
\usepackage{makeidx}
\usepackage[columns=2,initsep=0pt]{idxlayout}
\usepackage{csquotes}
\usepackage{mathrsfs}
\usepackage{amsmath}
\usepackage{amsthm}
\usepackage{thmtools}
\usepackage{amsfonts}
\usepackage{amssymb}
\usepackage{stmaryrd}
\usepackage{mathtools}
\usepackage{mathdots}
\usepackage{enumitem}
\usepackage{calc}
\usepackage{graphicx}
\usepackage[dvipsnames]{xcolor}
\usepackage[style=alphabetic,citestyle=alphabetic,maxcitenames=4,maxalphanames=4,maxbibnames=10,backref=false,hyperref=true,backend=biber,doi=false,isbn=false,giveninits=true]{biblatex}
\usepackage[colorlinks=true,linkcolor=black,urlcolor=black,citecolor=black,bookmarksnumbered]{hyperref}
\usepackage{tikz}
\usetikzlibrary{arrows,decorations.markings,chains,calc,matrix}
\usepackage{gnuplot-lua-tikz}

\bibliography{oriented}
\DeclareFieldFormat{postnote}{#1}	
\DeclareFieldFormat{multipostnote}{#1}
\DeclareFieldFormat[article,eprint,unpublished,incollection]{title}{\mkbibemph{#1}}
\DeclareFieldFormat[article,eprint,unpublished,incollection]{journaltitle}{#1}
\DeclareFieldFormat[incollection]{booktitle}{#1}

\makeindex

\AtBeginDocument{%

}

\declaretheoremstyle[spaceabove=\topsep,spacebelow=0pt,bodyfont=\normalfont]{scdef}
\declaretheoremstyle[spaceabove=\topsep,spacebelow=0pt,bodyfont=\itshape]{scthm}
\declaretheoremstyle[spaceabove=\topsep,spacebelow=0pt,headfont=\normalfont\itshape,notefont=\normalfont\itshape,notebraces={}{},headformat={\NAME\NOTE},postheadspace=1em,qed=\qedsymbol]{scprf}

\declaretheorem[style=scthm,numberwithin=section,name=Theorem,    refname={Theorem,Theorems},        Refname={Theorem,Theorems}]        {Thm}
\declaretheorem[style=scthm,sharenumber=Thm,     name=Lemma,      refname={Lemma,Lemmas},            Refname={Lemma,Lemmas}]            {Lem}
\declaretheorem[style=scthm,sharenumber=Thm,     name=Corollary,  refname={Corollary,Corollaries},   Refname={Corollary,Corollaries}]   {Cor}
\declaretheorem[style=scthm,sharenumber=Thm,     name=Proposition,refname={Proposition,Propositions},Refname={Proposition,Propositions}]{Prop}
\declaretheorem[style=scdef,sharenumber=Thm,     name=Definition, refname={Definition,Definitions},  Refname={Definition,Definitions}]  {Def}
\declaretheorem[style=scdef,sharenumber=Thm,     name=Remark,     refname={Remark,Remarks},          Refname={Remark,Remarks}]          {Rem}
\declaretheorem[style=scdef,sharenumber=Thm,     name=Remarks,    refname={Remark,Remarks},          Refname={Remark,Remarks}]          {Rems}
\declaretheorem[style=scdef,sharenumber=Thm,     name=Example,    refname={Example,Examples},        Refname={Example,Examples}]        {Ex}
\declaretheorem[style=scprf,unnumbered,          name=Proof]{Prf}

\setlist[enumerate,1]{label={(\roman*)}}
\newcommand{\itemnr}[1]{(\romannumeral #1\relax)}

\DeclareTextFontCommand{\red}{\bfseries\color{red}}
\DeclareTextFontCommand{\green}{\bfseries\color{ForestGreen}}

\makeatletter\def\reallynopagebreak{\par\nopagebreak\@nobreaktrue}\makeatother

\DeclareMathOperator{\Ad}{Ad}
\DeclareMathOperator{\ad}{ad}
\DeclareMathOperator{\id}{id}

\DeclareMathOperator{\pos}{pos}
\DeclareMathOperator{\codim}{codim}
\DeclareMathOperator{\GL}{GL}
\DeclareMathOperator{\SL}{SL}
\DeclareMathOperator{\PSL}{PSL}
\DeclareMathOperator{\SO}{SO}
\DeclareMathOperator{\PSO}{PSO}

\DeclareMathOperator{\Gr}{Gr}
\DeclareMathOperator{\Gro}{Gr^+}

\DeclareMathOperator{\Hom}{Hom}
\DeclareMathOperator{\Span}{span}
\DeclareMathOperator{\sgn}{sgn}
\DeclareMathOperator{\ord}{ord}
\DeclareMathOperator{\Sym}{Sym}

\newcommand{\dc}[1]{\llbracket #1\rrbracket}

\newcommand{\F}{\mathcal F}
\newcommand{\dd}{\mathrm d}

\newcommand{\fg}{\mathfrak g}
\newcommand{\fk}{\mathfrak k}
\newcommand{\fp}{\mathfrak p}
\newcommand{\fa}{\mathfrak a}
\newcommand{\fu}{\mathfrak u}
\newcommand{\fb}{\mathfrak b}
\newcommand{\fn}{\mathfrak n}
\newcommand{\fsl}{\mathfrak{sl}}

\newcommand{\aplusbar}{\overline{\fa^+}}

\newcommand{\RP}{\mathbb{R}\mathrm{P}}

\newcommand{\bR}{\mathbb{R}}

\newcommand{\bZ}{\mathbb{Z}}
\newcommand{\bN}{\mathbb{N}}
\newcommand{\bH}{\mathbb{H}}

\newcommand{\equalo}{\overset{+}{=}}

\newcommand{\C}{\mathcal C}
\newcommand{\I}{\mathsf{v}}
\newcommand{\Q}{\mathcal{Q}}
\newcommand{\bdry}{\partial_\infty\Gamma}
\newcommand{\ctheta}{\Delta\mathord\setminus\theta}
\newcommand{\ceta}{\Delta\mathord\setminus\eta}
\newcommand{\Mbar}{{\,\overline{\!M}}}

\newcommand{\ival}[2]{(\! ( #1,#2)\! ) }

\def\svdots{\vbox{\baselineskip=4pt \lineskiplimit=0pt \kern2pt \hbox{.}\hbox{.}\hbox{.}\vspace{1pt}}}

\newenvironment{smatrix}{\left|\begin{smallmatrix}}{\end{smallmatrix}\right|}

\begin{document}
\title{Domains of discontinuity in oriented flag manifolds}
\author{Florian Stecker \and Nicolaus Treib \vphantom{\thanks{The authors acknowledge support from the Klaus Tschira Foundation, the European Research Council under ERC-Consolidator grant 614733, the RTG 2229 grant of the German Research Foundation, and the U.S. National Science Foundation grants DMS 1107452, 1107263, 1107367 ``RNMS: Geometric Structures and Representation Varieties'' (the GEAR Network).}}}
\deffootnote[0em]{0em}{0cm}{} 
\maketitle
\deffootnote[1em]{1.5em}{1em}{\textsuperscript{\thefootnotemark}} 
\begin{abstract}
  We study actions of discrete subgroups $\Gamma$ of semi--simple Lie groups $G$ on associated oriented flag manifolds. These are quotients $G/P$, where the subgroup $P$ lies between a parabolic subgroup and its identity component. For Anosov subgroups $\Gamma\subset G$, we identify domains in oriented flag manifolds by removing a set obtained from the limit set of $\Gamma$, and give a combinatorial description of proper discontinuity and cocompactness of these domains. This generalizes analogous results of Kapovich--Leeb--Porti to the oriented setting. We give first examples of cocompact domains of discontinuity which are not lifts of domains in unoriented flag manifolds. These include in particular domains in oriented Grassmannians for Hitchin representations, which we also show to be nonempty. As a further application of the oriented setup, we give a new lower bound on the number of connected components of $B$--Anosov representations of a closed surface group into $\SL(n,\bR)$.
\end{abstract}
\tableofcontents
\section{Introduction}
Since their introduction by Labourie \cite{Labourie}, Anosov representations have been established as a meaningful generalization of convex cocompact representations into Lie groups $G$ of higher rank. While the original definition was restricted to representations of surface groups into $\PSL(n,\bR)$, it has been extended to representations of any word hyperbolic group $\Gamma$ into any semi--simple Lie group $G$ in \cite{GuichardWienhardDomains}. Furthermore, the original more dynamical definition involving the geodesic flow has since been shown to be equivalent to several simpler descriptions in \cite{KapovichLeebPortiMorseActions}, \cite{KapovichLeebPortiMorseLemma}, \cite{GueritaudGuichardKasselWienhard} and \cite{BochiPotrieSambarino} (see \cite{KapovichLeebPortiCharacterizations} for an overview of different equivalent definitions).

A central feature of Anosov representations is the existence of a continuous, equivariant boundary map from the Gromov boundary $\bdry$ to a flag manifold $\F_\theta = G/P_\theta$. Here, $P_\theta$ is a parabolic subgroup defined by a subset $\theta \subset \Delta$ of the simple restricted roots, assumed to be invariant under the opposition involution.

Anosov representations have multiple desirable properties. They form an open subset of the representation variety, so any small deformation of an Anosov representation is still Anosov. Moreover, every Anosov representation is a quasiisometric embedding of $\Gamma$, equipped with the word metric, into $G$. The most interesting feature for us, however, is that they give rise to ``nice'' actions of $\Gamma$ on some flag manifolds $\F_\eta$, where $\eta\subset\Delta$ is another subset of the simple restricted roots. Such actions were investigated in \cite{GuichardWienhardDomains} and \cite{KapovichLeebPortiFlagManifolds}. For many choices of a flag manifold $\F_\eta$ they found cocompact domains of discontinuity $\Omega \subset \F_\eta$, open $\Gamma$--invariant subsets on which $\Gamma$ acts properly discontinuously with compact quotient $\Gamma \backslash \Omega$. Following the Ehresmann--Thurston principle, the existence of such domains allows to interpret Anosov representations as holonomy representations of geometric structures on the quotient. For example, the Hitchin component in $\PSL(3,\bR)$, for a closed surface $\Sigma$, is homeomorphic to the moduli space of (marked) convex projective structures on $\Sigma$ \cite{ChoiGoldmanRP2}.

In \cite{KapovichLeebPortiFlagManifolds} these domains are described by \emph{balanced ideals} or \emph{balanced thickenings}. A balanced ideal is a certain subset of the Weyl group $W$ of $G$. It specifies a submanifold for every point in $\xi(\bdry)$ which has to be removed from $\F_\eta$ to obtain a cocompact domain of discontinuity $\Omega \subset \F_\eta$.

For some choices of $\theta$ and $\eta$ there are no balanced ideals, so this construction does not produce any cocompact domains of discontinuity. But in many cases, for example if we are interested in geometric structures, it can also be useful to get a cocompact domain of discontinuity in some finite cover of $\F_\eta$ instead. This is the situation we investigate here.

\subsection{Domains of discontinuity}

Our main goal in this paper is to generalize the construction of cocompact domains of discontinuity from \cite{KapovichLeebPortiFlagManifolds} to actions on \emph{oriented flag manifolds}. Before giving more details on what we mean by that, let us start with a motivating example.

First, let $\rho \colon \Gamma \to \SL(3,\bR)$ be a Hitchin representation of the fundamental group $\Gamma$ of a closed surface $\Sigma$ with genus $g \geq 2$. These are continuous deformations of the holonomy representations $\Gamma \to \SO_0(2,1)$ of hyperbolic structures on $\Sigma$, composed with the standard embedding into $\SL(3,\bR)$. We consider the action of $\Gamma$ on $\RP^2$ defined by this representation. It is well--known that it preserves an open convex subset $\Omega \subset \RP^2$ (see \cite{ChoiGoldmanRP2}) and that the quotient of $\Omega$ by this action is homeomorphic to $\Sigma$, endowing $\Sigma$ with a so--called convex $\RP^2$ structure. In the simple case of a Fuchsian representation, i.e. if the image of $\rho$ is contained in $\SO_0(2,1)$, then $\Omega$ is an open disc, the Klein model of hyperbolic space.

As a Hitchin representation, $\rho$ is Anosov. Its limit map $\xi \colon \bdry \to \F_\varnothing$ maps to the complete flags comprising a point on the boundary of $\Omega$ and the line tangent to $\partial\Omega$ at this point. Note that $\Omega$ can be reconstructed from the limit map: It is the complement of the union of all these tangent lines.

\begin{figure}[htb]
  \centering
  \includegraphics[height=7cm]{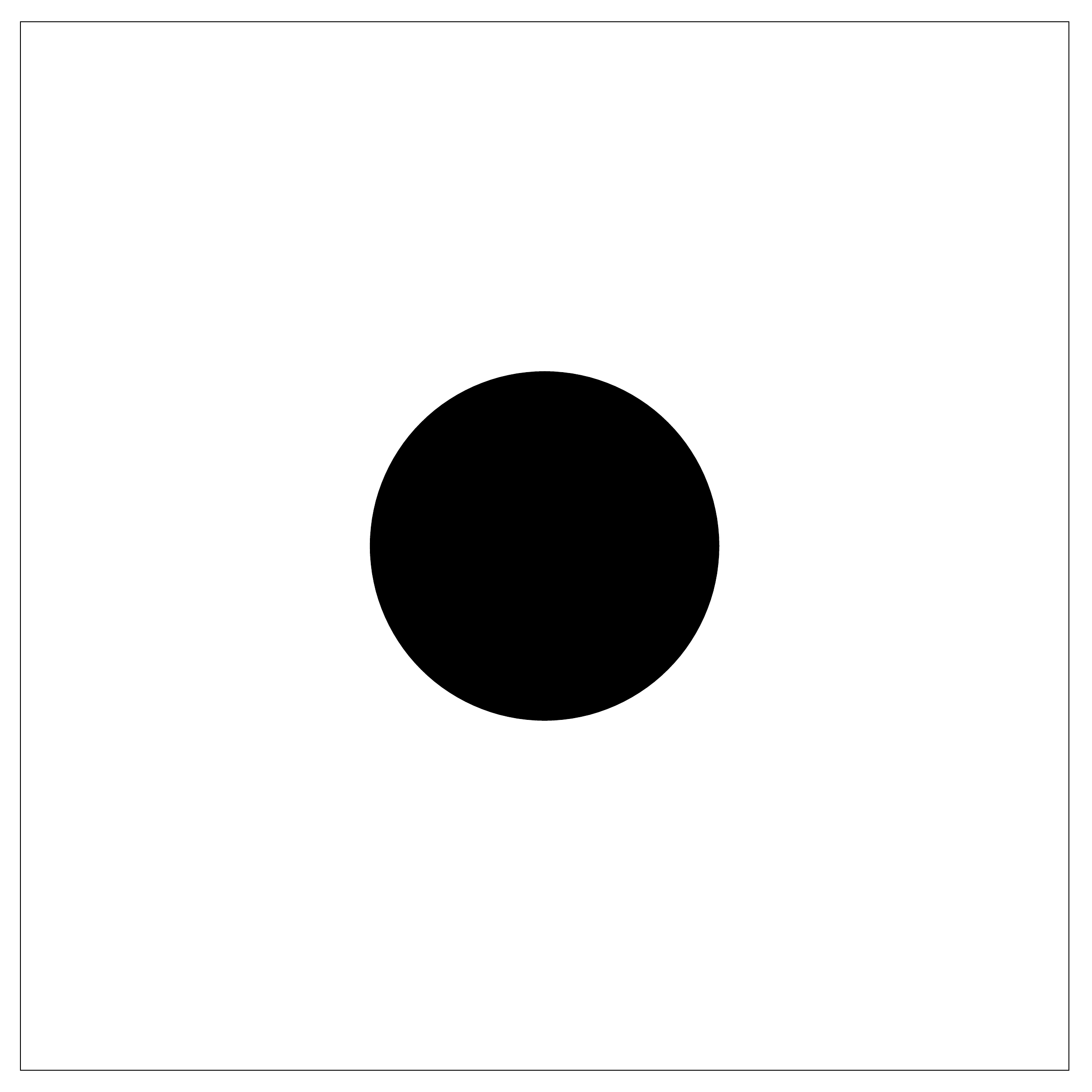}
  \includegraphics[height=7cm]{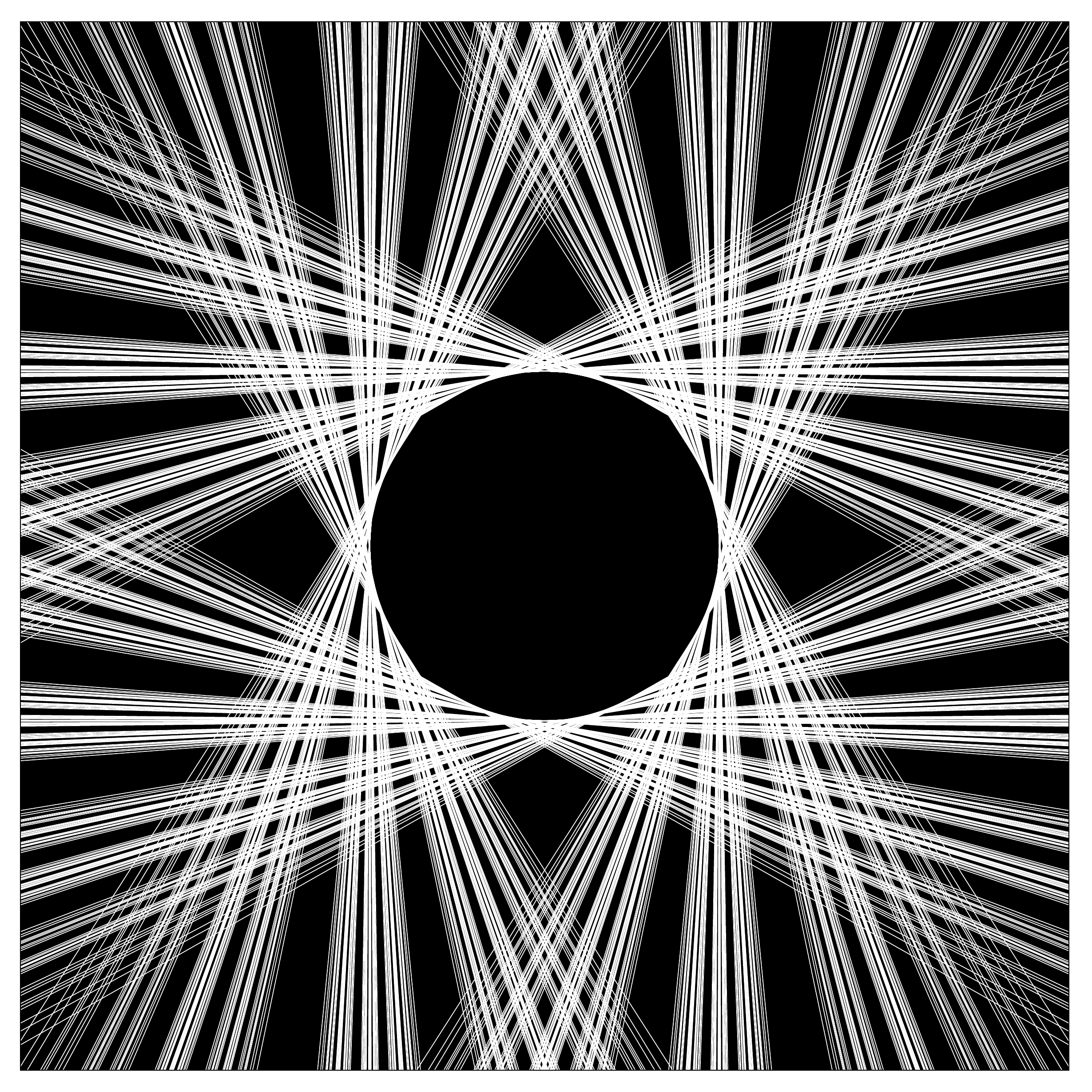}
  \caption{The domain of discontinuity (in black) in $\RP^2$ for a closed surface group $\Gamma$ (left) and a surface with boundary (right).}
  \label{fig:projective_pictures}
\end{figure}

If we take $\Sigma$ to be a surface with boundary, a few things behave differently. Now $\Gamma = \pi_1(\Sigma)$ is a free group. The holonomy representations $\Gamma \to \SO_0(2,1) \to \SL(3,\bR)$ of convex cocompact hyperbolic structures are again Anosov, as are small deformations, but they can lose the Anosov property if we deform too far. Assume that $\rho \colon \Gamma \to \SL(3,\bR)$ is close to a convex cocompact hyperbolic holonomy and therefore Anosov. The image of the limit map $\xi \colon \bdry \to \F_\varnothing$ is now a Cantor set of flags (see \autoref{fig:projective_pictures} right). The maximal domain of discontinuity for this action is once more obtained from the limit map by removing the union of these lines. It consists of a countable number of open quadrilaterals and a central ``big'' connected component containing hyperbolic space, whose quotient is a non--compact surface with open ends. This action has no cocompact domain of discontinuity in $\RP^2$ \cite{Stecker}.

\begin{figure}[htb]
  \centering
  \includegraphics[height=7cm]{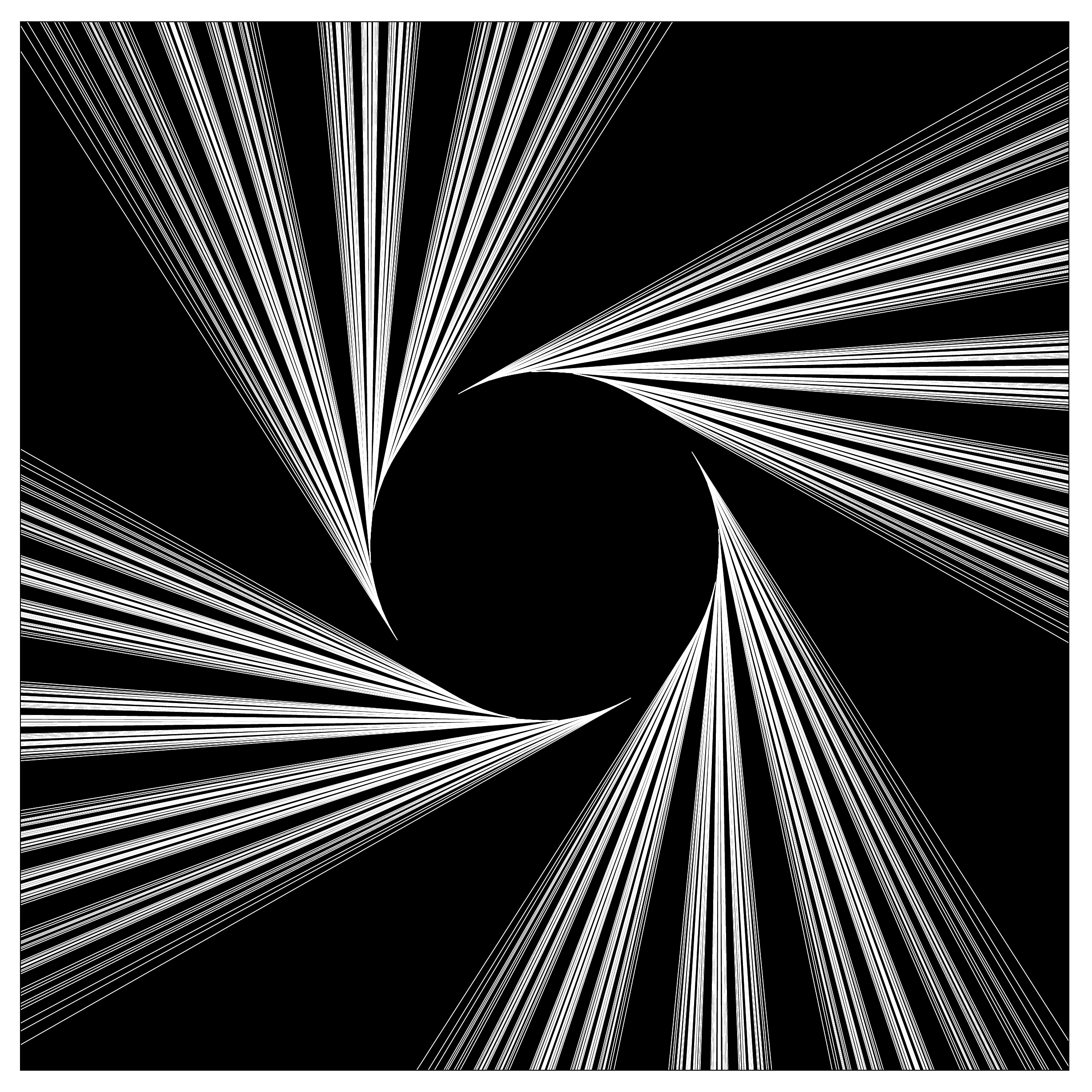}
  \begin{tikzpicture}[gnuplot]
    \draw [shading=ball,ball color=black] (3,3) circle (3);
    \path (-0.5,-0.5) -- (6,6);
    \input{images/schottky_group3d}
  \end{tikzpicture}
  \caption{The domain of discontinuity in oriented lines in $\bR^3$, shown in an affine chart (left) and in $S^2$ (right).}
  \label{fig:sphere}
\end{figure}

Surprisingly, a cocompact domain of discontinuity for this representation exists in the space of oriented lines, which is simply the double cover $S^2$ of $\RP^2$. This was first observed by Choi and can be found in \cite{ChoiGoldmanSpacetimes}. We can describe the domain as follows. First, there are two copies of the ``big'' component from the unoriented case. For every point in the limit set of these two copies, we have to remove half of the tangent great circle, where the choices are made consistently so the half circles are all disjoint (see \autoref{fig:sphere}). As a result, the two regions are joined by a countable number of ``strips'', one for every gap in the limit set. In the quotient, these strips become tubes connecting the two copies of the surface with open ends. It is therefore homeomorphic to the double of a surface with boundary and in particular compact.

As it turns out, this example is far from the only one where passing to a finite cover of the flag manifold leads to new cocompact domains of discontinuity. For instance, the above can be generalized to higher dimensions. Let $\rho \colon \Gamma \to \PSL(4n+3,\bR)$ be a small deformation of the composition of a convex cocompact representation into $\PSL(2,\bR)$ and the irreducible representation. Then there is a cocompact domain of discontinuity in $S^{4n+2}$. It is obtained by removing the spherical projectivizations of $2n\mathord{+}2$--dimensional (half--dimensional) half--spaces.

In fact, this construction also works if $\Gamma$ is a closed surface group and $\rho$ is a Hitchin representation. In $S^2$ it will only give the two lifts of the domain in $\RP^2$, but in higher dimensional $S^{4n+2}$ we get cocompact domains which do not exist in the unoriented case. This family of examples was independently found by Danciger, Gu\'eritaud and Kassel.

In this paper, we try to gain a better understanding of this behavior by generalizing the construction in \cite{KapovichLeebPortiFlagManifolds} to actions on oriented flag manifolds. The above examples will be special cases of this. An example of an oriented flag (for $\SL(n,\bR)$) is a sequence of nested subspaces with orientations assigned to some of these spaces. More generally, for a semi--simple Lie group $G$, we consider as oriented flag manifolds the homogeneous spaces $G/P$ where $P$ is a proper subgroup of $G$ containing the identity component $B_0$ of a minimal parabolic subgroup.

Unfortunately, this inconspicuous change leads to quite a bit of additional setup. We have to replace the Weyl group $W = N_K(\fa) / Z_K(\fa)$ by a different group, the ``extended Weyl group'' $\widetilde W = N_K(\fa)/Z_K(\fa)_0$. Then the oriented flag manifolds $\F_R$ are parametrized by certain subgroups $R\subset \widetilde W$ which we call their \emph{oriented parabolic types} (\autoref{def:types}). In order to construct cocompact domains of discontinuity in oriented flag manifolds, the limit map of an Anosov representation needs to lift to some oriented flag manifold. Let us assume that $\F_R = G/P_R$ is a finite cover of the unoriented flag manifold $\F_\theta$. Then we call a $P_\theta$--Anosov representation \emph{$P_R$--Anosov} if its limit map $\xi$ lifts continuously and equivariantly to $\F_R$. There is in fact a unique maximal choice of such an oriented flag manifold $\F_R$ that one can lift $\xi$ to (\autoref{prop:minimal_type}). To such a lift, we can associate its \emph{transversality type}, which is the relative position of two flags in its image (\autoref{def:relative_position}). It can be represented by an element $w_0 \in \widetilde W$ which is a lift of the longest element in $W$.

Finally, there is a natural partial order on $\widetilde W$ which is an analogue of the Bruhat order on the Weyl group, and we can define the notion of a balanced ideal also in this setting. An ideal $I \subset \widetilde W$ is a subset such that whenever $w \in I$ and $w' \leq w$, then $w' \in I$. The term \emph{balanced} refers to the action of the transversality type $w_0 \in \widetilde W$ on $\widetilde W$: This is an order--reversing involution, and $I$ is balanced if every orbit has two elements, exactly one of which is contained in $I$. Our main theorem is then

\begin{Thm}	\label{thm:domains_introduction}
  Let $\Gamma$ be a non--elementary word--hyperbolic group and $G$ a connected, semi--simple, linear Lie group. Let $\rho \colon \Gamma \to G$ be a $P_\theta$--Anosov representation with limit map $\xi \colon \bdry \to \F_\theta$ which admits a $\rho$--equivariant continuous lift $\widehat{\xi} \colon \bdry \to \F_R$. Furthermore, let $w_0 \in \widetilde W$ represent the transversality type of $\widehat{\xi}$ and $I\subset \widetilde W$ be a $w_0$--balanced, $R$--left invariant and $S$--right invariant ideal, where $S\subset \widetilde W$ is another oriented parabolic type. Define $\mathcal{K} \subset \F_S$ by
  \[ \mathcal{K} = \bigcup_{x\in\bdry} \bigcup_{w\in I} \{ [gw] \in \F_S \mid g\in G, [g] = \widehat{\xi}(x) \}. \]
  Then $\mathcal{K}$ is $\Gamma$--invariant and closed, and $\Gamma$ acts properly discontinuously and cocompactly on the domain $\Omega = \F_S \setminus \mathcal{K}$.
\end{Thm}

If we take $\F_R$ and $\F_S$ to be unoriented flag manifolds (so in particular $\F_R = \F_\theta$), then we recover the original theorem from \cite{KapovichLeebPortiFlagManifolds}. The $R$--left and $S$--right invariant balanced ideals in $\widetilde W$ are in this case just lifts of balanced ideals in the Weyl group $W$.

An application of this theorem concerns the action of Hitchin representations on Grassmannians. If $\rho \colon \Gamma \to \PSL(n,\bR)$ is a Hitchin representation, it defines an action of $\Gamma$ on the Grassmannian $\Gr(k,n)$ of $k$--dimensional subspaces of $\bR^n$. If $n$ is even and $k$ is odd, then there is a balanced ideal in $W$ and therefore a cocompact domain of discontinuity $\Omega \subset \Gr(k,n)$. For odd $n \geq 5$, no cocompact domains exist in any $\Gr(k,n)$ \cite{Stecker}. However, we find oriented Grassmannians admitting cocompact domains of discontinuity in these cases. More precisely, there is a cocompact domain of discontinuity in the oriented Grassmannian $\Gro(k,n)$ whenever $n$ is odd and $k(n+k+2)/2$ is even (see \autoref{prop:grassmannians}).

Note that, like \cite{KapovichLeebPortiFlagManifolds}, we are only concerned with domains of discontinuity which are open subsets of the oriented flag manifold. In contrast, there are recent works by Danciger, Gueritaud and Kassel \cite{DancigerGueritaudKasselProjective,DancigerGueritaudKasselPseudoRiemannian} and Zimmer \cite{Zimmer} constructing cocompact domains of discontinuity in $\RP^n$ and in $\mathbb H^{p,q}$ and proving that their existence is equivalent to the Anosov property. These domains are closed subsets.

\subsection{Connected components of Anosov representations}

Another use of the theory we develop concerns connected components of Anosov representations, based on the following proposition.

\begin{Prop}[\autoref{prop:clopen}]
  The set of $P_R$--Anosov representations is open and closed in the space of all $P_\theta$--Anosov representations.
\end{Prop}

Consequently, this notion can be used to distinguish connected components of Anosov representations: If two $P_\theta$--Anosov representations lie in the same connected component, then both the maximal choice of $\F_R$ that their boundary maps can be lifted to and the transversality type of their boundary maps must agree (\autoref{cor:connected_components}). We apply this fact in \autoref{sec:block_embeddings} to $B$--Anosov representations of closed surface groups into $\SL(n,\bR)$ for odd $n$. Namely, we consider block embeddings constructed by composing a Fuchsian representation into $\SL(2,\bR)$ with irreducible representations into $\SL(k,\bR)$ and $\SL(n-k,\bR)$. For different choices of block sizes, we show that these representations lie in different connected components of $B$--Anosov representations. Together with an observation by Thierry Barbot and Jaejeong Lee, which is explained in \cite[Section 4.1]{KimKuessner}, this leads to the following lower bound for the number of connected components.
\begin{Prop}[\autoref{cor:number_of_components}] \label{prop:connected_components}
  Let $\Gamma=\pi_1(\Sigma)$, where $\Sigma$ is a closed surface of genus $g\geq 2$, and let $n$ be odd. Then the space $\Hom_{\textnormal{$B$--Anosov}}(\Gamma,\SL(n,\bR))$ has at least $2^{2g-1}(n-1) +1$ connected components.
\end{Prop}

\subsection{Organization of the paper}

In \autoref{sec:parabolics}, we recall basic facts about parabolic subgroups and the Weyl group and fix the notation we use. We also give a refined version of the Bruhat decomposition.

In \autoref{sec:oriented_relative_positions}, we introduce oriented flag manifolds and explain how to describe relative positions of oriented flags in terms of the group $\widetilde W$. Moreover, we introduce the Bruhat order on $\widetilde W$ and derive a combinatorial description of it in \autoref{prop:Bruhat_order_quotient}. We also deal with the fact that in $\widetilde W$, there is no unique longest element, and multiple choices can be valid.

In \autoref{sec:orientable_Anosov}, we define $P_R$--Anosov representations and examine their properties. In particular, we show that they form a union of connected components of the set of $P_\theta$--Anosov representations, where $\theta$ is the subset of simple restricted roots contained in $R$. We also prove that starting with a $P_\theta$--Anosov representation, there is a unique maximal choice of oriented flag manifold $\F_R$ that one can lift its boundary map to.

In \autoref{sec:domains_of_discontinuity}, we prove \autoref{thm:domains_introduction}.

In \autoref{sec:examples_of_balanced_ideals}, we give some examples of balanced ideals for the group $\SL(n,\bR)$. In particular, we give a complete description of the Bruhat order on the extended Weyl group of $\SL(3,\bR)$, contrast it with the unoriented case, and list the balanced ideals it admits. Furthermore, we give a geometric interpretation of oriented relative positions.

We turn to applications in \autoref{sec:examples_of_representations}. One such application is a list of oriented Grassmannians admitting cocompact domains of discontinuity for Hitchin representations. We then discuss generalized Schottky representations, obtaining as a special case the domain from our motivating example (\autoref{fig:sphere}).

\autoref{sec:block_embeddings} concerns connected components of $B$--Anosov representations in $\SL(n,\bR)$, for odd $n$. We give a lower bound on the number of connected components by showing that certain block embeddings of Fuchsian representations must land in different components, proving \autoref{prop:connected_components}.

A number of technical results were included as an appendix to improve readability: \autoref{sec:B0B0action} is concerned with orbits and orbit closures of the $B_0\times B_0$--action on $G$. A combinatorial description of such orbit closures with respect to the refined Bruhat decomposition is derived. \autoref{sec:compact_homogeneous} contains the results necessary to relate expansivity of a group action to cocompactness. This relation is a key component of the cocompactness part of \autoref{thm:domains_introduction}.

We thank Jean--Philippe Burelle, Jeff Danciger and Anna Wienhard for helpful discussions and suggestions. We are greatly indebted to Beatrice Pozzetti and Johannes Horn for pointing out mistakes in preliminary versions and for suggestions that helped improve the structure of this paper.
\section{Parabolic subgroups}\label{sec:parabolics}
Let $G$ be a connected semi--simple Lie group with finite center and $\fg$ its Lie algebra. Choose a maximal compact subgroup $K \subset G$. Let $\fk \subset \fg$ be the Lie algebra of $K$ and $\fp = \fk^\perp$ its orthogonal complement with respect to the Killing form. Choose a maximal abelian subalgebra $\fa$ in $\fp$ and let $\fa^*$ be its dual space.\index{a@$\fa$ maxmal abelian subalgebra of $\fp$}\index{p@$\fp = \fk^\perp$}\index{k@$\fk$ Lie algebra of $K$}\index{K@$K$ maximal compact subgroup of $G$}

For any $\alpha \in \fa^*$ let
\[\fg_\alpha = \{X \in \fg \mid \forall H \in \fa \colon [H,X] = \alpha(H) X\}\]
and let $\Sigma \subset \fa^*$ be the set of restricted roots, that is the set of $\alpha \neq 0$ such that $\fg_\alpha \neq 0$. $\Sigma$ is in general not reduced, i.e. there can be $\alpha \in \Sigma$ with $2\alpha \in \Sigma$ (but no other positive multiples except $2$ or $1/2$). Choose a simple system $\Delta \subset \Sigma$ (a basis of $\fa^*$ such that every element of $\Sigma$ can be written as a linear combination in $\Delta$ with only non--negative or only non--positive integer coefficients). Let $\Sigma^\pm$ be the corresponding positive and negative roots, $\Sigma_0$ the indivisible roots (the roots $\alpha \in \Sigma$ with $\alpha/2 \not\in\Sigma$) and let $\Sigma_0^\pm = \Sigma^\pm \cap \Sigma_0$. Let $\fa^+ = \{X \in \fa \mid \alpha(X) > 0 \,\forall \alpha \in \Delta \}$ and $\aplusbar = \{X \in \fa \mid \alpha(X) \geq 0 \,\forall \alpha \in \Delta \}$.\index{g_a@$\fg_\alpha$ root space}\index{S_@$\Sigma$ set of restricted roots}\index{D_@$\Delta$ simple restricted roots}\index{S_0@$\Sigma_0$ indivisible restricted roots}\index{S_+-@$\Sigma^\pm$ positive/negative restricted roots}\index{S_0+-@$\Sigma_0^\pm = \Sigma_0 \cap \Sigma^\pm$}\index{a+@$\fa^+ = \{X \in \fa \mid \alpha(X) > 0 \,\forall \alpha \in \Delta \}$}

The above choices are always possible and such a triple $(K, \fa, \Delta)$ is unique up to conjugation in $G$ (see e.g. \cite[Theorem 2.1]{Helgason} and \cite[Theorems 2.63, 6.51, 6.57]{Knapp}).

For any proper subset $\theta \subsetneq \Delta$ define the Lie algebras
\[\fn = \bigoplus_{\alpha \in \Sigma^+} \fg_\alpha, \qquad \fn^- = \bigoplus_{\alpha \in \Sigma^-} \fg_\alpha, \qquad \fb = \bigoplus_{\alpha \in \Sigma^+ \cup \{0\}} \fg_\alpha, \qquad \fp_\theta = \!\!\! \bigoplus_{\alpha \in \Sigma^+ \cup \,\Span(\theta)} \!\!\! \fg_\alpha.\]
The $\fp_\theta \subset \fg$ are the standard parabolic subalgebras and $\fb = \fp_\varnothing$ is the minimal standard parabolic subalgebra. A parabolic subalgebra is a subalgebra which is conjugate to a standard parabolic subalgebra.\index{n@$\fn,\fn^-$ nilpotent subalgebras}\index{b@$\fb$ minimal parabolic subalgebra}\index{pt_@$\fp_\theta$ parabolic subalgebra}

Let $A, N, N^- \subset G$ be the connected subgroups with Lie algebras $\fa, \fn, \fn^-$. The exponential map of $G$ restricts to diffeomorphisms $\exp \colon \fa \to A$, $\exp \colon \fn \to N$ and $\exp \colon \fn^- \to N^-$. Let $P_\theta = N_G(\fp_\theta) \subset G$ be the (standard) parabolic subgroups and $B = P_\varnothing$ the (standard) minimal parabolic subgroup. The Lie algebra of $P_\theta$ is $\fp_\theta$. A parabolic subgroup is typically defined to be a subgroup conjugate to some $P_\theta$. However, when we write ``parabolic subgroup'' here, we will just mean $P_\theta$ for some $\theta \subsetneq \Delta$. Using the Iwasawa decomposition $G = KAN$, the minimal parabolic can also be described as $B = Z_K(\fa)AN$, with $Z_K(\fa)$ being the centralizer of $\fa$ in $K$.\index{N@$N,N^-$ unipotent subgroups of $G$}\index{A@$A$ subgroup with Lie algebra $\fa$}\index{Pt_@$P_\theta$ parabolic subgroup}\index{B@$B$ minimal parabolic subgroup}\index{B0@$B_0$ identity component of $B$}

The quotients
\[ \F_\theta := G/P_\theta \]
are compact $G$--homogeneous spaces and are called (partial) flag manifolds.

Define the groups\index{W@$W = N_K(\fa)/Z_K(\fa)$ Weyl group}\index{W_@$\widetilde W = N_K(\fa)/Z_K(\fa)_0$ ext. Weyl group}\index{Mbar@$\Mbar = Z_K(\fa)/Z_K(\fa)_0$}
\[W = N_K(\fa)/Z_K(\fa), \quad \widetilde W = N_K(\fa)/Z_K(\fa)_0, \quad \Mbar = Z_K(\fa)/Z_K(\fa)_0,\]
where $Z_K(\fa)_0$ is the identity component of the centralizer of $\fa$ in $K$. The group $W$ is the Weyl group of $G$, and we call $\widetilde W$ the extended Weyl group. The set of simple roots $\Delta$ can be realized as a subset of $W$ by identifying $\alpha \in \Delta$ with the Killing orthogonal reflection along $\ker\alpha$ in $\fa$. Then $\Delta \subset W$ is a generating set. We denote by $\ell \colon W \to \bN$ the word length with respect to the generating set $\Delta$. There is a unique longest element $w_0 \in W$. Conjugation by $w_0$ preserves $\Delta$ and defines a map $\iota \colon \Delta \to \Delta$ which is called the opposition involution. If we write $\ell(w)$ for $w \in \widetilde W$ we mean the length $\ell(\pi(w))$ of its projection to $W$.\index{l@$\ell$ word length in $W$}\index{i_@$\iota$ opposition involution}

Though the groups $W$, $\widetilde W$ and $\Mbar$ are not really subgroups of $G$, we can often simplify notation by pretending they are. For example, if $H \subset G$ is a subgroup containing $Z_K(\fa)_0$, we will write $H \cap \Mbar$ as a shorthand for $(H \cap Z_K(\fa)) / Z_K(\fa)_0 \subset \Mbar$, or $wH$ for the coset $nH$ where $n \in N_K(\fa)$ is any choice of representative for $w \in \widetilde W$.

See \autoref{sec:SLnR} for a description of $\Delta$, $\widetilde W$, $\Mbar$ etc. in the case $G = \SL(n,\bR)$.

An central feature of $W$ for us is the Bruhat decomposition \cite[Theorem 7.40]{Knapp}
\[G = \bigsqcup_{w\in W} BwB.\]
Analogous to this, we have the following statement for $B_0$ and $\widetilde W$. A proof is given as \autoref{prop:bruhat_appendix} in the appendix.
\begin{Prop}[Refined Bruhat decomposition]\label{prop:bruhat}
  $G$ decomposes into the disjoint union
  \[G = \bigsqcup_{w \in \widetilde W} B_0 w B_0.\]
\end{Prop}
We write $L_g$ and $R_g$ for left and right multiplication by $g \in G$. For subsets $T_1, \dots, T_n$ of a group we write $\langle T_1, \dots, T_n \rangle$ for the smallest subgroup containing all of them. When speaking about quotients of groups, we usually denote equivalence classes by square brackets. To avoid confusion, we use double brackets $\dc{\ }$ for all quotients of the group $\widetilde W$. It should always be clear from the context which quotient we are referring to.\index{Lg@$L_g$ left multiplication by $g$}\index{Rg@$R_g$ right multiplication by $g$}\index{<>@$\langle T \rangle$ group generated by $T$}\index{[]@$\dc\cdot$ equivalence class in $\widetilde W$}

Our attention in this paper will be restricted to semi--simple Lie groups $G$ such that the group $\Mbar$ is finite abelian and consists entirely of involutions. This holds for all $G$ which are linear, i.e. isomorphic to a closed subgroup of some $\GL(n,\bR)$ (see \cite[Theorem 7.53]{Knapp} and note that all connected linear Lie groups have a complexification). Also, every linear semi--simple Lie group has a finite center \cite[Proposition 7.9]{Knapp}. All our arguments work equally well for Lie groups which are not linear, as long as their center is finite and $\Mbar$ is finite abelian and consists of involutions. These assumptions on $\Mbar$ do not appear to be essential for our theory, but they significantly simplify several arguments, e.g. the statement and proof of \autoref{lem:doublecoset_generator}.

If $\Mbar$ is trivial, the class of oriented flag manifolds (see \autoref{def:oriented_flag_manifold}) reduces to ordinary flag manifolds. Our theory gives nothing new in this case. In particular, this happens whenever $G$ is a complex Lie group, since their minimal parabolic subgroups are connected.
\section{Oriented relative positions}	\label{sec:oriented_relative_positions}
\subsection{Oriented flag manifolds} \label{sec:oriented_flag_manifolds}
Let $B$ be the minimal parabolic subgroup as defined \autoref{sec:parabolics} and $B_0$ its identity component. Note that a proper closed subgroup $B_0 \subset P \subsetneq G$ containing $B_0$ has a parabolic Lie algebra and is thus a union of connected components of a parabolic subgroup.
\begin{Def}	\label{def:oriented_flag_manifold}
  Let $B_0 \subset P \subsetneq G$ be a proper closed subgroup containing $B_0$. We call such a group \emph{(standard) oriented parabolic subgroup} and the quotient $G/P$ an \emph{oriented flag manifold}.
\end{Def}
\begin{Ex}
  Let $G = \SL(n,\bR)$ be the special linear group. Then $B_0$ is the set of upper triangular matrices with positive diagonal entries $\lambda_1,\ldots,\lambda_n$. The space $G/B_0$ can be identified with the space of complete oriented flags, i.e. complete flags with a choice of orientation in every dimension. An example of a closed subgroup $B_0 \subsetneq P \subsetneq G$ is the group of upper triangular matrices where $\lambda_1$ and $\lambda_2$ are allowed to be negative, while the remaining entries are positive. The space $G/P$ identifies with the space of complete flags with a choice of orientation on every component except the $1$-dimensional one. In this way, all partial flag manifolds with a choice of orientation on a subset of the components of the flags can be obtained. However, we can also consider e.g. the group $P' = \langle B_0, -1 \rangle$ if $n$ is even. Its corresponding oriented flag manifold $G/P'$ is the space of complete oriented flags up to simultaneously changing the orientation on every odd-dimensional component.
\end{Ex}
The parabolic subgroups of $G$ are parametrized by proper subsets $\theta$ of $\Delta$. We want a similar description for oriented parabolics. To define this, we first need a lift $\I \colon \Delta \to \widetilde W$ of our generating set $\Delta \subset W$, which behaves well in the sense that $\I(\alpha) \in (P_\alpha)_0$ for all $\alpha \in \Delta$ (see \autoref{lem:Ptheta_cells}).

The construction requires some Lie theory. If $\langle \cdot, \cdot \rangle$ is the Killing form on $\fg$ and $\Theta$ the Cartan involution which is $1$ on $\fk$ and $-1$ on $\fp$, then $\|X\|^2 = - \langle X, \Theta X \rangle$ defines a norm on $\fg$. Its restriction to $\fa$ is just $\langle X,X \rangle$. For $\alpha \in \fa^*$ let $H_\alpha \in \fa$ be its dual with respect to $\langle\cdot,\cdot\rangle$, i.e. $\langle H_\alpha, X \rangle = \alpha(X)$ for all $X \in \fa$. We use the norm on $\fa^*$ defined by $\|\alpha\|^2 = \langle H_\alpha, H_\alpha \rangle$.\index{T_@$\Theta$ Cartan involution on $\fg$}\index{Ha_@$H_\alpha$ dual of $\alpha$}

\begin{Def}\label{def:generators}
  For every $\alpha \in \Delta$ choose a vector $E_\alpha \in \fg_\alpha$ such that $\|E_\alpha\|^2 = 2 \|\alpha\|^{-2}$. Then we define
  \[\I(\alpha) = \exp \left(\frac{\pi}{2}\left(E_\alpha + \Theta E_\alpha\right)\right) \!.\]
  By \cite[Proposition 6.52(c)]{Knapp} this is in $N_K(\fa)$ and acts on $\fa$ as a reflection along the kernel of the simple root $\alpha$. We will regard $\I(\alpha)$ as an element of $\widetilde W = N_K(\fa)/Z_K(\fa)_0$.
\end{Def}\index{v(a_)@$\I(\alpha)$ lift of $\alpha$ to $\widetilde W$}
\begin{Rems}\label{rem:choice_of_r}\ \reallynopagebreak
  \begin{enumerate}
  \item $\I(\alpha) \in \widetilde W$ is almost independent of the choice of $E_\alpha$: If $\dim \fg_\alpha > 1$, then the set of admissible $E_\alpha$ is connected. Since $\widetilde W$ is discrete and $\I(\alpha)$ depends continuously on $E_\alpha$ this means that $\I(\alpha) \in \widetilde W$ is independent of $E_\alpha$. In particular, we get the same $\I(\alpha)$ when substituting $E_\alpha$ by $- E_\alpha$, so $\I(\alpha) = \I(\alpha)^{-1}$. On the other hand, if $\dim \fg_\alpha = 1$ then $\I(\alpha)$ need not be of order $2$, and there can be two different choices for $\I(\alpha)$, which are inverses of each other. If they do not coincide, $\I(\alpha)$ is of order $4$ since $\I(\alpha)^2$ acts trivially on $\fa$ and is therefore contained in $\Mbar$. By our assumption of $G$ being linear, $\Mbar$ consists of involutions. In the group $\SL(n,\bR)$ for example, $\I(\alpha)$ is of order $4$ for all simple restricted roots $\alpha$, while in $\SO_0(p,q), \ p < q,$ the image of the ``last'' simple root $\alpha_p$ is of order $2$.
  \item For every $\alpha \in \Delta$, we have $\pi(\I(\alpha)) = \alpha \in W$, so the projection of $\I(\Delta)$ to $W$ is just the usual generating set $\Delta$. In fact, $\I(\Delta)$ also generates the group $\widetilde W$ (this can be seen using \autoref{lem:Ptheta_cells} with $\theta = \Delta$).
  \end{enumerate}
\end{Rems}
Now we can define our objects parametrizing oriented parabolic subgroups.
\begin{Def}\label{def:types}
  Let $\theta \subsetneq \Delta$ and $\Mbar_\theta = \langle \I(\theta) \rangle \cap \Mbar$. Let $\Mbar_\theta \subset E \subset \Mbar$ be a subgroup. Then we call the group $R = \langle \I(\theta), E \rangle \subset \widetilde W$ an \emph{oriented parabolic type}.\index{Mbart_@$\Mbar_\theta = \langle \I(\theta) \rangle \cap \Mbar$}\index{t_@$\theta,\eta$ subsets of $\Delta$}\index{R@$R \kern-0.5mm = \langle \I(\theta), E \rangle$ oriented parabolic type}\index{S@$S = \langle \I(\eta), F \rangle$ oriented parabolic type}\index{E@$E, F$ subgroups of $\Mbar$}
\end{Def}
\begin{Rems}\ \reallynopagebreak	\label{rem:types}
  \begin{enumerate}
  \item This definition does not depend on the choices involved in $\I$ (see \autoref{rem:choice_of_r}\itemnr{1}).
  \item For every oriented parabolic type $R$, there is a unique pair $(\theta, E)$ with $\theta \subsetneq \Delta$, $\Mbar_\theta \subset E \subset \Mbar$, and $R = \langle \I(\theta), E \rangle$. In fact, using \autoref{lem:Enormal} below, we can recover $\theta$ and $E$ from $R$ by
    \[R \cap \Mbar = \langle \I(\theta) \rangle E \cap \Mbar = \Mbar_\theta E = E\]
    and
    \[\pi(R) \cap \Delta = \pi(\langle \I(\theta) \rangle) \cap \Delta = \langle \theta \rangle \cap \Delta = \theta\]
    where $\pi$ is the projection from $\widetilde W$ to $W$.
  \end{enumerate}
\end{Rems}
\begin{Prop}\label{prop:parabolics}
  The map
  \[\{\text{oriented parabolic types}\} \to \{\text{oriented parabolic subgroups}\}\]
  mapping $R$ to $P_R = B_0 R B_0$ is a bijection. Its inverse maps $P$ to $P \cap \widetilde W$. We will call $P \cap \widetilde W$ the \emph{type of $P$}.
\end{Prop}\index{PR@$P_R$ oriented parabolic of type $R$}
\begin{Def}
  Let $P_R$ be the oriented parabolic of type $R = \langle \I(\theta), E \rangle$. Then we write
  \[\F_R = G/P_R, \qquad \F_\theta = G/P_\theta\]
  for the associated oriented and unoriented flag manifolds.\index{FR@$\F_R = G/P_R$ oriented flag manifold}\index{Ft_@$\F_\theta = G/P_\theta$ unoriented flag manifold}
\end{Def}
The remainder of \autoref{sec:oriented_flag_manifolds} is a proof of \autoref{prop:parabolics}. We first need a few lemmas.

\begin{Lem}\label{lem:Enormal}
  Let $\alpha \in \Delta$ and $w \in \widetilde W$ such that $\pi(w)$ and $\alpha$ are commuting elements of $W$. Then $w \, \I(\alpha) \, w^{-1} \in \{\I(\alpha), \I(\alpha)^{-1}\} \subset \widetilde W$. In particular, this holds for any $w \in \Mbar$. As a consequence, for any $\theta \subset \Delta$ and any subgroup $E \subset \Mbar$
  \[\langle \I(\theta), E \rangle = \langle \I(\theta) \rangle \, E = E \, \langle \I(\theta) \rangle.\]
\end{Lem}
\begin{Prf}
  We compute, using that $\Ad_w$ commutes with the Cartan involution,
  \[\textstyle w\,\I(\alpha)\,w^{-1} = \exp \left(\frac{\pi}{2} \left(\Ad_w E_\alpha + \Theta \Ad_w E_\alpha\right)\right)\!.\]
  Since $\Ad_w$ preserves $\|\cdot\|$ and the root $\alpha$ is preserved by $w$ this just corresponds to a different choice of $E_\alpha \in \fg_\alpha$ in the definition of $\I(\alpha)$, so $w\,\I(\alpha)\,w^{-1}$ must be either $\I(\alpha)$ or $\I(\alpha)^{-1}$ by \autoref{rem:choice_of_r}\itemnr{1}. So in particular $m \, \langle \I(\theta) \rangle \, m^{-1} \subset \langle \I(\theta) \rangle$ for any $m \in \Mbar$ and $\theta \subset \Delta$, which shows the second statement.
\end{Prf}
\begin{Lem}\label{lem:B_0R}
  Let $R,S$ be oriented parabolic types and $w \in \widetilde W$. Then
  \[B_0 R B_0 w B_0 S B_0 = B_0 R w S B_0.\]
\end{Lem}
\begin{Prf}
  Let $R = \langle \I(\theta), E \rangle$. We first prove $B_0 w' B_0 w B_0 \subset B_0 R w B_0$ for all $w \in \widetilde W$ and $w' \in R$ by induction on $\ell(w')$. If $\ell(w') = 0$, then $w' \in \Mbar$, so $B_0 w' B_0 w B_0 = B_0 w' w B_0 \subset B_0 R w B_0$. If $\ell(w') > 0$ then we can find $\alpha \in \theta$ and $s = \I(\alpha)$ with $w' = w'' s$ and $\ell(w') = \ell(w'') + 1$. So by \autoref{lem:doublecoset_additivity2}
  \begin{align*}
    B_0 w' B_0 w B_0 &= B_0 w'' s B_0 w B_0 = B_0  w'' B_0 s B_0 w B_0 \\
                     &\subset B_0 w'' B_0 w B_0 \cup B_0 w'' B_0 s w B_0 \cup B_0 w'' B_0 s^2 w B_0,
  \end{align*}
  which is in $B_0 R w B_0$ by the induction hypothesis, since $s,w'' \in R$. So $B_0 R B_0 w B_0 \subset B_0 R w B_0$. By the same argument $B_0 S B_0 w^{-1} B_0 \subset B_0 S w^{-1} B_0$, and thus $B_0 w B_0 S B_0 \subset B_0 w S B_0$. Together, this shows the lemma.
\end{Prf}
\begin{Lem}\label{lem:PR_decomposition}
  If $R$ is an oriented parabolic type, then $B_0 R B_0$ is a closed subgroup of $G$.
\end{Lem}
\begin{Prf}
  Closedness follows from \autoref{prop:combinatorial_bruhat_order} as $A_w \subset R$ for every $w \in R$. This is because we can write $w = w' m$ with $w' \in \langle \I(\theta) \rangle$ and $m \in E$, and then $A_w = A_{w'}m \subset R$. To see that $B_0 R B_0$ is a subgroup we take $w, w' \in R$ and need to prove that $B_0 w B_0 w' B_0 \subset B_0 R B_0$. But this follows from \autoref{lem:B_0R} (with $S = 1$)
\end{Prf}
\begin{Lem}\label{lem:Ptheta_cells}
  Let $\theta \subsetneq \Delta$. Then $P_\theta \cap \widetilde W = \langle \I(\theta), \Mbar \rangle$ and $(P_\theta)_0 \cap \widetilde W = \langle \I(\theta) \rangle$.
\end{Lem}
\begin{Prf}
  Since $P_\theta$ is $B$--invariant from both sides, it is a union of Bruhat cells, so $P_\theta \cap \widetilde W = \pi^{-1}(P_\theta \cap W)$. Recall that $P_\theta = N_G(\fp_\theta)$, so $w \in W$ is in $P_\theta$ if and only if $\Ad_w \fp_\theta \subset \fp_\theta$. This holds if and only if $w$ preserves $\Sigma_0^+ \cup \Span(\theta)$. A simple computation shows that this is equivalent to $\Psi_w \subset \Span(\theta)$, which in turn is equivalent to $w \in \langle \theta \rangle \subset W$ by \autoref{lem:properties_of_Psi}. This proves the first equality.

  For the second one, note that $\I(\alpha) \in (P_{\alpha})_0 \subset (P_\theta)_0$ by \autoref{lem:doublecoset_generator} for every $\alpha \in \theta$, so $\langle \I(\theta) \rangle \subset (P_\theta)_0 \cap \widetilde W$. By \autoref{lem:PR_decomposition}, $P_{\langle \I(\theta)\rangle}$ is a closed subgroup of $G$ and $P_{\langle \I(\theta) \rangle} \subset (P_\theta)_0$. But by the preceding paragraph, $P_\theta = P_{\langle \I(\theta), \Mbar \rangle} = P_{\langle \I(\theta) \rangle} \Mbar$ is a union of finitely many copies of $P_{\langle \I(\theta) \rangle}$. This is only possible if $P_{\langle \I(\theta) \rangle} = (P_\theta)_0$.
\end{Prf}
\begin{Prf}[of \autoref{prop:parabolics}]
  \autoref{lem:PR_decomposition} shows that $P_R = B_0 R B_0$ is a closed subgroup containing $B_0$ for every oriented parabolic type $R$. On the other hand, the Lie algebra of such a subgroup $P$ contains $\fb$ and is therefore of the form $\fp_\theta$ for some $\theta \subset \Delta$ \cite[Proposition 7.76]{Knapp}. So $(P_\theta)_0 \subset P \subset P_\theta$ and, by \autoref{lem:Ptheta_cells}, $\langle \I(\theta) \rangle \subset P \cap \widetilde W \subset \langle \I(\theta) \rangle \Mbar$. Let $E = P \cap \Mbar$. Then $\Mbar_\theta \subset E \subset \Mbar$ and $P \cap \widetilde W = \langle \I(\theta) \rangle E$ is an oriented parabolic type.

  So the maps in both directions are well--defined. It is clear by \autoref{prop:bruhat} that they are inverses of each other.
\end{Prf}
\subsection{Relative positions}
Let $P_R$ and $P_S$ be the oriented parabolic subgroups of types $R = \langle \I(\theta), E \rangle$ and $S = \langle \I(\eta), F \rangle$ and let $\F_R$, $\F_S$ be the oriented flag manifolds.
\begin{Def}	\label{def:relative_position}
  The \emph{set of relative positions} is the quotient
  \[\widetilde W_{R,S} = G \backslash (\F_R \times \F_S),\]
  where $G$ acts diagonally on $\F_R \times \F_S$. The projection
  \[\pos_{R,S} \colon \F_R \times \F_S \to \widetilde W_{R,S}\]
  is called the \emph{relative position map}.\index{posRS@$\pos_{R,S}$ relative position map}
\end{Def}
\begin{Ex}
  Consider the group $G=\SL(2,\bR)$ and $R=S=\{1\}$, so that both $\F_R$ and $\F_S$ are identified with $S^1$, the space of oriented lines in $\bR^2$. Then there are two $2$--dimensional and two $1$--dimensional $G$--orbits in the space $S^1\times S^1$. The $2$--dimensional orbits consist of all transverse pairs $(v,w)$ defining a positively or negatively oriented basis of $\bR^2$. The $1$--dimensional orbits consist of pairs $(v,\pm v)$.
\end{Ex}
The relative positions admit a combinatorial description in the framework of the preceding sections. This is the main reason why we consider the parabolic types as subgroups of $\widetilde W$, and the reason for the notation $\widetilde{W}_{R,S}$. When we write $\widetilde W_{R,S}$ in the following, we will usually regard it as $R \backslash \widetilde W / S$ and we will write double brackets $\dc{\,\cdot\,}$ for equivalence classes in these quotients.\index{W_RS@$\widetilde W_{R,S} = R \backslash \widetilde W / S$ oriented relative positions}
\begin{Prop}\label{prop:position_combinatorics}
  There are bijections
  \[R \backslash \widetilde W / S \to P_R \backslash G / P_S \to G \backslash (\F_R \times \F_S).\]
  The first map is induced by the inclusion of $N_K(\fa)$ into $G$ and the second by mapping $g \in G$ to $([1], [g]) \in \F_R \times \F_S$. In particular, $\widetilde W_{R,S}$ is a finite set.
\end{Prop}
\begin{Prf}
  It is clear from the definitions that the maps are well--defined and the second map is injective. It is surjective since $G$ acts transitively on $\F_R$ so every pair in $\F_R \times \F_S$ can thus be brought into the form $([1],[g])$ by the diagonal action. To see that the first map is surjective, let $P_R g P_S \in P_R \backslash G / P_S$. By \autoref{prop:bruhat} $g \in B_0wB_0$ for some $w \in \widetilde W$. Then $\dc{w} \in R \backslash \widetilde W / S$ maps to $[g]$.

  To prove injectivity of the first map, let $w, w' \in \widetilde W$ with $P_R w P_S = P_R w' P_S$. Since $P_R = B_0RB_0$ and $P_S = B_0 S B_0$ by \autoref{lem:PR_decomposition}, we can write $w' \in P_R w P_S = B_0 R B_0 w B_0 S B_0$. By \autoref{lem:B_0R}, $B_0RB_0 w B_0 S B_0 = B_0 R w S B_0$, and by \autoref{prop:bruhat} this implies $w' \in RwS$, proving injectivity.
\end{Prf}
A particularly important subset of relative positions consists of transverse or ``maximally generic'' ones. It will be studied in more detail in \autoref{sec:involutions}.
\begin{Def}	\label{def:transverse_position}
  A element $w_0 \in \widetilde{W}$ is called \emph{transverse} if it is a lift of the longest element of $W$. The set of transverse elements in $\widetilde{W}$ will be denoted by $T$. \\
  A relative position $\dc{w_0} \in \widetilde{W}_{R,S}$ is called \emph{transverse} if it is a projection of a transverse position in $\widetilde{W}$. The set of transverse positions in $\widetilde{W}_{R,S}$ will be denoted by $T_{R,S}$.
\end{Def}\index{T@$T$, $T_{R,S}$ transverse positions}
\subsection{Refined Schubert strata} \label{sec:Schubert_strata}
\begin{Def}
  Let $f \in \F_R$ and $\dc{w} \in \widetilde{W}_{R,S}$. Then we call the set
  \[ C_{\dc{w}}^{R,S}(f) := \{ f' \in \F_S \mid \pos_{R,S}(f,f') = \dc{w} \} \]
  of flags at position $\dc{w}$ with respect to $f$ a \emph{refined Schubert stratum}. We sometimes omit the superscript $R,S$ if it is clear from the context.\index{CwRS@$C_{\dc{w}}^{R,S}(f)$ refined Schubert stratum}
\end{Def}
Every refined Schubert stratum admits the following simple (but in general not injective) parametrization by $R$ and the unipotent subgroup $N \subset G$:
\begin{Lem}	\label{lem:Bruhat_cells_projection}
  Let $w \in \widetilde W$. Then
  \[ C^{R,S}_{\dc{w}}([1]) = NR[w] \subset \F_S. \]
\end{Lem}
\begin{Prf}
  Let $g\in G$ such that $[g] \in C^{R,S}_{\dc{w}}([1])$. Then by \autoref{lem:B_0R} we have
  \[ g \in P_R w P_S = B_0 R w P_S. \]
  Using the Iwasawa decomposition $B_0 = NAZ_K(\fa)_0$ and the fact that both $A$ and $Z_K(\fa)$ are normalized by $N_K(\fa)$, this implies
  \[ B_0 R w P_S = NAZ_K(\fa)_0 R w P_S = NRwAZ_K(\fa)_0P_S = NRwP_S. \qedhere\]
\end{Prf}
If $R = S$ and $\dc{w_0}\in\widetilde{W}_{R,R}$ is a transverse position, then $C_{\dc{w_0}}(f)$ is a cell, parametrized as follows.
\begin{Lem}\label{lem:transverse_parametrization}
  Assume that $\iota(\theta) = \theta$ and let $w_0 \in T$ such that $w_0 E w_0^{-1} = E$. Define
  \[\fn_\theta^- = \!\!\!\bigoplus_{\alpha \in \Sigma^- \setminus \, \Span (\theta)} \fg_\alpha.\]
  Then the map
  \[\varphi \colon \fn_\theta^- \to C^{R,R}_{\dc{w_0}}([w_0^{-1}]), \qquad X \mapsto [e^X]\]
  is a diffeomorphism.
\end{Lem}\index{nt_@$\fn_\theta^-$ subalgebra of $\fn^-$}
\begin{Prf}
  Let $N_\theta^- \subset G$ be the connected subgroup with Lie algebra $\fn_\theta^-$. As a subgroup of $N^-$ its exponential map $\fn_\theta^- \to N_\theta^-$ is a diffeomorphism. So it suffices to show that the projection map $\widetilde\varphi \colon N_\theta^- \to \F_R$ is a diffeomorphism onto $C_{\dc{w_0}}([w_0^{-1}])$.

  First we verify that $C_{\dc{w_0}}([w_0^{-1}]) = \{[n] \mid n \in N_\theta^- \}$. We have
  \begin{align*}
    C_{\dc{w_0}}([w_0^{-1}]) &= \{f \in \F_R \mid \pos_{R,R}([1], w_0 f) = \dc{w_0} \} \\
                             &= \{f \in \F_R \mid \exists p \in G \colon [1] = [p], \  w_0 f = [pw_0]\} = \{[w_0^{-1} p w_0] \mid p \in P_R\}
  \end{align*}
  and $w_0 N_\theta^- w_0^{-1} \subset (P_{\iota(\theta)})_0 \subset P_R$, so it remains to show that $w_0^{-1}P_Rw_0 \subset N_\theta^- P_R$. As a consequence of the Langlands decomposition \cite[Propositions 7.82(a) and 7.83(d)]{Knapp} we can write $P_\theta = N_\theta^+ Z_G(\fa_\theta)$ with $N_\theta^+ = w_0 N_{\iota(\theta)}^- w_0^{-1} = w_0 N_\theta^- w_0^{-1}$ and $\fa_\theta = \bigcap_{\beta \in \theta} \ker \beta$. Now $\Ad_{w_0}$ preserves $\fa_\theta$ and thus $w_0^{-1} Z_G(\fa_\theta) w_0 = Z_G(\fa_\theta)$, hence $w_0^{-1}P_\theta w_0 = N_\theta^- Z_G(\fa_\theta)$. As $N_\theta^-$ is connected, we even get $w_0^{-1}(P_\theta)_0w_0 = N_\theta^- Z_G(\fa_\theta)_0$ and therefore
  \[w_0^{-1}P_R w_0 = w_0^{-1} (P_\theta)_0 E w_0 = w_0^{-1} (P_\theta)_0 w_0 E = N_\theta^- Z_G(\fa_\theta)_0 E \subset N_\theta^- P_R.\]
  To prove injectivity of $\widetilde\varphi$, let $n, n' \in N_\theta^-$ with $[n] = [n']$. Then $n^{-1}n' \in N_\theta^- \cap P_R = \{1\}$ by \cite[Proposition 7.83(e)]{Knapp}, so $\widetilde \varphi$ is injective. To see that $\widetilde \varphi$ is regular, we observe that $\fn_\theta^-$ is composed of the root spaces of roots in $\Sigma^- \setminus \Span (\theta)$ while $\fp_\theta$ has the root spaces $\Sigma^+ \cup \Span (\theta)$. So $\fg = \fn_\theta^- \oplus \fp_\theta$ and $D_1\widetilde \varphi \colon \fn_\theta^- \to \fg / \fp_\theta$ is an isomorphism. By equivariance we then see that $\widetilde\varphi$ is a diffeomorphism onto its image.
\end{Prf}
\subsection{The Bruhat order}	\label{sec:Bruhat_order}
Again, let $P_R$ and $P_S$ be the oriented parabolic subgroups of types $R = \langle \I(\theta), E \rangle$ and $S = \langle \I(\eta), F \rangle$ and let $\F_R$, $\F_S$ be the oriented flag manifolds.
\begin{Def}
  The \emph{Bruhat order} on $\widetilde W_{R,S} = G \backslash (\F_R \times \F_S)$ is the inclusion order on closures, i.e.
  \[G(f_1,f_2) \leq G(f_1',f_2') \quad \Leftrightarrow \quad \overline{G(f_1,f_2)} \subset \overline{G(f_1',f_2')}.\]
\end{Def}\index{<=@$\leq$ Bruhat order on $\widetilde W_{R,S}$}
In other words, if we have sequences of flags $f_n \to f \in \F_R$ and $f_n' \to f'\in \F_S$ with $\pos_{R,S}(f_n,f_n')$ constant, then
\[\pos_{R,S}(f,f') \leq \pos_{R,S}(f_n,f_n').\]
The Bruhat order thus encodes the genericity of a pair of flags.

If we view the set of relative positions as the double quotient $R \backslash \widetilde W / S$ via the correspondence in \autoref{prop:position_combinatorics}, the Bruhat order is given as follows.
\begin{Lem}
  Let $w,w' \in \widetilde{W}$ and denote by $\dc{w},\dc{w'} \in \widetilde{W}_{R,S}$ the equivalence classes they represent. Then we have
  \[\dc{w} \leq \dc{w'} \quad \Leftrightarrow \quad \overline{P_RwP_S} \subset \overline{P_Rw'P_S}.\]
\end{Lem}
\begin{Prf}
  Recall that $w$ represents the orbit $G([1],[w])$. The inequality $\dc{w} \leq \dc{w'}$ is equivalent to the existence of a sequence $(g_n) \in G^\bN$ such that $[g_n] \to [1]$ in $\F_R$ and $[g_nw'] \to [w]$ in $\F_S$. This means that there exist sequences $(p_n) \in P_R^\bN$ and $(q_n) \in P_S^\bN$ such that $g_np_n \to 1$ and $g_nw'q_n \to w$. Writing
  \[ g_np_np_n^{-1}w'q_n \to w\]
  shows that $p_n^{-1}w'q_n \to w$.\\
  Conversely, if such sequences $(p_n)$ and $(q_n)$ are given, we can simply take $g_n = p_n$ in the description above.
\end{Prf}
The following lemma shows how the Bruhat order on $\widetilde W$ relates to that on the quotients $\widetilde W_{R,S} = R \backslash \widetilde W / S$.
\begin{Lem}\label{lem:order_projection}
  Let $R \subset R'$ and $S \subset S'$ be oriented parabolic types. In this lemma, we write $\dc{w}$ for the equivalence class of $w \in \widetilde W$ in $\widetilde W_{R,S}$ and $\dc{w}'$ for its equivalence class in $\widetilde W_{R',S'}$. Then for every $w_1, w_2 \in \widetilde W$
  \begin{enumerate}
  \item If $\dc{w_1} \leq \dc{w_2}$, then $\dc{w_1}' \leq \dc{w_2}'$.
  \item If $\dc{w_1}' \leq \dc{w_2}'$, then there exists $w_3 \in \widetilde W$ with $\dc{w_3}' = \dc{w_2}'$ and $\dc{w_1} \leq \dc{w_3}$.
  \item If $\dc{w_1}' \leq \dc{w_2}'$, then there exists $w_3 \in \widetilde W$ with $\dc{w_3}' = \dc{w_1}'$ and $\dc{w_3} \leq \dc{w_2}$.
  \end{enumerate}
\end{Lem}

\begin{Prf}
  If $\dc{w_1} \leq \dc{w_2}$ then $w_1 \in \overline{P_R w_2 P_S} \subset \overline{P_{R'} w_2 P_{S'}}$. As the last term is $P_{R'}$--left and $P_{S'}$--right invariant and closed, this implies $\overline{P_{R'}w_1P_{S'}} \subset \overline{P_{R'}w_2P_{S'}}$, hence \itemnr{1}.

  The assumption in \itemnr{2} is equivalent to $w_1 \in \overline{P_{R'} w_2 P_{S'}}$. By \autoref{lem:B_0R} and \autoref{lem:PR_decomposition}, $P_{R'}w_2P_{S'} = B_0 R' w_2 S' B_0$, so there exist $r \in R'$ and $s \in S'$ such that $w_1 \in \overline{B_0 r w_2 s B_0}$. So $w_3 = rw_2s$ satisfies the properties we want.

  In part \itemnr{3}, as $\dc{w_1}' \leq \dc{w_2}'$ there is a sequence $(g_n) \in G^\bN$ such that $[g_n] \to [1]$ in $\F_{R'}$ and $[g_n w_2] \to [w_1]$ in $\F_{S'}$. Passing to a subsequence, we can also assume that $[g_n] \to f_1 \in \F_R$ and $[g_n w_2] \to f_2 \in \F_S$. Let $\dc{w_3} = \pos_{R,S}(f_1, f_2)$. Then $\dc{w_3} \leq \dc{w_2}$ and $\dc{w_3}' = \pos_{R',S'}(\pi_{R'}(f_1), \pi_{S'}(f_2)) = \pos_{R',S'}([1], [w_1]) = \dc{w_1}'$.
\end{Prf}
The Bruhat order on $\widetilde W$ is defined by orbit closures of the $B_0 \times B_0$--action. \Autoref{prop:combinatorial_bruhat_order} describes this in combinatorial terms. Combining this with \autoref{lem:order_projection} allows us to also describe the Bruhat order on $\widetilde W_{R,S}$ combinatorially. Essentially, we get everything lower than $\dc{w}$ in the Bruhat order by deleting or squaring letters in a suitable reduced word for $w$.
\begin{Prop}	\label{prop:Bruhat_order_quotient}
  For any $w \in \widetilde W$ choose $\alpha_1, \dots, \alpha_k \in \Delta$ and $m \in \Mbar$ such that $w = \I(\alpha_1) \dots \I(\alpha_k) \, m$ and that this is a reduced word, meaning $k = \ell(m)$. Then define\index{Aw@$A_w$ set of elements $\leq w$}
  \[A_w = \{\I(\alpha_1)^{i_1} \cdots \I(\alpha_k)^{i_k} m \mid i_1,\dots,i_k \in \{0,1,2\}\} \subset \widetilde W.\]
  This is independent of the choice of reduced word. The Bruhat order on $\widetilde W$ is given by
  \begin{equation}
     \label{eq:bruhat_combinatorial}
    w' \leq w \quad \Leftrightarrow \quad w' \in A_w.
  \end{equation}
  On $\widetilde W_{R,S}$ it is given by
\begin{equation}
    \label{eq:bruhat_combinatorial_quotient}
    \dc{w'} \leq \dc{w} \quad \Leftrightarrow \quad w' \in R A_w S \quad \Leftrightarrow \quad w' \in \!\! \bigcup_{r\in R,\;s \in S} \!\! A_{rws}.
  \end{equation}
  The Bruhat order on $\widetilde{W}_{R,S}$ is a partial order.
  \end{Prop}
\begin{Prf}
  The well--definedness of $A_w$ and \eqref{eq:bruhat_combinatorial} hold by \autoref{prop:combinatorial_bruhat_order}. For the general case \eqref{eq:bruhat_combinatorial_quotient}, if $w' \in A_{rws}$ or $w' \in r^{-1} A_w s^{-1}$ for some $r \in R$, $s \in S$, then $w' \leq r w s$ or $rw's \leq w$ in $\widetilde W$, respectively. By \autoref{lem:order_projection}\itemnr{1} both inequalities imply $\dc{w'} \leq \dc{w}$. Conversely, if $\dc{w'} \leq \dc{w}$, then $w' \leq rws$ and $r'w's' \leq w$ for some $r,r' \in R$ and $s,s' \in S$ by \autoref{lem:order_projection}\itemnr{2} and \autoref{lem:order_projection}\itemnr{3}, so $w' \in A_{rws}$ and $w' \in r^{-1}A_w s^{-1}$.

  The Bruhat order is transitive and reflexive by its definition. To show antisymmetry, suppose $\dc{w} \leq \dc{w'} \leq \dc{w}$ in $\widetilde W_{R,S}$. Then $RA_wS = RA_{w'}S$. Let $\mathcal L_{w} \subset RA_wS$ be the subset of elements which are maximal in $R A_w S = \bigcup_{r,s}A_{rws}$ with respect to $\ell$. Then every element of $\mathcal L_w$ is also maximal in $A_{rws}$ for some $r \in R$, $s \in S$. But the unique longest element of $A_{rws}$ is $rws$, since squaring and deleting letters only reduces $\ell$. So $\mathcal L_w \subset RwS$ and, since $\mathcal L_w = \mathcal L_{w'} \neq \varnothing$, $RwS = Rw'S$.
\end{Prf}
The following characterization of the Bruhat order on $\widetilde W$ will also be useful later:
\begin{Lem}\label{lem:folding_order}
  Let $w, w' \in \widetilde W$ with $\ell(w') = \ell(w) + 1$. Let
  \[Q = \{w \, \I(\alpha)^{\pm 1}\, w^{-1} \mid w \in \widetilde W, \alpha \in \Delta\} \subset \widetilde W\]
  be the set of conjugates of the standard generators or their inverses. Then
  \[w \leq w' \quad \Longleftrightarrow \quad \exists q \in Q \colon w = q w'. \]
\end{Lem}
\begin{Prf}
  The implication `$\Rightarrow$' follows from \autoref{prop:combinatorial_bruhat_order} by choosing $q$ of the form
  \[q = \I(\alpha_1) \dots \I(\alpha_{i-1}) \I(\alpha_i)^{\pm 1} \I(\alpha_{i-1})^{-1} \dots \I(\alpha_1)^{-1}\]
  for some $i$. For the other direction, assume that $w = qw'$ and write
  \[w' = \I(\alpha_1) \dots \I(\alpha_k) \, m\]
  for some $\alpha_1, \dots, \alpha_k \in \Delta$ with $k = \ell(w')$ and $m \in \Mbar$. Then $\pi(w') = \alpha_1 \dots \alpha_k$ and by the strong exchange property of Coxeter groups \cite[Theorem 1.4.3]{BjoernerBrenti}
  \[ \pi(w) = \pi(q) \pi(w') = \alpha_1 \dots \widehat{\alpha_i} \dots \alpha_k \]
  for some $i$, so $\pi(q) = (\alpha_1 \dots \alpha_{i-1}) \alpha_i (\alpha_1 \dots \alpha_{i-1})^{-1}$. Set $c = \I(\alpha_1) \dots \I(\alpha_{i-1}) \in \widetilde W$. Then $c^{-1} q c \in \pi^{-1}(\alpha_i) \cap Q = \{\I(\alpha_i)^{\pm 1}\}$ by \autoref{lem:Enormal}. So
  \[w = qw' = c\,\I(\alpha_i)^{\pm 1} c^{-1} w' = \I(\alpha_1) \dots \I(\alpha_{i-1}) \I(\alpha_i)^{1\pm 1} \I(\alpha_{i+1}) \dots \I(\alpha_k) \, m \leq w', \]
  where the inequality at the end follows by \autoref{prop:combinatorial_bruhat_order}.
\end{Prf}
The following lemmas will be useful when calculating with relative positions.
\begin{Lem}	\label{lem:position_rightmult}\ \reallynopagebreak
  \begin{enumerate}
  \item $Z_K(\fa)$ normalizes the subgroups $P_R, P_S$. Consequently, the (finite abelian) group $\Mbar/E$ acts on $\F_R$ by right multiplication, and this action is simply transitive on fibers.
  The analogous statement holds for $\Mbar/F$ acting on $\F_S$.\\
  Furthermore, $\Mbar/E$ acts on $\widetilde{W}_{R,S}$ by left multiplication and $\Mbar/F$ acts on $\widetilde{W}_{R,S}$ by right multiplication. Both of these actions preserve the Bruhat order.
  \item For any $f_1\in \F_R$, $f_2 \in \F_S$, $m_1 \in \Mbar/E$, $m_2 \in \Mbar/F$, right multiplication by $m_1$ and $m_2$ has the following effect on relative positions:
    \[ \pos_{R,S}(R_{m_1}(f_1),R_{m_2}(f_2)) = m_1^{-1} \pos_{R,S}(f_1,f_2) m_2 \]
  \end{enumerate}
\end{Lem}
\begin{Prf}\ \reallynopagebreak
  \begin{enumerate}
  \item It follows e.g. from \autoref{lem:Enormal} that $\Mbar$ normalizes $R = \langle \I(\theta), E \rangle$ and $S = \langle \I(\eta), F \rangle$. Furthermore, $Z_K(\fa)$ normalizes $B_0$ and thus also $P_R = B_0 R B_0$ and $P_S = B_0 S B_0$. This implies that the actions of $\Mbar$ on $\F_R$ and $\F_S$ by right multiplication and on $\widetilde{W}_{R,S}$ by left and right multiplication are well--defined. Since $E$ resp. $F$ acts trivially, we obtain the induced actions of $\Mbar/E$ resp. $\Mbar/F$. The action of $\Mbar/E$ on $\F_R$ is simply transitive on each fiber (over $\F_\theta$) since $\Mbar \cap R = E$ (see \autoref{rem:types}\itemnr{2}); in the same way, $\Mbar/F$ acts simply transitively on fibers of $\F_S$.

  The actions on $\widetilde W_{R,S}$ preserve the Bruhat order since, for $m,m' \in \Mbar$,
    \begin{align*}
      \dc{w} \leq \dc{w'} & \Leftrightarrow \overline{P_RwP_S} \subset \overline{P_Rw'P_S} \Leftrightarrow m\overline{P_RwP_S}m' \subset m\overline{P_Rw'P_S}m' \\
                          & \Leftrightarrow \overline{P_Rmwm'P_S} \subset \overline{P_Rmw'm'P_S} \Leftrightarrow \dc{mwm'} \leq \dc{mw'm'}.
    \end{align*}
  \item Let $\pos_{R,S}(f_1,f_2) = \dc{w} \in \widetilde{W}_{R,S}$. This means that there exists some $g \in G$ such that $g(f_1,f_2) = ([1],[w])$. It follows that
    \[ m_1^{-1}g(R_{m_1}(f_1),R_{m_2}(f_2)) = m_1^{-1}([m_1],[wm_2]) = ([1],[m_1^{-1}wm_2]). \]
    So we obtain $\pos_{R,S}(R_{m_1}(f_1),R_{m_2}(f_2)) = \dc{m_1^{-1}wm_2}$.
    \qedhere
  \end{enumerate}
\end{Prf}
\begin{Cor}	\label{cor:position_rightmult}
  Let $f\in \F_R$, $w \in \widetilde{W}$ and $m \in \Mbar$. Then we have
  \[ R_m(C_{\dc{w}}^{R,S}(f)) = C_{\dc{wm}}^{R,S}(f) = C_{\dc{w}}^{R,S}(R_{wmw^{-1}}(f)). \]
\end{Cor}
\begin{Prf}
  From the previous lemma we obtain
  \[\pos_{R,S}(f,R_{m^{-1}}(f')) = \dc{w} \Leftrightarrow \pos_{R,S}(f,f') = \dc{wm} \Leftrightarrow \pos_{R,S}(R_{wmw^{-1}}(f), f') = \dc{w}. \qedhere\]
\end{Prf}
We close this section with an inequality for relative positions that will play an important role in the proof of cocompactness of domains of discontinuity in \autoref{sec:cocompactness}. It can be read as a triangle inequality if the position $w_0$ is a transverse one.
\begin{Lem}	\label{lem:triangle_ineq}
  Let $w_0,w_1,w_2 \in \widetilde{W}$ with $w_0 R w_0^{-1} = R$. Assume there are $f_1,f_2 \in \F_R$ and $f_3 \in \F_S$ such that
  \begin{align*}
    \pos_{R,R}(f_1,f_2) & = \dc{w_0} \in \widetilde{W}_{R,R}, \\
    \quad \pos_{R,S}(f_1,f_3) & = \dc{w_1} \in \widetilde{W}_{R,S}, \\
    \quad \pos_{R,S}(f_2,f_3) & = \dc{w_2} \in \widetilde{W}_{R,S}.
  \end{align*}
  Then
  \[ \dc{w_1} \geq \dc{w_0w_2} \]
  in $\widetilde{W}_{R,S}$.
\end{Lem}
\begin{Prf}
  Using the $G$--action on pairs, we can assume that $(f_1,f_2) = ([1],[w_0])$. Then, since $\pos_{R,S}(f_2,f_3) = \pos_{R,S}([w_0],f_3) = \dc{w_2}$, \autoref{lem:Bruhat_cells_projection} implies that $f_3$ has a representative in $G$ of the form $w_0urw_2$ for some $u \in N$ and $r\in R$. We want to find elements $g_n\in G$ such that
  \[ g_n(f_1,f_3) = g_n([1],[w_0urw_2]) \xrightarrow{n\to\infty} ([1],[w_0w_2]). \]
  Let $(A_n) \in \aplusbar^\bN$ be a sequence with $\alpha(A_n) \to \infty$ for all $\alpha \in \Delta$ and $g_n = w_0r^{-1}e^{-A_n}w_0^{-1}$. Then $g_n \in P_R$ since $A$ and $R$ are normalized by $w_0$. Observe that $w_2^{-1}r^{-1}e^{A_n}rw_2 \in A \subset P_S$, since $A$ is normalized by all of $\widetilde{W}$. Then $g_n$ stabilizes $[1] \in \F_R$, and we calculate
  \begin{align*}
    g_n[w_0urw_2] &= [(w_0r^{-1}e^{-A_n}w_0^{-1})w_0urw_2(w_2^{-1}r^{-1}e^{A_n}rw_2)] \\
                  &= [w_0r^{-1}e^{-A_n}ue^{A_n}rw_2] \xrightarrow{n\to\infty} [w_0w_2],
  \end{align*}
  where we used that $e^{-A_n}u e^{A_n} \xrightarrow{n\to\infty} 1$.
\end{Prf}
\subsection{Transverse positions}	\label{sec:involutions}
In this section, we consider the action of transverse elements $w_0\in T\subset\widetilde{W}$ (\autoref{def:transverse_position}) on the set of relative positions $\widetilde{W}_{R,S}$. We want to identify those $w_0$ acting as order--reversing involutions. There are in general multiple choices which depend on the oriented parabolic type $R$. Now we also assume that $\iota(\theta) = \theta$, where $R = \langle \I(\theta), E \rangle$.

In contrast to the unoriented setting, there can be multiple transverse positions. Indeed, there are always $|\Mbar|$ transverse positions in $\widetilde{W}$. However, the situation is less obvious in double quotients $\widetilde{W}_{R,S}$: In the example in \autoref{sec:ideal_spheres}, there is a unique transverse position although the projection to unoriented relative positions is nontrivial. We first observe that every transverse position is maximal in the Bruhat order, thus any two of them are incomparable.

\begin{Lem}	\label{lem:transverse_maximal}
  Let $\dc{w} \in T_{R,S}$ be a transverse position. Then $\dc{w}$ is maximal in the Bruhat order.
\end{Lem}

\begin{Prf}
  Since $\dc{w}$ is transverse, we may assume that $\ell(w)$ is maximal. Let $\dc{w'} \in \widetilde{W}_{R,S}$ be such that $\dc{w}\leq\dc{w'}$. By \autoref{prop:Bruhat_order_quotient}, this implies the following. There exist $\tilde r\in R, s\in S$, and we can write $\tilde rw's = \I(\alpha_1)\ldots \I(\alpha_k)m$, where $\I(\alpha_i), \ \alpha_i \in \Delta$ are our preferred choice of generators for $\widetilde W$, $k = \ell(\tilde rw's)$, and $m\in \Mbar$. Furthermore, a word representing $w$ is obtained from $\tilde rw's = \I(\alpha_1)\ldots \I(\alpha_k)m$ by squaring or deleting some of the $\I(\alpha_i)$. However, if any letters were indeed squared or deleted, $\ell(w)$ would by strictly smaller and thus not maximal, so the two words must be equal.
\end{Prf}

The next lemma gives a criterion for an element of $T$ to act as an order--reversing involution on relative positions.
\begin{Lem}\label{lem:involution_order_reversing}
  Let $w_0 \in T\subset\widetilde W$. Then for any $w, w' \in \widetilde W$, $w \leq w'$ implies $w_0 w' \leq w_0 w$. If $w_0 E w_0^{-1} = E$ and $w_0^2 \in E$, then $w_0$ acts as an order--reversing involution on $\widetilde W_{R,S}$.
\end{Lem}
\begin{Prf}
  Assume $w \leq w'$. Then by \autoref{prop:combinatorial_bruhat_order} there exists a sequence
  \[ w = w_1 \leq \dots \leq w_k = w' \]
  with $\ell(w_{i+1}) = \ell(w_i) + 1$. So by \autoref{lem:folding_order} $w_i = q_i w_{i+1}$ for some $q_i \in Q$. Therefore,
  \[ w_0 w_{i+1} = w_0 q_i^{-1} w_i = q_i' w_0 w_i \]
  for $q_i' = w_0 q_i^{-1} w_0^{-1} \in Q$. By \cite[Corollary 2.3.3]{BjoernerBrenti}, $\ell(w_0w) = \ell(w_0)-\ell(w)$ for any $w\in\widetilde W$, so
  \[ \ell(w_0 w_i) = \ell(w_0) - \ell(w_i) = \ell(w_0) - \ell(w_{i+1}) + 1 = \ell(w_0 w_{i+1}) + 1, \]
  so $w_0w' = w_0w_k \leq \dots \leq w_0 w_1 = w_0w$ by the same lemma.

  We now show that $w_0$ normalizes $R$ if it normalizes $E$, and thus the action of $w_0$ on $\widetilde{W}_{R,S} = R \backslash \widetilde{W} / S$ by left--multiplication is well--defined. The induced action of $w_0$ on reduced roots is given by $\iota$ and $\theta$ is $\iota$--invariant. Moreover, \autoref{rem:choice_of_r}\itemnr{1} implies that for every $\alpha \in \theta$, we have $w_0 \I(\alpha) w_0^{-1} = \I(\iota(\alpha))$ or $w_0 \I(\alpha) w_0^{-1} = \I(\iota(\alpha))^{-1}$. Therefore, $\langle \I(\theta) \rangle$ is normalized by $w_0$. Since $R = \langle \I(\theta),E \rangle$, the same is true for $R$.

  If in addition we have $w_0^2\in E$, the induced action on $\widetilde{W}_{R,S}$ is an involution. It is an easy consequence of \autoref{lem:order_projection} that the action on this quotient still reverses the order.
\end{Prf}
\begin{Rem}
  For $w_0 \in T$ the condition $w_0 R w_0^{-1}$ is equivalent to $\iota(\theta) = \theta$ and $w_0 E w_0^{-1}$. Moreover, $w_0^2 \in E$ is equivalent to $w_0^2 \in R$.
\end{Rem}
\begin{Ex}
  Consider $G=\SL(3,\bR)$ with its maximal compact $K=\SO(3,\bR)$ and $\fa = \left\{ \begin{pmatrix} \lambda_1 \\ & \lambda_2 \\ & & -\lambda_1 - \lambda_2 \end{pmatrix} \mid \lambda_1,\lambda_2 \in \bR \right\}$. The extended Weyl group $\widetilde{W} = N_K(\fa)/Z_K(\fa)_0$ consists of all permutation matrices $A$ with determinant 1 -- i.e. all matrices with exactly one $\pm 1$ entry per line and row and all other entries $0$, such that $\det(A)=1$. The transverse positions are
  \[ \begin{pmatrix} & & 1 \\ & -1 \\ 1 \end{pmatrix}, \begin{pmatrix} & & -1 \\ & -1 \\ -1 \end{pmatrix}, \begin{pmatrix} & & -1 \\ & 1 \\ 1 \end{pmatrix}, \begin{pmatrix} & & 1 \\ & 1 \\ -1 \end{pmatrix}. \]
  The first two of these are actually involutions in $\widetilde{W}$, so the condition $w_0^2 \in E$ is empty. The condition $w_0 E w_0^{-1} = E$ is a symmetry condition on $E$, similar to the condition $\iota(\theta)=\theta$. It does not depend on the choice of lift: Any other $w_0'$ is of the form $w_0' = w_0m$ for some $m\in \Mbar$, and $\Mbar$ is abelian. In \autoref{prop:R-Anosov_implies_Rhat-Anosov}, we will show that we can always assume it to hold in our setting.\\
  The last two are not involutions in $\widetilde{W}$. The smallest possible choice of $E$ containing their square is
  \[ E = \left\{ \begin{pmatrix} 1 \\ & 1 \\ & & 1 \end{pmatrix},\begin{pmatrix} -1 \\ & 1 \\ & & -1 \end{pmatrix} \right\}. \]
\end{Ex}
\subsection{Ideals}
\begin{Def}	\label{def:ideal}
  Let $(X,\leq)$ be a set equipped with a partial order. Then a subset $I\subset X$ is called an \emph{ideal} if for every $x\in I$ and $y \in X$ with $y\leq x$, we have $y\in I$.
\end{Def}
In the case of the Bruhat order on the extended Weyl group, an ideal corresponds to a $G$--invariant closed subset of $\F_R \times \F_S$. Intuitively, if a specific relative position is contained in the ideal, then all less generic relative positions are contained in it as well. Note that ``less generic'' is somewhat more subtle in the oriented case (a discussion of examples is found in \autoref{sec:examples_of_balanced_ideals}).

The existence of an involution $w_0$ on $\widetilde{W}_{R,S}$ allows us to define the following properties of an ideal $I\subset\widetilde{W}_{R,S}$. They play a crucial role in the description of properly discontinuous and cocompact group actions of oriented flag manifolds in \autoref{sec:proper_discontinuity} and \autoref{sec:cocompactness}.

\begin{Def}
  Let $I\subset\widetilde{W}_{R,S}$ be an ideal and $w_0 \in T\subset\widetilde{W}$ satisfying $w_0Ew_0^{-1} = E$ and $w_0^2 \in E$. Then
  \begin{itemize}
  \item $I$ is called $w_0$--\emph{fat} if $x\not\in I$ implies $w_0 x \in I$.
  \item $I$ is called $w_0$--\emph{slim} if $x \in I$ implies $w_0 x \not\in I$.
  \item $I$ is called $w_0$--\emph{balanced} if it is fat and slim.
  \end{itemize}
\end{Def}

Observe that there can be no $w_0$--balanced ideal if $w_0$ has a fixed point. Conversely, if $w_0$ has no fixed points, there will be $w_0$--balanced ideals by the following lemma. For the case of the Weyl group, this is proved in \cite[Proposition 3.29]{KapovichLeebPortiFlagManifolds}.

\begin{Lem}\label{lem:ideal_existence}
  Let $X$ be a partially ordered set and $\sigma \colon X \to X$ an order--reversing involution without fixed points. Then every minimal $\sigma$--fat ideal and every maximal $\sigma$--slim ideal is $\sigma$--balanced.
\end{Lem}

\begin{Prf}
  The two statements are equivalent by replacing an ideal $I$ by $X \setminus \sigma(I)$. So assume that $I \subset X$ is a minimal $\sigma$--fat ideal which is not $\sigma$--balanced. Choose a maximal element $x \in I \cap \sigma(I) \neq \varnothing$ and let $I' = I \setminus \{x\}$. If $I'$ is an ideal, it is clearly $\sigma$--fat, contradicting minimality of $I$. So $I'$ is not an ideal. Then there exist $x_1 \leq x_2$ with $x_2 \in I'$ but $x_1 \not\in I'$. So $x_1 = x$ since $I$ is an ideal. Furthermore $\sigma(x_2) \leq \sigma(x_1) = \sigma(x)$ and $\sigma(x) \in I$, so $\sigma(x_2) \in I$ and therefore $x_2 \in I \cap \sigma(I)$. Since $x$ is maximal in $I \cap \sigma(I)$ and $x\leq x_2$, this implies $x_2 = x = x_1$, a contradiction.
\end{Prf}

Using ideals in $\widetilde{W}_{R,S}$ and refined Schubert strata, we can define the following map $\Q$ assigning to each flag in $\F_R$ a subset of $\F_S$. It is the centerpiece of our construction of domains of discontinuity in \autoref{sec:domains_of_discontinuity}.

\begin{Lem}	\label{ex:equivariant_map_K}
  Let $\F_R$ and $\F_S$ be two oriented flag manifolds, and let $I\subset\widetilde{W}_{R,S}$ be an ideal. Then the map
  \begin{align*}
    \Q \colon \F_R & \to \mathcal P(\F_S) \\
    f & \mapsto \bigcup\limits_{\dc{w}\in I} C_{\dc{w}}(f)
  \end{align*}\index{Q@$\Q \colon \F_R \to \C(\F_S)$ equivariant map}%
  is $G$--equivariant with image in $\C(\F_S)$, the set of closed subsets of $\F_S$.
\end{Lem}

\begin{Prf}
  Observe that for any element $g\in G$ satisfying $[g]=f \in \F_R$ and any relative position $\dc{w}\in\widetilde{W}_{R,S}$, we have $C_{\dc{w}}(f) = gC_{\dc{w}}([1])$; in other words, the map $f \mapsto C_{\dc{w}}(f)$ from $\F_R$ to subsets of $\F_S$ is equivariant. By definition of the Bruhat order on $\widetilde{W}_{R,S}$, the closure of $C_{\dc{w}}(f)$ is given by
  \[ \overline{C_{\dc{w}}(f)} = \bigcup\limits_{\dc{w'} \leq \dc{w}} C_{\dc{w'}}(f). \]
  In particular, if $I\subset \widetilde{W}_{R,S}$ is an ideal, then $\bigcup_{\dc{w}\in I} C_{\dc{w}}(f)$ is closed for any $f\in \F_R$.
\end{Prf}
\section{\texorpdfstring{$P_R$}{PR}--Anosov representations} \label{sec:orientable_Anosov}
Let $P_R$ be the oriented parabolic of type $R = \langle \I(\theta), E \rangle$, with $\iota(\theta) = \theta$. Moreover, let $w_0 \in T\subset\widetilde W$ be a transverse position. Let $\Gamma$ be a finitely generated group and $\rho \colon \Gamma \to G$ a representation. We denote the Cartan projection by\index{m_@$\mu \colon G \to \aplusbar$ Cartan projection}
\[\mu \colon G \to \aplusbar.\]
It maps $g \in G$ to the unique element $\mu(g) \in \aplusbar$ with $g \in K \exp(\mu(g)) K$. We need the following notions in order to define Anosov representations.
\begin{Def} Let $\Gamma$ be a word hyperbolic group, $\bdry$ its Gromov boundary and $\xi \colon \bdry \to \F_\theta$ a map.
  \begin{enumerate}
  \item A sequence $(g_n) \in G^\bN$ is called \emph{$P_\theta$--divergent} if
    \[\alpha(\mu(g_n)) \to \infty \qquad \forall \alpha \in \ctheta.\]
    The representation $\rho$ is \emph{$P_\theta$--divergent} if for every divergent sequence $\gamma_n \to \infty$ in $\Gamma$ its image $\rho(\gamma_n)$ is $P_\theta$--divergent.
  \item Two elements $x,y \in \F_\theta$ are called \emph{transverse} if the pair $(x,y) \in \F_\theta \times \F_\theta$ lies in the unique open $G$--orbit. Equivalently, their relative position is represented by the longest element of $W$.

    The map $\xi$ is called \emph{transverse} if for every pair $x\neq y \in \bdry$, the images $\xi(x),\xi(y) \in \F_\theta$ are transverse.
  \item $\xi$ is called \emph{dynamics--preserving} if, for every element $\gamma \in \Gamma$ of infinite order, its unique attracting fixed point $\gamma^+ \in \bdry$ is mapped to an attracting fixed point of $\rho(\gamma)$.
  \end{enumerate}
\end{Def}
\begin{Def}[{\cite[Theorem 1.3]{GueritaudGuichardKasselWienhard}}]	\label{def:Anosov}
  The representation $\rho \colon \Gamma \to G$ is \emph{$P_\theta$--Anosov} if $\Gamma$ is word hyperbolic, $\rho$ is $P_\theta$--divergent and there is a continuous, transverse, dynamics--preserving, $\rho$--equivariant map $\xi \colon \bdry \to \F_\theta$ called \emph{limit map} or \emph{boundary map}.\index{x_@$\xi,\widehat\xi$ limit maps}\index{G_@$\Gamma$ hyperbolic group}\index{dG_@$\bdry$ Gromov boundary of $\Gamma$}
\end{Def}
We extend the notion of an Anosov representation to oriented parabolic subgroups by requiring that the limit map lifts to the corresponding oriented flag manifold.
\begin{Def}	\label{def:R-orientable_Anosov}
  Assume that $\Gamma$ is non--elementary. The representation $\rho \colon \Gamma \to G$ is \emph{$P_R$--Anosov} if it is $P_\theta$--Anosov with limit map $\xi \colon \bdry \to \F_\theta$ and there is a continuous, $\rho$--equivariant lift $\widehat \xi \colon \bdry \to \F_R$ of $\xi$. Such a map $\widehat \xi$ will be called a limit map or boundary map of $\rho$ as a $P_R$--Anosov representation. The relative position $\pos_{R,R}(\widehat\xi(x),\widehat\xi(y))$ for $x \neq y \in \bdry$ is its \emph{transversality type}.
\end{Def}
We should verify that the transversality type is in fact well--defined.
\begin{Lem}	\label{lem:transversality_type}
  The relative position $\pos_{R,R}(\widehat\xi(x),\widehat\xi(y))$ in the above definition does not depend on the choice of $x$ and $y$.
\end{Lem}
\begin{Prf}
  By \cite[8.2.I]{Gromov}, there exists a dense orbit in $\bdry \times \bdry$. By equivariance of $\widehat{\xi}$, the relative position $\dc{w_0}$ of pairs in this orbit is constant. It is a transverse position because this orbit contains (only) pairs of distinct points. An arbitrary pair $(x,y)$ of distinct points in $\bdry$ can be approximated by pairs in the dense orbit, so by continuity of $\widehat{\xi}$, we have
  \[ \pos_{R,R}(\widehat\xi(x),\widehat\xi(y)) \leq \dc{w_0}. \]
  But $\pos_{R,R}(\widehat\xi(x),\widehat\xi(y))$ is a transverse position, thus equality holds by \autoref{lem:transverse_maximal}.
\end{Prf}
\begin{Rems} \label{rem:boundary_maps}\ \reallynopagebreak
  \begin{enumerate}
  \item The definition (apart from the transversality type) makes sense for elementary hyperbolic groups, but it is not a very interesting notion in this case: The boundary has at most two points. Consequently, after restricting to a finite index subgroup, the boundary map always lifts to the maximally oriented setting. Moreover, after restricting to the subgroup preserving the boundary pointwise, the lifted boundary map holds no additional information.
  \item In the oriented setting, the boundary map $\widehat{\xi}:\bdry\to\F_R$ is not unique: For any element $[m] \in \Mbar / E$, the map $R_m \circ \widehat{\xi}$ is also continuous and equivariant. This gives all possible boundary maps in $\F_R$:\\
    Since the unoriented boundary map $\xi:\bdry \to \F_\theta$ is unique \cite[Lemma 3.3]{GuichardWienhardDomains}, an oriented boundary map $\widehat{\xi}'$ must be a lift of it. But if $\widehat\xi'$ agrees with $R_m \circ \widehat{\xi}$ at a single point, it must agree everywhere by equivariance and continuity since any orbit is dense in $\bdry$ \cite[Proposition 4.2]{BenakliKapovich}. If the transversality type of $\widehat \xi$ was $\dc{w_0}$, then that of $R_m \circ \widehat{\xi}$ is $\dc{m^{-1}w_0m}$ by \autoref{lem:position_rightmult}.
  \end{enumerate}
\end{Rems}
The oriented flag manifold $\F_R$ in \autoref{def:R-orientable_Anosov} which is the target of the lift $\widehat{\xi}$ is not unique. However, there is a unique maximal choice of such a $\F_R$ (or equivalently, minimal choice of $R$), similar to the fact that an Anosov representation admits a unique minimal choice of $\theta$ such that it is $P_\theta$--Anosov.
\begin{Prop}	\label{prop:minimal_type}
  Let $\rho \colon \Gamma \to G$ be $P_\theta$--Anosov. Then there is a unique minimal choice of $E$ such that $\Mbar_\theta \subset E \subset \Mbar$ and $\rho$ is $P_R$--Anosov, where $R = \langle \I(\theta), E \rangle$.
\end{Prop}
\begin{Prf}
  Assume that there are two different choices $E_1$ and $E_2$ such that $\rho$ is both $P_{R_1}$-- and $P_{R_2}$--Anosov, where $R_i = \langle \I(\theta), E_i \rangle$. Let $E_3 = E_1 \cap E_2$. We will show that $\rho$ is also $P_{R_3}$--Anosov. To do so, we have to construct a boundary map into $\F_{R_3}$ from the boundary maps into $\F_{R_1}$ and $\F_{R_2}$.

  Let $\xi \colon \bdry \to \F_\theta$ be the boundary map of $\rho$ as a $P_\theta$--Anosov representation, and let $\xi_1 \colon \bdry \to \F_{R_1}$, $\xi_2 \colon \bdry \to \F_{R_2}$ be the two lifts we are given. Fix a point $x\in\bdry$, and let $F_x \in \F_{R_3}$ be a lift of $\xi(x) \in \F_\theta$. Denote by $\pi_1 \colon \F_{R_3} \to \F_{R_1}$ and $\pi_2 \colon \F_{R_3} \to \F_{R_2}$ the two projections. After right--multiplying $\xi_1$ with an element of $\Mbar/E_1$ and $\xi_2$ with an element of $\Mbar/E_2$, we may assume that $\pi_i(F_x) = \xi_i(x)$ for $i=1,2$. Set $\xi_3(x) := F_x$ and observe the following general property:

  For every point $y\in\bdry$, there is at most one flag $F_y \in \F_{R_3}$ satisfying $\pi_i(F_y) = \xi_i(y)$ for $i=1,2$. Indeed, if $gP_{R_3}$ and $hP_{R_3}$ satisfy $gP_{R_i} = hP_{R_i}$ for $i=1,2$, there are elements $p_i \in P_{R_i}$ such that $g = hp_1 = hp_2$. This implies that $h^{-1}g \in P_{R_1} \cap P_{R_2} = P_{R_3}$.

  By equivariance of $\xi_1,\xi_2$ and uniqueness of lifts to $\F_{R_3}$, we can extend $\xi_3$ equivariantly to a map $\xi_3 \colon \Gamma x \to \F_{R_3}$. It is a lift of both $\xi_1|_{\Gamma x}$ and $\xi_2|_{\Gamma x}$. Recall that the orbit $\Gamma x$ is dense in $\bdry$ for any choice of $x$ (\cite[Proposition 4.2]{BenakliKapovich}). Using the corresponding properties of $\xi_1$ and $\xi_2$, we now show that this map is continuous and extends continuously to all of $\bdry$.  Let $x_n \in \Gamma x$ and assume that $x_n \to x_\infty \in \bdry$. Then $\xi_i(x_n) \to \xi_i(x_\infty)$ for $i=1,2$. Therefore, there exist $m_n \in \Mbar/ (E_1\cap E_2)$ such that $\xi_3(x_n)m_n$ converges in $\F_{R_3}$. By injectivity of the map
  \[ \Mbar/E_1 \cap E_2 \to \Mbar/E_1 \times \Mbar/E_2 \]
  and convergence of $\xi_i(x_n), \ i=1,2$, $m_n$ must eventually be constant. Thus the limit $\xi_3(x_\infty) := \lim_{n\to\infty} \xi_3(x_n)$ exists and is the unique lift of $\xi_1(x_\infty),\xi_2(x_\infty)$ to $\F_{R_3}$.
\end{Prf}
The following proposition shows that given a $P_R$--Anosov representation $\rho$ of transversality type $\dc{w_0}$, we may always assume that $R$ is stable under conjugation by $w_0$. This appeared as an assumption in \autoref{sec:involutions} and plays a role later on when showing that balanced ideals give rise to cocompact domains of discontinuity.
\begin{Prop}	\label{prop:R-Anosov_implies_Rhat-Anosov}
  Let $w_0 \in T$ and $E' = E \cap w_0 E w_0^{-1}$, and let $R = \langle \I(\theta), E \rangle$ and $R' = \langle \I(\theta), E' \rangle$. Assume that $\rho \colon \Gamma \to G$ is $P_R$--Anosov with a limit map $\widehat \xi \colon \bdry \to \F_R$ of transversality type $\dc{w_0} \in \widetilde W_{R,R}$. Then $\rho$ is $P_{R'}$--Anosov.
\end{Prop}
\begin{Prf}
  Let $x \neq z \in \bdry$, and consider the images $\widehat{\xi}(x),\widehat{\xi}(z) \in \F_R$. We claim that there is a unique lift $\eta_x(z) \in \F_{R'}$ satisfying
  \begin{itemize}
  \item $\eta_x(z)$ projects to $\widehat{\xi}(z)$.
  \item $\pos_{R,R'}(\widehat{\xi}(x),\eta_x(z)) = \dc{w_0}$.
  \end{itemize}
  To show this, let us first fix a good representative in $G$ for $\widehat{\xi}(z)$: Since $\pos_{R,R}(\widehat{\xi}(x),\widehat{\xi}(z)) = \dc{w_0}$, there exists $h\in G$ such that
  \[ h(\widehat{\xi}(x),\widehat{\xi}(z)) = ([1],[w_0]). \]
  Then $h^{-1}w_0 =: g \in G$ represents the flag $\widehat{\xi}(z)$ and also satisfies
  \[ \pos_{R,R'}(\widehat{\xi}(x),[g]) = \dc{w_0}. \]
  Any lift of $\widehat{\xi}(z)$ into $\F_{R'}$ can be written as $[gm] \in \F_{R'}$ for some $m\in E$. By \autoref{lem:position_rightmult}, we have $\pos_{R,R'}(\widehat{\xi}(x),[gm]) = \dc{w_0m}$. We claim that $\dc{w_0m} = \dc{w_0} \in \widetilde{W}_{R,R'}$ implies that $m \in E'$ and therefore $[gm] = [g] \in \F_{R'}$, proving uniqueness of $\eta_x(z)$. Indeed, if $w_0m = rw_0 r'$ for some $r\in R, r'\in R'$, we obtain
  \[ m = w_0^{-1}rw_0 r' \in w_0^{-1} R w_0 \cdot R' \subset w_0^{-1} R w_0. \]
  Since
  \[ E \cap w_0^{-1} R w_0 = E \cap w_0^{-1} E w_0 = E', \]
  it follows that $m \in E'$ and the lift $[g] \in \F_{R'}$ is unique.\\
  This defines a map
  \[ \eta_x \colon \bdry \setminus \{x\} \to \F_{R'} \]
  which is continuous since $\widehat \xi$ is continuous. We will show that it is independent of the choice of $x$, i.e. if $y\neq z$ is another point, we have $\eta_x(z) = \eta_y(z)$. Let $\gamma \in \Gamma$ be an element of infinite order with fixed points $\gamma^\pm \in \bdry$ such that $x\neq \gamma^-$ and $y\neq \gamma^-$. Then we have
  \[ \dc{w_0} = \pos_{R,R'}(\widehat{\xi}(x),\eta_x(\gamma^-)) = \pos_{R,R'}(\rho(\gamma)^n\widehat{\xi}(x),\rho(\gamma)^n\eta_x(\gamma^-)) \]
  for every $n\in\bN$. Moreover, $\rho(\gamma)^n\widehat{\xi}(x) \to \widehat{\xi}(\gamma^+)$ and $\rho(\gamma)^n\eta_x(\gamma^-)$ is a lift of $\widehat{\xi}(\gamma^-)$. For every subsequence $n_k$ such that $\rho(\gamma)^{n_k}\eta_x(\gamma^-)$ is constant, it follows that
  \begin{equation}	\label{eq:Rhat_Anosov}
    \pos_{R,R'}(\widehat{\xi}(\gamma^+),\rho(\gamma)^{n_k}\eta_x(\gamma^-)) \leq \dc{w_0}.
  \end{equation}
  But as $\pos_{R,R}(\widehat{\xi}(\gamma^+),\widehat{\xi}(\gamma^-)) = \dc{w_0}$, the position in \eqref{eq:Rhat_Anosov} must be a transverse one, thus equality holds by \autoref{lem:transverse_maximal}. As seen before, this uniquely determines $\rho(\gamma)^{n_k}\eta_x(\gamma^-)$ among the lifts of $\widehat{\xi}(\gamma^-)$. Since the same holds for any subsequence $n_k$ such that $\rho(\gamma)^{n_k}\eta_x(\gamma^-)$ is constant, $\rho(\gamma)$ fixes $\eta_x(\gamma^-)$ and we obtain
  \begin{equation*}
    \pos_{R,R'}(\widehat{\xi}(\gamma^+),\eta_x(\gamma^-)) = \dc{w_0}.
  \end{equation*}
  Applying the same argument to $y \neq \gamma^-$ shows that $\eta_x(\gamma^-) = \eta_y(\gamma^-)$.\\
  Therefore, $\eta_x$ and $\eta_y$ are continuous functions on $\bdry \setminus \{x,y\}$ which agree on the dense subset of poles, hence they agree everywhere. We denote by
  \[ \eta \colon \bdry \to \F_{R'} \]
  the continuous function defined by $\eta(y) = \eta_x(y)$ for any choice of $x\neq y$. It is $\rho$--equivariant because $\eta(\gamma y) = \eta_{\gamma x}(\gamma y) \in \F_{R'}$ is the lift of $\widehat{\xi}(\gamma y)$ defined by
  \[ \pos_{R,R'}(\widehat{\xi}(\gamma x),\eta_{\gamma x}(\gamma y)) = \pos_{R,R'}(\rho(\gamma)\widehat{\xi}( x),\eta_{\gamma x}(\gamma y)) = \dc{w_0}, \]
  which is $\rho(\gamma) \eta_x(y) = \rho(\gamma) \eta(y)$.
\end{Prf}
\begin{Rem}
  It is worth noting that the independence of $\eta_x(z)$ of the point $x$ simplifies greatly if $\bdry$ is connected: If $x$ and $y$ can be connected by a path $x_t$ in $\bdry$, we consider the lifts $\eta_{x_t}(z)$ along the path. They need to be constant by continuity of $\widehat{\xi}$, so $\eta_x(z)$ and $\eta_y(z)$ agree.
\end{Rem}
\begin{Ex}
  Let us illustrate \autoref{prop:R-Anosov_implies_Rhat-Anosov} with an example. Let $G=\SL(n,\bR)$ and $\rho:\Gamma\to G$ a representation which is $P_\theta$--Anosov with $\theta = \{\alpha_2,\dots,\alpha_{n-2} \}$, so that we have a boundary map $\xi:\bdry\to\F_{1,n-1}$ into the space of partial flags comprising a line and a hyperplane. Assume that $\rho$ is $P_R$--Anosov, where $R=\langle \I(\theta), \I(\alpha_{n-1})^2 \rangle$. Then there is a boundary map $\widehat{\xi}$ into the space $\F_R$ of flags comprising an oriented line and an unoriented hyperplane. Let $x,z \in \bdry$ be two points as in the proof of the proposition. We can fix an orientation on $\widehat{\xi}(z)^{(n-1)}$ by requiring that $(\widehat{\xi}(x)^{(1)},\widehat{\xi}(z)^{(n-1)})$, written in this order, induces the standard orientation on $\bR^n$ (or the opposite orientation, depending on which element $w_0 \in \widetilde{W}$ we chose to represent the transversality type $\dc{w_0}$ of $\widehat{\xi}$). Doing so for all points $z\in\bdry$ extends the boundary map to a map into the space $\F_{R'}$ of flags comprising an oriented line and an oriented hyperplane.
\end{Ex}
As a consequence of the previous two propositions, the minimal oriented parabolic type associated to a $P_\theta$--Anosov representation automatically has certain properties.
\begin{Prop}	\label{prop:w02_in_E}
  Let $\theta\subsetneq\Delta$ be stable under $\iota$, $R = \left< \I(\theta),E \right>$ an oriented parabolic type, $w_0 \in T$ with $w_0 E w_0^{-1} = E$, and $\rho \colon \Gamma \to G$ be $P_R$--Anosov with transversality type $\dc{w_0} \in \widetilde W_{R,R}$. Then $w_0^2 \in E$.
\end{Prop}
\begin{Prf}
  Let $x\neq y \in \bdry$ be two points in the boundary. Then $\dc{w_0} = \pos_{R,R}(\widehat{\xi}(x),\widehat{\xi}(y))$ where $\widehat \xi \colon \bdry \to F_R$ is the limit map. Then $\pos_{R,R}(\widehat{\xi}(y),\widehat{\xi}(x)) = \dc{w_0^{-1}}$. As observed in the proof of \autoref{lem:transversality_type}, there is a dense orbit in $\bdry \times \bdry$; let $(a,b)$ be an element of this orbit. Since we can approximate both $(x,y)$ and $(y,x)$ by this orbit, continuity of $\widehat{\xi}$ implies that $\dc{w_0} \leq \pos_{R,R}(\widehat{\xi}(a),\widehat{\xi}(b))$ and $\dc{w_0^{-1}} \leq \pos_{R,R}(\widehat{\xi}(a),\widehat{\xi}(b))$. All of these are transverse positions, thus equality must hold in both cases by \autoref{lem:transverse_maximal} and we conclude $\dc{w_0} = \dc{w_0^{-1}} \in \widetilde{W}_{R,R}$. Since $w_0 R w_0^{-1} = R$, this implies $w_0^2 \in R \cap \Mbar = E$.
\end{Prf}
The final part of this chapter is aimed at distinguishing connected components of Anosov representations by comparing the possible lifts of the limit map.
\begin{Prop}\label{prop:clopen}
  The set of $P_R$--Anosov representations is open and closed in the space of $P_\theta$--Anosov representations $\Hom_{\text{$P_\theta$--Anosov}}(\Gamma, G) \subset \Hom(\Gamma, G)$.
\end{Prop}
To prove this proposition, we will make use of the following technical lemma. Choose an auxiliary Riemannian metric on $\F_\theta$ and equip $\F_R$ with the metric which makes the finite covering $\pi_R \colon \F_R \to \F_\theta$ a local isometry.
\begin{Lem}\label{lem:local_lift}
  Let $\widehat\xi \colon \bdry \to \F_R$ be a limit map of an $P_R$--Anosov representation and $\xi = \pi_R \circ \widehat\xi$. Then there exists $\delta > 0$ such that for every $x \in \bdry$
  \begin{enumerate}
  \item $\pi_R^{-1}(B_\delta(\xi(x))) = \displaystyle\bigsqcup_{[m] \in \Mbar/E} B_\delta(R_m(\widehat\xi(x)))$, and $\pi_R$ maps any of these components isometrically to $B_\delta(\xi(x))$,
  \item and the set $R_m(B_\delta(\widehat\xi(x))) \subset \F_R$ intersects $\widehat\xi(\bdry)$ if and only if $[m] = 1 \in \Mbar/E$.
  \end{enumerate}
\end{Lem}
\begin{Prf}
  By compactness of $\F_\theta$ there is an $\varepsilon > 0$ such that, for every $f \in \F_\theta$, the preimage of $B_\varepsilon(f)$ under $\pi_R$ is the disjoint union of $\varepsilon$--balls around the preimages of $x$. Together with the choice of metric on $\F_R$, this shows \itemnr{1} for any $\delta \leq \varepsilon$.

  Now for every $x \in \bdry$ define
  \[\mathcal R_x = \xi(\widehat\xi^{-1}(\F_R \setminus B_\varepsilon(\widehat\xi(x)))), \qquad \delta_x = \min\{\varepsilon, \ \textstyle\frac{1}{2}d(\xi(x), \mathcal R_x)\}. \]
  This is positive since $\mathcal R_x \subset \F_\theta$ is closed and $\xi$ is injective. By compactness there is a finite collection $x_1, \dots, x_m \in \bdry$ such that the sets $B_{\delta_{x_i}}(\xi(x_i)) \subset \F_\theta$ cover $\xi(\bdry)$. Let $\delta = \min_i \delta_{x_i}$. Then $U = B_\delta(\xi(x)) \subset B_\varepsilon(\xi(x))$ for every $x \in \bdry$, so $\pi_R^{-1}(U)$ decomposes into disjoint $\delta$--balls as in \itemnr{1}. One of these is $V = B_\delta(\widehat\xi(x))$, and it is indeed the only one intersecting $\widehat\xi(\bdry)$:

  If $y \in \bdry$ with $\widehat\xi(y) \in \pi_R^{-1}(U)$, then $\xi(y) \in U = B_{\delta}(\xi(x)) \subset B_{2\delta_{x_i}}(\xi(x_i))$ for some $i$. So
  \[d(\xi(x_i), \xi(y)) < 2\delta_{x_i} \leq d(\xi(x_i), \mathcal R_{x_i}),\]
  thus $\xi(y) \not\in \mathcal R_{x_i}$ or equivalently $\widehat\xi(y) \in B_\varepsilon(\widehat\xi(x_i))$. So $\widehat\xi(y) \in \pi_R^{-1}(U) \cap B_\varepsilon(\widehat\xi(x_i))$, which is exactly $V$.
\end{Prf}
\begin{Prf}[of \autoref{prop:clopen}]
  To show openness, let $\rho_0$ be $P_R$--Anosov with limit map $\widehat\xi_0 \colon \bdry \to \F_R$ and $\xi_0 = \pi_R \circ \widehat\xi_0$. Let $\delta$ be the constant from \autoref{lem:local_lift} for $\widehat\xi_0$. Choose $x_1, \dots, x_k \in \bdry$ such that $B_{\delta/4}(\xi_0(x_i))$ cover $\xi_0(\bdry)$ and let $U_i = B_{\delta/2}(\xi_0(x_i))$ and $V_i = B_{\delta/2}(\widehat\xi_0(x_i))$. For every $i$ we get a local section $s_i \colon U_i \to V_i$. If $U_i$ and $U_j$ intersect, then $s_i$ and $s_j$ coincide on the intersection, since $U_i \cup U_j$ is contained in a $\delta$--ball, of which only a single lift can intersect $\widehat\xi_0(\bdry)$, so $V_i$ and $V_j$ both have to be contained in this lift. Therefore, the $s_i$ combine to a smooth section $s \colon \bigcup_i U_i \to \bigcup_i V_i$.

  For every $\rho_1 \in \Hom_{\text{$P_\theta$--Anosov}}(\Gamma,G)$ which is close enough to $\rho_0$, there is a path $\rho_t \in \Hom_{\text{$P_\theta$--Anosov}}(\Gamma,G)$ connecting $\rho_0$ and $\rho_1$ such that $d_{C^0}(\xi_t, \xi_0) < \delta/4$ for every $t \in [0,1]$. This is because $\Hom_{\text{$P_\theta$--Anosov}}(\Gamma,G)$ is open and the limit map depends continuously on the representation \cite[Theorem 5.13]{GuichardWienhardDomains}. Then for every $x \in \bdry$ there is an $i$ such that
  \[d(\xi_1(x), \xi_0(x_i)) \leq d_{C^0}(\xi_1, \xi_0) + d(\xi_0(x), \xi_0(x_i)) < \delta/2,\]
  hence $\xi_1(x) \in U_i$. So $\xi_1(\bdry) \subset \bigcup_iU_i$ and we can define $\widehat\xi_1 = s \circ \xi_1$. This is a continuous lift of $\xi_1$. Note that also $\widehat\xi_0 = s \circ \xi_0$ and that we can equally define $\widehat\xi_t = s \circ \xi_t$ for every $t \in [0,1]$.

  To show $\rho_1$--equivariance of $\widehat\xi_1$, let $\gamma \in \Gamma$, $x \in \bdry$ and consider the curves
  \[ \alpha(t) = \rho_t(\gamma)^{-1} \widehat\xi_t(\gamma x), \quad \beta(t) = \widehat\xi_t(x).\]
  They are continuous and $\pi_R(\alpha(t)) = \rho_t(\gamma)^{-1} \pi_R(\widehat\xi_t(\gamma x)) = \xi_t(x) = \pi_R(\beta(t))$. Also $\alpha(0) = \rho_0(\gamma)^{-1}\widehat\xi_0(\gamma x) = \widehat\xi_0(x) = \beta(0)$ by $\rho_0$--equivariance of $\xi_0$. Therefore, the curves $\alpha$ and $\beta$ coincide, so in particular $\widehat\xi_1$ is $\rho_1$--equivariant.

  For closedness, let $\rho_n$ be a sequence of $P_R$--Anosov representations with limit maps $\widehat\xi_n$ converging to the $P_\theta$--Anosov representation $\rho$. Then the unoriented limit maps $\xi_n = \pi_R \circ \widehat\xi_n$ converge uniformly to $\xi$, the limit map of $\rho$. Let $\gamma \in \Gamma$ be an element of infinite order and $\gamma^-,\gamma^+ \in \bdry$ its poles. Since $\pi_R$ is a finite covering, up to taking a subsequence, we can assume that $\widehat\xi_n(\gamma^+)$ converges to a point we call $\widehat\xi(\gamma^+)$. First, we are going to show that there is a neighborhood of $\gamma^+$ in $\bdry$ on which the maps $\widehat\xi_n$ converge uniformly to some limit.

  As $\rho$ is Anosov, the points $\xi(\gamma^-),\xi(\gamma^+) \in \F_\theta$ are transverse fixed points of $\rho(\gamma)$. Since $\xi_n(\gamma^\pm) \to \xi(\gamma^\pm)$, we can find an $\epsilon > 0$ such that all elements of $\overline{B_\epsilon(\xi(\gamma^+))}$ are transverse to $\xi_n(\gamma^-)$ for sufficiently large $n$. In particular, $\rho_n(\gamma^k)$ restricted to this ball converges uniformly to (the constant function with value) $\xi_n(\gamma^+)$ as $k \to \infty$. After shrinking $\epsilon$, the preimage $\pi_R^{-1}(B_\epsilon(\xi(\gamma^+)))$ is a union of finitely many disjoint copies of $B_\epsilon(\xi(\gamma^+))$:
  \[ \pi_R^{-1}(B_\epsilon(\xi(\gamma^+))) = \bigsqcup_{[m] \in \Mbar/E} B_\epsilon(\widehat\xi(\gamma^+))m \]
  For $n$ large, $\widehat\xi_n(\gamma^+) \in B_\epsilon(\widehat\xi(\gamma^+))$. Furthermore, for large $k$, $\rho_n(\gamma^k)$ maps $B_\varepsilon(\xi(\gamma^+))$ into itself. Since $\rho_n(\gamma^k)B_\varepsilon(\widehat\xi(\gamma^+)) \subset \F_R$ is connected and contains $\widehat\xi_n(\gamma^+)$, it must be inside $B_\varepsilon(\widehat\xi(\gamma^+))$. So as $\rho_n(\gamma^k)|_{B_\varepsilon(\xi(\gamma^+))} \to \xi_n(\gamma^+)$ uniformly for $k \to \infty$, when seen as maps on $\F_R$, $\rho_n(\gamma^k)|_{B_\epsilon(\widehat\xi(\gamma^+))}$ also converges uniformly to $\widehat\xi_n(\gamma^+)$. Now choose $\delta > 0$ such that $\xi(B_\delta(\gamma^+)) \subset B_{\varepsilon/2}(\xi(\gamma^+))$ and $\gamma^- \not\in B_\delta(\gamma^+)$. We claim that then $\widehat\xi_n|_{B_\delta(\gamma^+)}$ converges uniformly to a lift of $\xi|_{B_\delta(\gamma^+)}$.

  To see this, let $n$ be large enough so that $d_{C^0}(\xi_n,\xi) < \epsilon / 2$. Then $\xi_n(B_\delta(\gamma^+)) \subset B_\epsilon(\xi(\gamma^+))$. Let $y \in \widehat\xi_n(B_\delta(\gamma^+))$ be any point, and let $m\in \Mbar$ be chosen such that $y \in B_\epsilon(\widehat\xi(\gamma^+))m$. It follows that $\rho_n(\gamma^k)(y) \xrightarrow{k\to\infty} \widehat\xi_n(\gamma^+)m$. So by $\rho(\Gamma)$--invariance and closedness of $\widehat\xi_n(\bdry)$, $\widehat\xi_n(\gamma^+)m \in \widehat\xi_n(\bdry)$, so $[m] = 1 \in \Mbar/E$ by transversality. Thus for all sufficiently large $n$, the image of $\widehat\xi_n|_{B_\delta(\gamma^+)}$ is entirely contained in $B_\epsilon(\widehat\xi(\gamma^+))$, so we can use the local section $s: B_\epsilon(\xi(\gamma^+)) \to B_\epsilon(\widehat\xi(\gamma^+))$ to write $\widehat\xi_n|_{B_\delta(\gamma^+)} = s \circ \xi_n|_{B_\delta(\gamma^+)}$. This proves the stated uniform convergence on $B_\delta(\gamma^+)$.

  Now we use local uniform convergence at $\gamma^+$ to obtain uniform convergence everywhere. For any point $y \in\bdry\setminus\{\gamma^-\}$ and any neighborhood $U \ni y$ whose closure does not contain $\gamma^-$, there is an integer $k(U)$ such that $\gamma^k(U) \subset B_\delta(\gamma^+)$ for all $k\geq k(U)$ \cite[Theorem 4.3]{BenakliKapovich}. Then, for $z\in U$,
  \[ \widehat\xi_n(z) = \rho_n(\gamma)^{-k(U)}\widehat\xi_n(\gamma^{k(U)}z) \xrightarrow{n\to\infty} \rho(\gamma)^{-k(U)}\widehat\xi(\gamma^{k(U)}z), \]
  so we get local uniform convergence on $\bdry \setminus \{\gamma^-\}$. Similarly, since $\widehat\xi_n = \rho_n(\zeta)^{-1} \circ \widehat\xi \circ \zeta$ for some $\zeta \in \Gamma$ with $\zeta \gamma^- \neq \gamma^-$, $\widehat \xi_n$ also converges uniformly in a neighborhood $\gamma^-$. So the maps $\widehat\xi_n$ converge uniformly to a limit $\widehat\xi$, which is continuous and equivariant.
\end{Prf}
From the previous proposition, we obtain the following two criteria to distinguish connected components of Anosov representations.
\begin{Cor}	\label{cor:connected_components}
  Let $\rho,\rho' \colon \Gamma \to G$ be $P_\theta$--Anosov. Furthermore, let $R,R' \subset \widetilde{W}$ be the minimal oriented parabolic types such that $\rho$ is $P_R$--Anosov and $\rho'$ is $P_{R'}$--Anosov (see \autoref{prop:minimal_type}). Assume that $\rho$ and $\rho'$ lie in the same connected component of $\Hom_{\text{$P_\theta$--Anosov}}(\Gamma,G)$. Then the types $R$ and $R'$ agree. Furthermore, if $\widehat{\xi},\widehat{\xi}' \colon \bdry \to \F_R$ are limit maps of $\rho,\rho'$ of transversality types $\dc{w_0},\dc{w_0'} \in \widetilde{W}_{R,R}$, then $\dc{w_0},\dc{w_0'}$ are conjugate by an element of $\Mbar$.
\end{Cor}
\begin{Prf}
  By \autoref{prop:clopen}, $\rho$ is also $P_{R'}$--Anosov and $\rho'$ is $P_R$--Anosov. If $R$ and $R'$ were not equal, either $R'$ would not be minimal for $\rho'$ or $R$ would not be minimal for $\rho$.\\
  By \autoref{rem:boundary_maps}\itemnr{2}, the transversality type of any limit map $\xi_\rho \colon \bdry \to \F_R$ of $\rho$ is conjugate to $\dc{w_0}$ by an element of $\Mbar$. By (the proof of) \autoref{prop:clopen}, $\rho$ also admits a limit map of transversality type $\dc{w_0'}$, so they must be conjugate by an element of $\Mbar$.
\end{Prf}
\section{Domains of discontinuity} \label{sec:domains_of_discontinuity}
In this section, we extend the machinery developed in \cite{KapovichLeebPortiFlagManifolds} to the setting of oriented flag manifolds (\autoref{def:oriented_flag_manifold}). More specifically, we show that their description of cocompact domains of discontinuity for the action of Anosov representations on flag manifolds can be applied with some adjustments to oriented flag manifolds. Our main result is the following theorem, which is analogous to \cite[Theorem 7.14]{KapovichLeebPortiFlagManifolds}:\ 

\begin{Thm} \label{thm:prop_disc_cocpt}
  Let $\Gamma$ be a non--elementary word hyperbolic group and $G$ a connected, semi--simple, linear Lie group. Furthermore, let $R,S \subset \widetilde{W}$ be oriented parabolic types and $w_0 \in T\subset\widetilde{W}$ a transverse position such that $w_0Rw_0^{-1} = R$ and $w_0^2 \in R$.\\
  Let $\rho \colon \Gamma \to G$ be a $P_R$-Anosov representation and $\xi: \bdry \to \F_R$ a limit map of transversality type $\dc{w_0} \in \widetilde{W}_{R,R}$. Assume that $I \subset \widetilde{W}_{R,S}$ is a $w_0$--balanced ideal, and define $\mathcal{K} \subset \F_S$ as
  \[ \mathcal{K} = \bigcup\limits_{x\in\bdry} \bigcup\limits_{\dc{w}\in I} C_{\dc{w}}(\xi(x)). \]
  Then $\mathcal{K}$ is $\Gamma$--invariant and closed, and $\Gamma$ acts properly discontinuously and cocompactly on the domain $\Omega = \F_S \setminus \mathcal{K}$.
\end{Thm}
We first observe that this theorem implies \autoref{thm:domains_introduction}: By \autoref{prop:R-Anosov_implies_Rhat-Anosov}, when starting with a $P_R$--Anosov representation such that $w_0 R w_0^{-1} \neq R$, then it is actually $P_{R'}$--Anosov, with $R' = R \cap w_0 R w_0^{-1}$. Then the conditions $w_0 R' w_0^{-1} = R'$ and by \autoref{prop:w02_in_E} also $w_0^2 \in R$ are automatically satisfied. Also, if a balanced ideal is invariant by $R$ from the left and $S$ from the right, then it is also invariant by $R'$ and gives a balanced ideal in $\widetilde W_{R',S}$.

A large part of the work required to prove this version, namely extending the Bruhat order to the extended Weyl group $\widetilde{W}$, was already done in and \ref{sec:oriented_relative_positions} and \autoref{sec:B0B0action}. We prove proper discontinuity and cocompactness of the action of $\Gamma$ on $\Omega$ separately in the following two subsections (Theorems \ref{thm:prop_disc} and \ref{thm:cocpt}). The part about cocompactness follows \cite{KapovichLeebPortiFlagManifolds} in all key arguments. Since oriented flag manifolds are not as established and well--studied as their unoriented counterparts, we reprove all the required technical lemmas in the setting of compact $G$-homogeneous spaces $X,Y$ and $G$--equivariant maps between $X$ and $\mathcal{C}(Y)$, the space of closed subsets of $Y$.
\subsection{Proper discontinuity}	\label{sec:proper_discontinuity}
Let $P_R$ and $P_S$ be oriented parabolic subgroups of types $R = \langle \I(\theta), E \rangle$ and $S = \langle \I(\eta), F \rangle$. Furthermore, let $w_0 \in T\subset\widetilde W$ be a transverse position. We assume that $\iota(\theta) = \theta$, $w_0 E w_0^{-1} = E$ and $w_0^2 \in E$, so that $w_0$ acts involutively on $\widetilde W_{R,S}$ (see \autoref{sec:involutions}).

The following definition of $w_0$--related limits is an oriented version of the one used for contracting sequences in \cite[Definition 6.1]{KapovichLeebPortiFlagManifolds}. The idea goes back to the study of discrete quasiconformal groups in \cite{GehringMartin}. Apart from the dependence on the choice of $w_0$, we will see later that pairs of such limits are not unique in this setting (\autoref{lem:dynamics_oriented}).

\begin{Def}
  Let $(g_n) \in G^\bN$ be a diverging sequence. A pair $F^-,F^+ \in \F_R$ is called a pair of \emph{$w_0$--related limits} of the sequence $(g_n)$ if
  \[ g_n|_{C_{\dc{w_0}}(F^-)} \xrightarrow{n\to\infty} F^+ \]
  locally uniformly.
\end{Def}
We will also make use of the notion of dynamically related points, which was given in \cite[Definition 1]{Frances}. It is easily proved that an action of a discrete group by homeomorphisms is proper if and only if there are no dynamically related points.
\begin{Def}
  Let $X$ be a topological space. Two points $x,y \in X$ are called \emph{dynamically related via the sequence $(g_n) \in \mathrm{Homeo}(X)^\bN$} if $(g_n)$ is a divergent sequence and there exists a sequence $x_n \to x$ such that
  \[ g_n(x_n) \to y. \]
\end{Def}
Using similar arguments as in the unoriented case, we can prove the following useful relative position inequality.
\begin{Lem}[{\cite[Proposition 6.5]{KapovichLeebPortiFlagManifolds}}]\label{lem:dynamical_relation}
  Let $(g_n) \in G^\bN$ be a sequence admitting a pair $F^\pm \in \F_R$ of $w_0$--related limits. Assume that $F,F' \in \F_S$ are dynamically related via $(g_n)$. Then
  \[ \pos_{R,S}(F^+,F') \leq w_0 \pos_{R,S}(F^-,F). \]
\end{Lem}
\begin{Prf}
  Let $(F_n) \in \F_S^\bN$ be a sequence such that $F_n \to F$ and $g_n F_n \to F'$. We pick elements $h_n \in G$ satisfying $F_n = h_n F$ and $h_n \to 1$. Writing $\dc{w} = \pos_{R,S}(F^-,F)$, it follows that there exists some $g\in G$ such that $g(F^-,F) = ([1],[w])$. Define $f \in \F_R$ as $f = [g^{-1}w_0]$, so that we obtain the following relative positions:
  \begin{itemize}
  \item $\pos_{R,R}(F^-,f) = \dc{w_0}$
  \item $\pos_{R,S}(F^-,F) = \dc{w}$
  \item $\pos_{R,S}(f,F) = \pos_{R,S}([w_0],[w]) = \dc{w_0w}$
  \end{itemize}
  In other words, $f$ is chosen such that $\pos_{R,R}(F^-,f) = \dc{w_0}$ and $\pos_{R,S}(f,F)$ is as small as possible. Then, since $h_n f \to f$, $f$ lies in the open set $C_{\dc{w_0}}(F^-)$ and $F^\pm$ are $w_0$--related limits, it follows that $g_n h_n f \to F^+$. Finally, observe that
  \[ \pos_{R,S}(g_n h_n f,g_n h_n F) = \pos_{R,S}(f,F) \]
  is constant. We thus obtain the following inequalities:
  \begin{equation*}
    \pos_{R,S}(F^+,F') \leq \pos_{R,S}(g_n h_n f, g_n h_n F) = \pos_{R,S}(f,F) = \dc{w_0w} \qedhere
  \end{equation*}
\end{Prf}
One consequence of this inequality is that being $w_0$--related limits is a symmetric condition.
\begin{Lem}[{\cite[(6.7)]{KapovichLeebPortiFlagManifolds}}]\label{lem:inverse_dynamics}
  If $(F^-, F^+)$ is a pair of $w_0$--related limits in $\F_R$ of a sequence $(g_n)$ then $(F^+, F^-)$ is a pair of $w_0$--related limits of $(g_n^{-1})$.
\end{Lem}
\begin{Prf}
  Let $F_n \to F$ be a convergent sequence in $C_{\dc{w_0}}(F^+) \subset \F_R$ and $g_{n_k}^{-1}F_{n_k} \to F' \in \F_R$ a convergent subsequence of $g_n^{-1}F_n$. This means $F$ is dynamically related to $F'$ via $(g_{n_k}^{-1})$ or equivalently $F'$ is dynamically related to $F$ via $(g_{n_k})$. So by \autoref{lem:dynamical_relation} (with $S = R$)
  \[\dc{w_0} = \pos_{R,R}(F^+, F) \leq w_0 \pos_{R,R}(F^-, F').\]
  Since $\dc{w_0}$ is maximal in the Bruhat order, this implies that $w_0 \pos_{R,R}(F^-, F') = \dc{w_0}$ by \autoref{lem:order_projection}. As $w_0$ induces an involution on $\widetilde{W}_{R,R}$, we obtain $\pos_{R,R}(F^-, F') = \dc{1}$, i.e. $g_{n_k}F_{n_k} \to F' = F^-$. By the same argument every subsequence of $g_n^{-1}F_n$ accumulates at $F^-$ and thus $g_n^{-1}F_n \to F^-$, which shows that $(F^+,F^-)$ are $w_0$--related limits of $(g_n^{-1})$.
\end{Prf}
In the unoriented case, a $P_\theta$--divergent sequence admits subsequences with unique attracting limits in $\F_\theta$. In the oriented case, however, this uniqueness is lost and all lifts of such a limit will be attracting on an open set.
\begin{Lem}	\label{lem:dynamics_oriented}
  Let $(g_n) \in G^\bN$ be a $P_\theta$--divergent sequence. Then there is a subsequence $(g_{n_k})$ admitting $|\Mbar/E|$ pairs of $w_0$--related limits in $\F_R$. More precisely, the action
  \begin{equation}	\label{eq:related_pair_action}
    \Mbar/E \times \F_R^2 \to \F_R^2, \quad ([m], (F^-, F^+)) \mapsto (R_{w_0mw_0^{-1}}(F^-), R_{m}(F^+))
  \end{equation}
  is simply transitive on the pairs of $w_0$--related limits of $(g_{n_k})$.
\end{Lem}
\begin{Prf}
  Observe that since $w_0 E w_0^{-1} = E$, conjugation by $w_0$ defines an action on $\Mbar/E$. In other words, the choice of the representative $m\in \Mbar$ in \eqref{eq:related_pair_action} does not matter.\\
  Let us first prove that $\Mbar/E$ acts simply transitively on the $w_0$--related limits of $(g_{n_k})$, assuming such limits exist. We know from \autoref{cor:position_rightmult} that
  \[C_{\dc{w_0}}(R_{w_0mw_0^{-1}}(F^-)) = C_{\dc{w_0m}}(F^-) = R_m(C_{\dc{w_0}}(F^-)).\]
  Because of this and since left and right multiplication commute, \eqref{eq:related_pair_action} restricts to an action on $w_0$--related limits of $(g_{n_k})$. It is free by the definition of $E$ and $\F_R$. For transitivity, let $F^\pm$ and ${F'}^\pm$ be two $w_0$--related limit pairs for $(g_{n_k})$. Then
  \[\bigcup_{[m] \in \Mbar/E} R_m(C_{\dc{w_0}}(F^-)) = \bigcup_{[m] \in \Mbar/E} C_{\dc{w_0m}}(F^-)\]
  is dense in $\F_R$ since its closure is all of $\F_R$ by \autoref{prop:Bruhat_order_quotient}. So for some $[m] \in \Mbar/E$, $R_m(C_{\dc{w_0}}(F^-))$ must intersect the open set $C_{\dc{w_0}}({F'}^-)$. On this intersection, $g_{n_k}$ converges locally uniformly to $R_m(F^+)$ and ${F'}^+$, so ${F'}^+ = R_m(F^+)$. By \autoref{lem:inverse_dynamics}, $({F'}^+, {F'}^-)$ and $(R_m(F^+), R_{w_0mw_0^{-1}}(F^-))$ are $w_0$--related limits for the sequence $(g_{n_k}^{-1})$, so on $C_{\dc{w_0}}({F'}^+) = C_{\dc{w_0}}(R_m(F^+))$ it locally uniformly converges to both ${F'}^-$ and $R_{w_0mw_0^{-1}}(F^-)$, so ${F'}^- = R_{w_0mw_0^{-1}}(F^-)$.

  What is left to show is the existence of $w_0$--related limits. This is done by an argument similar to \cite[Lemma 4.7]{GuichardKasselWienhard}. Decompose the sequence $g_n$ as $g_n = k_n e^{A_n} \ell_n$ with $k_n, \ell_n \in K$ and $A_n \in \aplusbar$. After taking a subsequence, we can assume that $k_n \to k$ and $\ell_n \to \ell$. We want to show that $F^- = [\ell^{-1} w_0^{-1}] \in \F_R$ and $F^+ = [k] \in \F_R$  are $w_0$--related limits of $(g_n)$. We use the following characterization of locally uniform convergence: For every sequence $F_n \to F$ converging inside $C_{\dc{w_0}}(F^-)$ we want to show that $g_n F_n$ converges to $F^+$. For sufficiently large $n$ the sequence $\ell_n F_n$ will be inside $C_{\dc{w_0}}(\ell F^-) = C_{\dc{w_0}}([w_0^{-1}])$, so by \autoref{lem:transverse_parametrization} we can write $\ell_n F_n = [e^{X_n}]$ with
  \[X_n \in \!\!\!\bigoplus_{\alpha \in \Sigma^- \setminus \Span (\theta)} \!\!\fg_\alpha\]
  converging to some $X$ from the same space. So
  \[g_n F_n = [k_n e^{A_n} e^{X_n}] = [k_n e^{A_n} e^{X_n} e^{-A_n}] = [k_n \exp(\Ad_{e^{A_n}} X_n)] = [k_n \exp(e^{\ad A_n} X_n)].\]
  If we decompose $X_n = \sum_\alpha X_n^\alpha$ into root spaces then
  \[e^{\ad A_n} X_n = \sum_{\alpha \in \Sigma^- \setminus \Span(\theta)} e^{\alpha(A_n)} X_n^\alpha.\]
  Now every $\alpha \in \Sigma^- \setminus \Span(\theta)$ can be written as a linear combination of simple roots with non--positive coefficients and with the coefficient of at least one simple root $\beta \in \ctheta$ being strictly negative. As $\beta(A_n) \to \infty$ by $P_\theta$--divergence, $\alpha(A_n)$ must converge to $-\infty$ and therefore $e^{\ad A_n} X_n$ goes to $0$. This implies $g_n F_n \to [k] = F^+$, so $g_n|_{C_{\dc{w_0}}(F^-)} \to F^+$ locally uniformly.
\end{Prf}
Let $\Gamma$ be a non--elementary word hyperbolic group and $G$ a connected, semi--simple, linear Lie group (see \autoref{sec:parabolics} for some remarks on these assumptions).
\begin{Lem}	\label{lem:limits_on_curve}
  Let $\rho:\Gamma \to G$ be a $P_R$--Anosov representation and let $\xi: \bdry \to \F_R$ be a continuous, equivariant limit map of transversality type $\dc{w_0}$. Then every sequence $(\gamma_n)$ of distinct elements admits a subsequence $(\gamma_{n_k})$ and points $x,y \in \bdry$ such that $\gamma_{n_k}|_{\bdry \setminus \{x\}} \to y$ locally uniformly. Moreover, for any such subsequence, $(\xi(x),\xi(y))$ is a pair of $w_0$--related limits for the sequence $(\rho(\gamma_{n_k}))$.
\end{Lem}
\begin{Prf}
  The first property is simply the fact that $\Gamma$ acts as a convergence group on $\bdry$ \cite[Lemma 1.11]{BowditchConvergence}.

  To simplify notation, we assume from now on that $\gamma_{n}|_{\bdry \setminus \{x\}} \to y$ locally uniformly. By \autoref{lem:dynamics_oriented}, there exists a subsequence $(\rho(\gamma_{n_k}))$ with $w_0$--related limits $F^\pm \in \F_R$. Then $F^-$ is a lift of the unique repelling limit $\pi(F^-) \in \F_\theta$, and we have $\pi(F^-) \in \pi(\xi(\bdry))$ (see the description of the boundary map in \cite[Theorem 5.3]{GueritaudGuichardKasselWienhard}). By right--multiplying $F^-$ with an element $m \in \Mbar/E$ and $F^+$ with $w_0^{-1}mw_0$ if necessary, we may assume that $F^- = \xi(x)$ for some $x\in\bdry$ (see \autoref{lem:dynamics_oriented}). For any $x\neq z\in\bdry$, we have $\pos_{R,R}(F^-,\xi(z)) = \pos_{R,R}(\xi(x),\xi(z)) = \dc{w_0}$. Since the $w_0$--related attracting limit $F^+ \in \F_R$ is characterized by
  \[ \rho(\gamma_{n_k})|_{C_{\dc{w_0}}(F^-)} \xrightarrow{n\to\infty} F^+, \]
  it follows that $\rho(\gamma_{n_k})(\xi(z)) \xrightarrow{n\to\infty} F^+$. As
  \[ \rho(\gamma_{n_k})(\xi(z)) = \xi(\gamma_{n_k} z) \xrightarrow{n\to\infty} \xi(y), \]
  we obtain $\xi(y) = F^+$.

  The same reasoning also shows that any subsequence $(\gamma_{n_k})$ of $(\gamma_n)$ has a further subsequence $(\gamma_{n_{k_l}})$ such that $(\xi(x),\xi(y))$ is a pair of $w_0$--related limits of $(\rho(\gamma_{n_{k_l}}))$. So it is in fact a pair of $w_0$--related limits of the whole sequence $(\rho(\gamma_n))
  $.
\end{Prf}
Recall from \autoref{sec:Bruhat_order} and \autoref{sec:involutions} that a subset $I\subset\widetilde{W}_{R,S}$ is an ideal if $\dc{w}\in I$ and $\dc{w'}\leq \dc{w}$ implies $\dc{w'}\in I$, and it is $w_0$--fat if $\dc{w} \not\in I$ implies $\dc{w_0w} \in I$.
\begin{Thm}	\label{thm:prop_disc}
  Let $\rho: \Gamma \to G$ be an $P_R$--Anosov representation and let $\xi: \bdry \to \F_R$ be a limit map of transversality type $\dc{w_0}$. Assume that $I \subset \widetilde{W}_{R,S}$ is a $w_0$--fat ideal, and define $\mathcal{K} \subset \F_S$ as
  \[ \mathcal{K} = \bigcup\limits_{x\in\bdry} \bigcup\limits_{\dc{w}\in I} C_{\dc{w}}(\xi(x)). \]
  Then $\mathcal{K}$ is $\Gamma$--invariant and closed, and $\Gamma$ acts properly discontinuously on the domain $\Omega = \F_S \setminus \mathcal{K}$.
\end{Thm}
\begin{Prf}
  $\Gamma$--invariance and closedness of $\mathcal{K}$ follows from \autoref{ex:equivariant_map_K} and \autoref{lem:properties_of_K}\itemnr{2}.

  Assume that the action of $\Gamma$ on $\Omega$ is not proper. Then there exist $F, F' \in \Omega$ which are dynamically related by some sequence $(\rho(\gamma_n))$. This sequence is $P_\theta$--divergent and by \autoref{lem:limits_on_curve}, a subsequence admits a pair of $w_0$--related limits of the form $\xi(x^\pm)$, where $x^\pm \in \bdry$. So \autoref{lem:dynamical_relation} shows that
  \begin{equation}\label{eq:position_inequality}
    \pos_{R,S}(\xi(x^+),F') \leq w_0 \pos_{R,S}(\xi(x^-),F).
  \end{equation}
  Since $F,F' \not\in \mathcal{K}$, neither $\pos_{R,S}(\xi(x^+), F')$ nor $\pos_{R,S}(\xi(x^-), F)$ can be in $I$. As $I$ is $w_0$--fat, this implies in particular $w_0 \pos_{R,S}(\xi(x^-), F) \in I$. But $I$ is an ideal, so \eqref{eq:position_inequality} implies $\pos_{R,S}(\xi(x^+), F') \in I$, a contradiction.
\end{Prf}
\subsection{Cocompactness}	\label{sec:cocompactness}
We now come to the cocompactness part of \autoref{thm:prop_disc_cocpt}. Owing to the fact that we want to apply everything to oriented flag manifolds, our setup here is more general than in \cite{KapovichLeebPortiFlagManifolds}. Nevertheless, all the key arguments from that paper still work. This includes in particular the idea of using expansion to prove cocompactness. The connection between (convex) cocompactness and expansion at the limit set was originally observed for Kleinian groups in \cite{Sullivan}.

We included a detailed discussion of these arguments in \autoref{sec:compact_homogeneous} and just state the result here. The following notion of an expanding action was introduced by Sullivan in \cite[§9]{Sullivan}.
\begin{Def}
  Let $Z$ be a metric space, $g$ a homeomorphism of $Z$ and $\Gamma$ a group acting on $Z$ by homeomorphisms.
  \begin{enumerate}
  \item $g$ is \emph{expanding at $z \in Z$} if there exists an open neighbourhood $z \in U \subset Z$ and a constant $c > 1$ (the expansion factor) such that
    \[d(gx, gy) \geq c\,d(x,y)\]
    for all $x, y \in U$.
  \item Let $A \subset Z$ be a subset. The action of $\Gamma$ on $Z$ is \emph{expanding at $A$} if for every $z \in A$ there is a $\gamma \in \Gamma$ which is expanding at $z$.
  \end{enumerate}
\end{Def}
Recall that for any compact metric space $Z$, we denote by $\C(Z)$ the set of compact subsets of $Z$. The key proposition will be the following:
\begin{Prop}[\autoref{prop:expansion_implies_cocompact_appendix}]	\label{prop:expansion_implies_cocompact}
  Let $\rho \colon \Gamma \to G$ be a representation of a discrete group $\Gamma$ into a Lie group $G$, $X$ and $Y$ compact $G$--homogeneous spaces and $\Q \colon X \to \mathcal C(Y)$ a $G$--equivariant map. Let $\Lambda \subset X$ be compact and $\Gamma$--invariant such that the action of $\Gamma$ on $X$ is expanding at $\Lambda$. Assume further that $\Q(\lambda) \cap \Q(\lambda') = \varnothing$ for all distinct $\lambda, \lambda' \in \Lambda$. Then $\Gamma$ acts cocompactly on $\Omega = Y \setminus \bigcup_{\lambda \in \Lambda} \Q(\lambda)$.
\end{Prop}

We apply this to the setting of Anosov representations and oriented flag manifolds, to get the main result of this section. Let $\Gamma$, $G$, $R$, $S$ and $w_0$ and $I$ be as in \autoref{thm:prop_disc_cocpt}. We need to show is that the action of a $P_R$--Anosov representation is expanding at the image of its limit map. This follows easily from the analogous statement in the unoriented setting (see \cite[Theorem 1.1\itemnr{2}]{KapovichLeebPortiCharacterizations}), but for convenience we also give a direct proof here.
\begin{Prop}	\label{prop:Anosov_expanding}
  Let $\rho:\Gamma \to G$ be a $P_R$--Anosov representation and $\xi: \bdry \to \F_R$ a limit map of transversality type $\dc{w_0} \in \widetilde W_{R,R}$. Then the action of $\Gamma$ on $\F_R$ is expanding at $\xi(\bdry)$.
\end{Prop}
\begin{Prf}
  First note that expansion does not depend on the choice of Riemannian metric on $\F_R$. To get rid of some constants, we will use a $K$--invariant metric for this proof.

  Fix $x \in \bdry$. Since $\Gamma$ is a non--elementary hyperbolic group, $x$ is a conical limit point, meaning there is a sequence $(\gamma_n) \in \Gamma^\bN$ and distinct points $b, c \in \bdry$ such that $\gamma_n|_{\bdry \setminus \{x\}} \to b$ locally uniformly and $\gamma_n x \to c$ (see \cite[Proposition 1.13]{BowditchConvergence} and \cite[Theorem 1A]{Tukia}). We will show that $\rho(\gamma_n)$ is expanding at $\xi(x)$ for (some) large enough $n$.

  Decompose $\rho(\gamma_n^{-1}) = k_n a_n \ell_n$ with $k_n,\ell_n \in K$ and $a_n = e^{A_n} \in \exp(\aplusbar)$. After replacing $(\gamma_n)$ by a subsequence, we may assume that $k_n$ converges to $k$ and $\ell_n$ converges to $\ell$. As in the proof of \autoref{lem:dynamics_oriented} we see that $([\ell^{-1}w_0^{-1}],[k])$ is a pair of $w_0$--related limits for $(\rho(\gamma_n^{-1}))$. By \autoref{lem:limits_on_curve} and \autoref{lem:inverse_dynamics}, $(\xi(b), \xi(x))$ are also $w_0$--related limits for $(\rho(\gamma_n^{-1}))$. Using the action from \autoref{lem:dynamics_oriented} and possibly modifying the KAK--decomposition accordingly by some $m \in Z_K(\fa)$, we can thus assume that $\xi(b) = [\ell^{-1}w_0^{-1}]$ and $\xi(x) = [k]$.

  Since $b \neq c$, $\xi(c)$ is contained in $C_{\dc{w_0}}(\xi(b))$. This is an open set, so we can choose $\delta > 0$ such that $B_{2\delta}(\xi(c)) \subset C_{\dc{w_0}}(\xi(b))$. As a first step, we prove that $a_n$ is contracting at
  \[\mathcal{C} = \F_R \setminus N_\delta(\F_R \setminus C_{\dc{w_0}}([w_0^{-1}])).\]
  By \autoref{lem:transverse_parametrization}, $X \mapsto [e^X]$ is a diffeomorphism from $\fn_\theta^-$ to $C_{\dc{w_0}}([w_0^{-1}])$. Choose a scalar product on $\fn_\theta^-$ which makes the root spaces orthogonal. It defines a Riemannian metric on $\fn_\theta^-$ and therefore also on $C_{\dc{w_0}}([w_0^{-1}])$. On the compact subset $\mathcal C \subset C_{\dc{w_0}}([w_0^{-1}])$, it is comparable to the $K$--invariant metric on $\F_R$ up to a constant $C$. Now $a_n [\exp(X)] = [a_n \exp(X) a_n^{-1}] = [\exp(e^{\ad A_n}X)]$ for every $X \in \fn_\theta^-$, so the action of $a_n$ on $C_{\dc{w_0}}([w_0^{-1}])$ translates to the linear action of $e^{\ad A_n}$ on $\fn_\theta^-$. So for every $z \in \C$
  \[\|D_za_n\| \leq C \|e^{\ad A_n}\|_{\fn_\theta^-} \leq C e^{-\min\{\alpha(A_n) \mid \alpha \in \Delta \setminus \theta\}}.\]
  Here $\|\cdot\|_{\fn_\theta^-}$ is the operator norm with respect to the chosen norm on $\fn_\theta^-$. The second inequality holds since every root $\beta \in \Sigma^- \mathord\setminus\langle\theta\rangle$ has a strictly negative coefficient for at least one $\alpha \in \ctheta$ and nonpositive coefficients for all $\alpha'\in\Delta$, so $\beta(A_n) \leq -\alpha(A_n)$.

  Let $n$ be large enough such that $\|D_za_n\| < 1$ for all $z \in \C$.
  By a standard argument on the lengths of curves, this implies that every $z \in \C$ has a small neighborhood and a constant $\kappa < 1$ such that $d(a_nv,a_nw) \leq \kappa\,d(v,w)$ for all $v,w$ in this neighborhood. Clearly, this means that $a_n^{-1}$ is expanding on $a_n\C$. As our metric is $K$--invariant, $\rho(\gamma_n) = \ell_n^{-1} a_n^{-1} k_n^{-1}$ is thus expanding at $k_na_n\C = \rho(\gamma_n^{-1})\ell_n^{-1}\C$.

  But since $\ell_n\xi(\gamma_n x) \to \ell \, \xi(c)$ we can, by taking $n$ large enough, assume that $\ell_n\xi(\gamma_n x) \in B_\delta(\ell\,\xi(c))$. Since $B_{2\delta}(\ell\,\xi(c)) \subset C_{\dc{w_0}}([w_0^{-1}])$, we have $B_\delta(\ell\,\xi(c)) \subset \C$. This shows $\xi(x) \in \rho(\gamma_n^{-1})\ell_n^{-1}\C$, so $\rho(\gamma_n)$ is expanding at that point, which is what we wanted to show.
\end{Prf}
\begin{Thm}	\label{thm:cocpt}
  Let $\rho:\Gamma \to G$ be an $P_R$--Anosov representation and $\xi: \bdry \to \F_R$ a limit map of transversality type $\dc{w_0} \in \widetilde W_{R,R}$. Assume that $I \subset \widetilde{W}_{R,S}$ is a $w_0$--slim ideal, and define the set $\mathcal{K} \subset \F_S$ as
  \[ \mathcal{K} = \bigcup\limits_{x\in\bdry} \bigcup\limits_{\dc{w}\in I} C_{\dc{w}}(\xi(x)). \]
  Then $\mathcal{K}$ is $\Gamma$--invariant and closed, and $\Gamma$ acts cocompactly on the domain $\Omega = \F_S \setminus \mathcal{K}$.
\end{Thm}
\begin{Prf}
  $\Gamma$--invariance and closedness of $\mathcal{K}$ follows from \autoref{ex:equivariant_map_K} and \autoref{lem:properties_of_K}\itemnr{2}. As discussed in \autoref{ex:equivariant_map_K}, the map
  \begin{align*}
    \Q \colon \F_R & \to \C(\F_S) \\
    f & \mapsto \bigcup\limits_{\dc{w} \in I} C_{\dc{w}}(f)
  \end{align*}
  is $G$--equivariant. Moreover, $\xi(\bdry)$ is compact and the action of $\Gamma$ on $\F_R$ is expanding at $\xi(\bdry)$ since $\rho$ is $P_R$--Anosov (\autoref{prop:Anosov_expanding}). If we can show that $\Q(\xi(x)) \cap \Q(\xi(x')) = \varnothing$ for $x\neq x'$, cocompactness will follow from \autoref{prop:expansion_implies_cocompact}.\\
  Let $y' \in \Q(\xi(x'))$ be any point. Since $\xi$ has transversality type $\dc{w_0}$ and by definition of $\Q$, we have the relative positions
  \begin{align*}
    \pos_{R,R}(\xi(x),\xi(x')) & = \dc{w_0}, \\
    \pos_{R,S}(\xi(x'),y') & =: \dc{w} \in I.
  \end{align*}
  By \autoref{lem:triangle_ineq}, this implies that $\pos_{R,S}(\xi(x),y') \geq \dc{w_0 w}$. Since $\dc{w_0 w} \not\in I$ by $w_0$--slimness and $I$ is an ideal, $y' \not\in \Q(\xi(x))$.
\end{Prf}
\section{\texorpdfstring{Examples of balanced ideals}{Examples of balanced ideals}}	\label{sec:examples_of_balanced_ideals}
In this section, we will describe explicitly the Bruhat order on $\widetilde W$ and the possible balanced ideals for the group $G = \SL(3,\bR)$. These examples already show how passing from $W$ to $\widetilde W$ vastly increases the number of balanced ideals and therefore the possibilities to build cocompact domains of discontinuity. We have no classification of all balanced ideals, so explicit examples of balanced ideals will be restricted to low dimensions and some special cases in higher dimension.
\subsection{\texorpdfstring{Reduction to $R = \{1,w_0^2\}$ and $S = \{1\}$}{Reduction to R = \{1,w0 2\} and S = \{1\}}}
In applications, we are usually given a fixed representation $\rho \colon \Gamma \to G$. If it is Anosov, then by \autoref{prop:minimal_type}, there is a unique minimal oriented parabolic type $R = \langle \I(\theta), E \rangle$ for $\rho$. If we also fix one of the possible lifts of the boundary map to $\F_R$, then we get a transversality type $\dc{w_0} \in \widetilde W_{R,R}$. To apply \autoref{thm:prop_disc_cocpt} and find cocompact domains of discontinuity in a given flag manifold $\F_S$, we need to look for $w_0$--balanced ideals in $\widetilde W_{R,S}$. Note that this notion of $w_0$--balanced only depends on the equivalence class $\dc{w_0}$.

To enumerate all balanced ideals as we want to do in this section, a different approach is more convenient: We first determine the set $T\subset\widetilde W$ of transverse positions. Then, we want to list all $w_0$--balanced ideals in $\widetilde W_{R,S}$ for $w_0 \in T$ and all possible oriented parabolic types $R = \langle \I(\theta), E \rangle$ and $S = \langle \I(\eta), F \rangle$. For this to be well--defined, $w_0$ must act as an involution on $\widetilde W_{R,S}$, which by \autoref{sec:involutions} happens if $\iota(\theta) = \theta$, $w_0 E w_0^{-1} = E$ and $w_0^2 \in E$.

Note that the smallest $E$ satisfying these conditions is $E = \{ 1,w_0^2 \}$. The following lemma implies that when listing all possible $w_0$--balanced ideals, one can restrict to the minimal choice $R = \{ 1,w_0^2 \}$ and $S = \{ 1 \}$.
\begin{Lem}
  Let $R$ and $S$ be oriented parabolic types as above, and consider the projection $\pi \colon \widetilde W_{\{1,w_0^2\},\{1\}} \to \widetilde W_{R,S}$. Assume that $I \subset \widetilde{W}_{R,S}$ is a $w_0$--balanced ideal. Then $\pi^{-1}(I) \subset \widetilde{W}_{\{1,w_0^2\},\{1\}}$ is a $w_0$--balanced ideal as well.
\end{Lem}
\begin{Prf}
  By \autoref{lem:order_projection}\itemnr{1}, $\pi^{-1}(I)$ is again an ideal. Let $\dc{w} \in \widetilde{W}_{\{1,w_0^2\},\{1\}}$, and recall that $w_0$ acts by left multiplication on both $\widetilde{W}_{\{1,w_0^2\},\{1\}}$ and $\widetilde{W}_{R,S}$, satisfying
  \[ \pi(\dc{w_0  w}) = w_0  \pi(\dc{w}). \]
  Therefore, we obtain the following equivalences:
  \[ \dc{w} \in \pi^{-1}(I) \Leftrightarrow \pi(\dc{w}) \in I \Leftrightarrow w_0  \pi(\dc{w}) \not\in I \Leftrightarrow \dc{w_0  w} \not\in \pi^{-1}(I) \qedhere \]
\end{Prf}
By this lemma, every $w_0$--balanced ideal in $\widetilde{W}_{R,S}$ is obtained by projecting a $R$--left invariant and $S$--right invariant $w_0$--balanced ideal of $\widetilde{W}_{\{1,w_0^2\},\{1\}}$. We can further reduce the number of $w_0$ we have to consider by observing that choices of $w_0$ conjugate by an element in $\Mbar$ lead to essentially the same balanced ideals:
\begin{Lem}
  Let $I \subset \widetilde{W}_{\{1,w_0^2\},\{1\}}$ be a $w_0$--balanced ideal and $m\in \Mbar$. Then $mI$ is a $mw_0m^{-1}$--balanced ideal.
\end{Lem}
\begin{Prf}
  $mI$ is again an ideal by \autoref{lem:position_rightmult}\itemnr{1}. It is $mw_0m^{-1}$--balanced because
  \[ \dc{w} \in mI \Leftrightarrow \dc{m^{-1}w} \in I \Leftrightarrow w_0 \dc{m^{-1}w} \not\in I \Leftrightarrow (mw_0m^{-1}) \dc{w} \not\in mI. \qedhere \]
\end{Prf}
Given a $w_0$--balanced ideal $I$ and an element $m\in \Mbar$ such that $mw_0m^{-1} = w_0$, $mI$ is again $w_0$--balanced. The cocompact domains obtained via \autoref{thm:prop_disc_cocpt} for $I$ and $mI$ are in general different. In contrast to this, the action of $\Mbar$ by right--multiplication is easy to describe: An ideal $I$ is $w_0$--balanced if and only if $Im$ is. Moreover, by \autoref{lem:position_rightmult}, the domain will simply change by global right--multiplication with $m$, i.e. by changing some orientations.
\subsection{\texorpdfstring{The extended Weyl group of $\SL(n,\bR)$}{The extended Weyl group of SL(n,R)}}	\label{sec:SLnR}
Let $G=\SL(n,\bR)$ with maximal compact $K=\SO(n,\bR)$ and $\fa \subset \fsl(n,\bR)$ the set of diagonal matrices with trace $0$. Then $\Sigma = \{\lambda_i - \lambda_j \mid i \neq j\} \subset \fa^*$, where $\lambda_i \colon \fa \to \bR$ is the $i$--th diagonal entry. Choose the simple system $\Delta$ consisting of all roots $\alpha_i \coloneqq \lambda_i - \lambda_{i+1}$ with $i \leq i \leq n-1$. Then $B_0$ is the subgroup of upper triangular matrices with positive diagonal. The group $Z_K(\fa)$ is the group of diagonal matrices with $\pm 1$ entries and $\det =1$. Its identity component is trivial, so $\Mbar = Z_K(\fa)$. The extended Weyl group $\widetilde{W} = N_K(\fa)$ consists of all permutation matrices with determinant 1 -- i.e. all matrices with exactly one $\pm 1$ entry per line and row and all other entries $0$, such that $\det=1$.

A generating set $\I(\Delta)$ in the sense of \autoref{def:generators} is given by
\[\I(\alpha_i) = \begin{pmatrix}I_{i-1} & & & \\ & & -1 & \\ & 1 & & \\ & & & I_{n-i-1}\end{pmatrix}.\]
The transverse positions $T \subset \widetilde W$ are antidiagonal matrices with $\pm 1$ entries. The number of $-1$ entries has to be even if $n$ is equal to $0$ or $1$ mod 4, and odd otherwise. In one formula, it has the same parity as $(n-1)n/2$.

The group $\Mbar$ is generated by diagonal matrices with exactly two $-1$ entries and the remaining entries $+1$. Conjugating $w_0 \in T$ by such an element negates the two lines and the two columns corresponding to the minus signs. This yields the following standard representatives for equivalence classes in $T$ under conjugation by $\Mbar$:
\begin{enumerate}
\item If $n$ is odd, the $(n-1)/2$--block in the upper right corner can be normalized to have $+1$--entries.
\item If $n$ is even, the $(n-2)/2$--block in the upper right corner can be normalized to have $+1$--entries.
\end{enumerate}
If $w_0,w_0' \in T$ of this form are different, they are not conjugate by an element of $\Mbar$.
\subsection{\texorpdfstring{Balanced ideals for $\SL(3,\bR)$}{Balanced ideals for SL(3,R)}}	\label{sec:SL3R}
Let $G = \SL(3,\bR)$ and $\Delta = \{\alpha_1, \alpha_2\}$ the set of simple roots, viewed as their associated reflections. These generate the Weyl group
\[W = \langle \alpha_1,\alpha_2 \mid \alpha_1^2 = \alpha_2^2 = (\alpha_1\alpha_2)^3 = 1 \rangle = \{1, \alpha_1, \alpha_2, \alpha_1\alpha_2, \alpha_2\alpha_1, \alpha_1\alpha_2\alpha_1\}.\]
The Bruhat order on $W$ is just the order by word length and there is a unique longest element $w_0 = \alpha_1\alpha_2\alpha_1$ which acts on $W$ from the left, reversing the order (see \autoref{fig:SL3_unoriented_Bruhat_order}).

\begin{center}
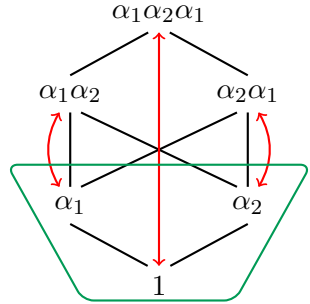

  \begin{tikzpicture}
    \matrix (mtop) [matrix of math nodes,ampersand replacement=\&]{
      \alpha_1\alpha_2\alpha_1 \\
    };

    \matrix (m) [matrix of math nodes,ampersand replacement=\&,row sep=1cm,yshift=-1.8cm]{
      \alpha_1\alpha_2 \&
      \hspace{1cm} \&
      \alpha_2\alpha_1 \\
      \alpha_1 \&
      \hspace{1cm} \&
      \alpha_2 \\
    };

    \matrix (mbot) [matrix of math nodes,ampersand replacement=\&,yshift=-3.6cm]{
      1 \\
    };

    \draw[black,thick] (mtop-1-1.240) -- (m-1-1.north);
    \draw[black,thick] (mtop-1-1.300) -- (m-1-3.north);
    \draw[black,thick] (m-1-3.240) -- (m-2-1.60);
    \draw[black,thick] (m-1-1.300) -- (m-2-3.120);
    \draw[black,thick] (m-1-1.south) -- (m-2-1.north);
    \draw[black,thick] (m-1-3.south) -- (m-2-3.north);
    \draw[black,thick] (m-2-1.south) -- (mbot-1-1.120);
    \draw[black,thick] (m-2-3.south) -- (mbot-1-1.60);

    \draw[red,thick] (mtop-1-1.south) edge [<->] (mbot-1-1.north);
    \draw[red,thick] (m-1-1.240) edge [<->,bend right] (m-2-1.120);
    \draw[red,thick] (m-1-3.300) edge [<->,bend left] (m-2-3.60);

    \draw[ForestGreen,rounded corners,thick] (-2,-2) -- (2,-2) -- (1,-3.8) -- (-1,-3.8) -- cycle;
  \end{tikzpicture}
  \captionof{figure}{The Weyl group of $\SL(3,\bR)$. The black lines indicate the Bruhat order, in the sense that a line going downward from $x$ to $y$ means that $x$ covers $y$ in the Bruhat order. The red arrows show the involution induced by $w_0$. The subset surrounded by the green box is the only balanced ideal.}
  \label{fig:SL3_unoriented_Bruhat_order}
\end{center}

There is only one $w_0$--balanced ideal in this case, which is indicated by the green box in \autoref{fig:SL3_unoriented_Bruhat_order}.

Since $|\Mbar| = 4$, each of the 6 elements of $W$ has 4 preimages in $\widetilde W$, corresponding to different signs in the permutation matrix. The Bruhat order on $\widetilde W$ can be determined using \autoref{prop:Bruhat_order_quotient} and is shown in \autoref{fig:SL3_Bruhat_order}. See \autoref{sec:geometric_interpretation} for a geometric interpretation in terms of oriented flags.

\begin{center}
  \begin{tikzpicture}
    \matrix (mtop) [matrix of math nodes,ampersand replacement=\&,row sep=1cm]{
      \begin{smatrix} && \!\!-1 \\ &  \!\!-1 \\    \!\!-1 \end{smatrix} \&
      \begin{smatrix} && 1      \\ &  \!\!-1 \\    1      \end{smatrix} \&
      \begin{smatrix} && 1      \\ &  1      \\    \!\!-1 \end{smatrix} \&
      \begin{smatrix} && \!\!-1 \\ &  1      \\    1      \end{smatrix} \\
    };

    \matrix (m) [matrix of math nodes,ampersand replacement=\&,row sep=2cm,yshift=-3.5cm]{
      \begin{smatrix} && 1      \\    1      \\ &  1      \end{smatrix} \&
      \begin{smatrix} && \!\!-1 \\    \!\!-1 \\ &  1      \end{smatrix} \&
      \begin{smatrix} && 1      \\    \!\!-1 \\ &  \!\!-1 \end{smatrix} \&
      \begin{smatrix} && \!\!-1 \\    1      \\ &  \!\!-1 \end{smatrix} \&
      \hspace{1cm} \&
      \begin{smatrix} &  1      \\ && \!\!-1 \\    \!\!-1 \end{smatrix} \&
      \begin{smatrix} &  1      \\ && 1      \\    1      \end{smatrix} \&
      \begin{smatrix} &  \!\!-1 \\ && 1      \\    \!\!-1 \end{smatrix} \&
      \begin{smatrix} &  \!\!-1 \\ && \!\!-1 \\    1      \end{smatrix} \\
      \begin{smatrix} &  1      \\    1      \\ && \!\!-1 \end{smatrix} \&
      \begin{smatrix} &  1      \\    \!\!-1 \\ && 1      \end{smatrix} \&
      \begin{smatrix} &  \!\!-1 \\    \!\!-1 \\ && \!\!-1 \end{smatrix} \&
      \begin{smatrix} &  \!\!-1 \\    1      \\ && 1      \end{smatrix} \&
      \hspace{1cm} \&
      \begin{smatrix}    1      \\ && \!\!-1 \\ &  1      \end{smatrix} \&
      \begin{smatrix}    \!\!-1 \\ && 1      \\ &  1      \end{smatrix} \&
      \begin{smatrix}    1      \\ && 1      \\ &  \!\!-1 \end{smatrix} \&
      \begin{smatrix}    \!\!-1 \\ && \!\!-1 \\ &  \!\!-1 \end{smatrix} \\
    };

    \matrix (mbot) [matrix of math nodes,ampersand replacement=\&,row sep=1cm,yshift=-7cm]{
      \begin{smatrix}    1      \\ &  1      \\ && 1      \end{smatrix} \&
      \begin{smatrix}    \!\!-1 \\ &  1      \\ && \!\!-1 \end{smatrix} \&
      \begin{smatrix}    1      \\ &  \!\!-1 \\ && \!\!-1 \end{smatrix} \&
      \begin{smatrix}    \!\!-1 \\ &  \!\!-1 \\ && 1      \end{smatrix} \\
    };

    \draw[red,thick] (mtop-1-1.south) -- (m-1-2.north);
    \draw[red,thick] (mtop-1-1.south) -- (m-1-4.north);
    \draw[red,thick] (mtop-1-2.south) -- (m-1-1.north);
    \draw[red,thick] (mtop-1-2.south) -- (m-1-3.north);
    \draw[red,thick] (mtop-1-3.south) -- (m-1-1.north);
    \draw[red,thick] (mtop-1-3.south) -- (m-1-3.north);
    \draw[red,thick] (mtop-1-4.south) -- (m-1-2.north);
    \draw[red,thick] (mtop-1-4.south) -- (m-1-4.north);

    \draw[ForestGreen,thick] (mtop-1-1.south) -- (m-1-6.north);
    \draw[ForestGreen,thick] (mtop-1-1.south) -- (m-1-8.north);
    \draw[ForestGreen,thick] (mtop-1-2.south) -- (m-1-7.north);
    \draw[ForestGreen,thick] (mtop-1-2.south) -- (m-1-9.north);
    \draw[ForestGreen,thick] (mtop-1-3.south) -- (m-1-6.north);
    \draw[ForestGreen,thick] (mtop-1-3.south) -- (m-1-8.north);
    \draw[ForestGreen,thick] (mtop-1-4.south) -- (m-1-7.north);
    \draw[ForestGreen,thick] (mtop-1-4.south) -- (m-1-9.north);

    \draw[blue,thick] (m-1-6.240) -- (m-2-1.60);
    \draw[blue,thick] (m-1-6.240) -- (m-2-2.60);
    \draw[blue,thick] (m-1-7.240) -- (m-2-1.60);
    \draw[blue,thick] (m-1-7.240) -- (m-2-2.60);
    \draw[blue,thick] (m-1-8.240) -- (m-2-3.60);
    \draw[blue,thick] (m-1-8.240) -- (m-2-4.60);
    \draw[blue,thick] (m-1-9.240) -- (m-2-3.60);
    \draw[blue,thick] (m-1-9.240) -- (m-2-4.60);

    \draw[ForestGreen,thick] (m-1-1.300) -- (m-2-6.120);
    \draw[ForestGreen,thick] (m-1-1.300) -- (m-2-7.120);
    \draw[ForestGreen,thick] (m-1-2.300) -- (m-2-6.120);
    \draw[ForestGreen,thick] (m-1-2.300) -- (m-2-7.120);
    \draw[ForestGreen,thick] (m-1-3.300) -- (m-2-8.120);
    \draw[ForestGreen,thick] (m-1-3.300) -- (m-2-9.120);
    \draw[ForestGreen,thick] (m-1-4.300) -- (m-2-8.120);
    \draw[ForestGreen,thick] (m-1-4.300) -- (m-2-9.120);

    \draw[red,thick] (m-1-1.south) -- (m-2-1.north);
    \draw[red,thick] (m-1-1.south) -- (m-2-4.north);
    \draw[red,thick] (m-1-2.south) -- (m-2-2.north);
    \draw[red,thick] (m-1-2.south) -- (m-2-3.north);
    \draw[red,thick] (m-1-3.south) -- (m-2-2.north);
    \draw[red,thick] (m-1-3.south) -- (m-2-3.north);
    \draw[red,thick] (m-1-4.south) -- (m-2-1.north);
    \draw[red,thick] (m-1-4.south) -- (m-2-4.north);

    \draw[red,thick] (m-1-6.south) -- (m-2-6.north);
    \draw[red,thick] (m-1-6.south) -- (m-2-9.north);
    \draw[red,thick] (m-1-7.south) -- (m-2-7.north);
    \draw[red,thick] (m-1-7.south) -- (m-2-8.north);
    \draw[red,thick] (m-1-8.south) -- (m-2-7.north);
    \draw[red,thick] (m-1-8.south) -- (m-2-8.north);
    \draw[red,thick] (m-1-9.south) -- (m-2-6.north);
    \draw[red,thick] (m-1-9.south) -- (m-2-9.north);

    \draw[red,thick] (m-2-1.south) -- (mbot-1-3.north);
    \draw[red,thick] (m-2-1.south) -- (mbot-1-2.north);
    \draw[red,thick] (m-2-2.south) -- (mbot-1-1.north);
    \draw[red,thick] (m-2-2.south) -- (mbot-1-4.north);
    \draw[red,thick] (m-2-3.south) -- (mbot-1-3.north);
    \draw[red,thick] (m-2-3.south) -- (mbot-1-2.north);
    \draw[red,thick] (m-2-4.south) -- (mbot-1-1.north);
    \draw[red,thick] (m-2-4.south) -- (mbot-1-4.north);

    \draw[ForestGreen,thick] (m-2-6.south) -- (mbot-1-1.north);
    \draw[ForestGreen,thick] (m-2-6.south) -- (mbot-1-3.north);
    \draw[ForestGreen,thick] (m-2-7.south) -- (mbot-1-2.north);
    \draw[ForestGreen,thick] (m-2-7.south) -- (mbot-1-4.north);
    \draw[ForestGreen,thick] (m-2-8.south) -- (mbot-1-1.north);
    \draw[ForestGreen,thick] (m-2-8.south) -- (mbot-1-3.north);
    \draw[ForestGreen,thick] (m-2-9.south) -- (mbot-1-2.north);
    \draw[ForestGreen,thick] (m-2-9.south) -- (mbot-1-4.north);

  \end{tikzpicture}
  
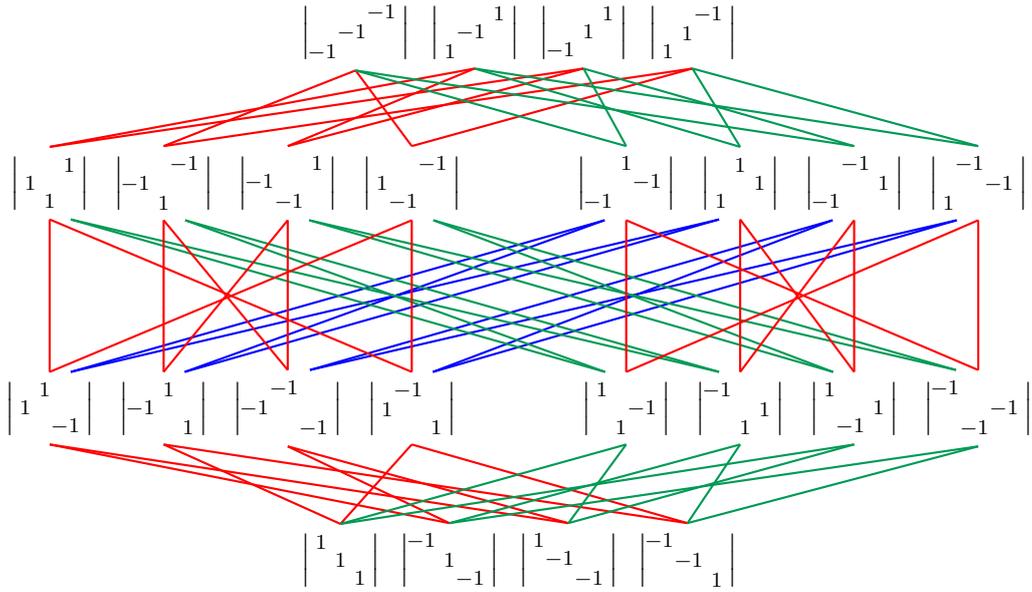
\captionof{figure}{The Bruhat order on the extended Weyl group of $\SL(3,\bR)$. Different colors are used purely for better visibility.}
  \label{fig:SL3_Bruhat_order}
\end{center}

To find balanced ideals, we first list the possible transverse positions $w_0 \in T$. These are
\[ \begin{pmatrix} & & 1 \\ & -1 \\ 1 \end{pmatrix}, \begin{pmatrix} & & -1 \\ & -1 \\ -1 \end{pmatrix}, \begin{pmatrix} & & -1 \\ & 1 \\ 1 \end{pmatrix}, \begin{pmatrix} & & 1 \\ & 1 \\ -1 \end{pmatrix}. \]
The first two of these are conjugate by $\Mbar$, as are the last two. So we have to distinguish two cases.
\begin{enumerate}
\item $w_0 = \begin{psmallmatrix} & & \;1 \\ & -1 \\ 1 \end{psmallmatrix}$. Since $w_0^2 = 1$, the minimal choice of $E$ (or $R$) is the trivial group, so $\widetilde W_{R,S} = \widetilde W_{\{1\},\{1\}} = \widetilde W$. The involution induced by $w_0$ acts on $\widetilde{W}$ in the following way (each dot represents the corresponding matrix from \autoref{fig:SL3_Bruhat_order}):

  \begin{center}
    \begin{tikzpicture}
      \matrix (mtop) [matrix of nodes,ampersand replacement=\&,column sep=0.5cm]{
        { } \& { } \& { } \& { } \\
      };

      \foreach \x in {1,2,3,4} \draw[fill=black] (mtop-1-\x) circle (0.1cm);

      \matrix (m) [matrix of nodes,ampersand replacement=\&,row sep=1cm,yshift=-2cm,column sep=0.5cm]{
        { } \& { } \& { } \& { } \&
        \hspace{1cm} \&
        { } \& { } \& { } \& { } \\
        { } \& { } \& { } \& { } \&
        \hspace{1cm} \&
        { } \& { } \& { } \& { } \\
      };

      \foreach \x in {1,2,3,4,6,7,8,9} {
        \foreach \y in {1,2} \draw[fill=black] (m-\y-\x) circle (0.1cm);
      }

      \matrix (mbot) [matrix of nodes,ampersand replacement=\&,column sep=0.5cm,yshift=-4cm]{
        { } \& { } \& { } \& { } \\
      };

      \foreach \x in {1,2,3,4} \draw[fill=black] (mbot-1-\x) circle (0.1cm);

      \draw[<->,red,thick] (mtop-1-1) -- (mbot-1-2);
      \draw[<->,red,thick] (mtop-1-2) -- (mbot-1-1);
      \draw[<->,red,thick] (mtop-1-3) -- (mbot-1-4);
      \draw[<->,red,thick] (mtop-1-4) -- (mbot-1-3);
      \draw[<->,red,thick] (m-1-1) -- (m-2-2);
      \draw[<->,red,thick] (m-1-2) -- (m-2-1);
      \draw[<->,red,thick] (m-1-3) -- (m-2-4);
      \draw[<->,red,thick] (m-1-4) -- (m-2-3);
      \draw[<->,red,thick] (m-1-6) -- (m-2-7);
      \draw[<->,red,thick] (m-1-7) -- (m-2-6);
      \draw[<->,red,thick] (m-1-8) -- (m-2-9);
      \draw[<->,red,thick] (m-1-9) -- (m-2-8);
    \end{tikzpicture}
  \end{center}

  We obtain the following balanced ideals:
  \begin{enumerate}
  \item The lift of the unoriented balanced ideal contains all relative positions in the bottom half of the picture.
  \item Ideals containing two positions from the third level and everything below these two positions in the Bruhat order. The possible pairs of positions from the third level that can be chosen are (ideals which are equivalent by right-multiplication by $\Mbar$ are in curly brackets): $\{ (1,4), (2,3) \}, \{ (5,8), (6,7)\} $ as well as $\{ (1,7), (2,8), (3,6), (4,5)\}, \{(1,8), (2,7), (3,5), (4,6)\}$. In the following picture, we drew the examples $(1,4)$ in red and $(1,7)$ in green.

    \begin{center}
      \begin{tikzpicture}
        \matrix (mtop) [matrix of nodes,ampersand replacement=\&,column sep=0.5cm]{
          { } \& { } \& { } \& { } \\
        };

        \foreach \x in {1,2,3,4} \draw[fill=black] (mtop-1-\x) circle (0.1cm);

        \matrix (m) [matrix of nodes,ampersand replacement=\&,row sep=1cm,yshift=-2cm,column sep=0.5cm]{
          { } \& { } \& { } \& { } \&
          \hspace{1cm} \&
          { } \& { } \& { } \& { } \\
          { } \& { } \& { } \& { } \&
          \hspace{1cm} \&
          { } \& { } \& { } \& { } \\
        };

        \foreach \x in {1,2,3,4,6,7,8,9} {
          \foreach \y in {1,2} \draw[fill=black] (m-\y-\x) circle (0.1cm);
        }

        \matrix (mbot) [matrix of nodes,ampersand replacement=\&,column sep=0.5cm,yshift=-4cm]{
          { } \& { } \& { } \& { } \\
        };

        \foreach \x in {1,2,3,4} \draw[fill=black] (mbot-1-\x) circle (0.1cm);

        \draw[ForestGreen,rounded corners] ([shift={(-0.4,-0.4)}]mbot-1-1.south west) -- ([shift={(-0.4,-0.4)}]m-2-4.south west) -- ([shift={(-0.4,-0.4)}]m-2-1.south west) -- ([shift={(-0.4,0.4)}]m-1-1.north west) -- ([shift={(0.4,0.4)}]m-1-1.north east) -- ([shift={(0.4,-0.2)}]m-2-1.south east) -- ([shift={(-0.25,-0.2)}]m-2-3.south west) -- ([shift={(-0.25,0.4)}]m-2-3.north west) -- ([shift={(-0.4,0.4)}]m-2-8.north west) -- ([shift={(-0.4,0.4)}]m-1-8.north west) -- ([shift={(0.4,0.4)}]m-1-8.north east) -- ([shift={(0.4,-0.4)}]m-2-8.south east) -- ([shift={(0.4,-0.4)}]m-2-6.south east) -- ([shift={(0.4,-0.4)}]mbot-1-4.south east) -- cycle;

        \draw[red,rounded corners] ([shift={(-0.3,-0.3)}]mbot-1-1.south west) -- ([shift={(-0.3,-0.3)}]m-2-4.south west) -- ([shift={(-0.3,-0.3)}]m-2-1.south west) -- ([shift={(-0.3,0.3)}]m-1-1.north west) -- ([shift={(0.3,0.3)}]m-1-1.north east) -- ([shift={(0.3,-0.1)}]m-2-1.south east) -- ([shift={(-0.3,-0.1)}]m-2-4.south west) -- ([shift={(-0.3,0.3)}]m-1-4.north west) -- ([shift={(0.3,0.3)}]m-1-4.north east) -- ([shift={(0.3,0.3)}]m-2-4.north east) -- ([shift={(0.3,0.3)}]m-2-9.north east) -- ([shift={(0.3,-0.3)}]m-2-9.south east) -- ([shift={(0.3,-0.3)}]m-2-6.south east) -- ([shift={(0.3,-0.3)}]mbot-1-4.south east) -- cycle;
      \end{tikzpicture}
    \end{center}

  \item Ideals containing one relative position from the third level and everything except its $w_0$-image from the second level and below. There are 8 balanced ideals of this type, determined by the element on the third level. Right multiplication by $\Mbar$ identfies the first 4 and the last 4 ideals.
  \end{enumerate}

  In total, we find 21 $w_0$--balanced ideals, which form 7 equivalence classes with respect to right--multiplication with $\Mbar$. Let us emphasize again that the balanced ideals in (b) and (c) are not lifts of balanced ideals from the unoriented setting. If a representation satisfies the prerequisites of \autoref{thm:prop_disc_cocpt}, we therefore obtain new cocompact domains of discontinuity in oriented flag manifolds. For example, we will apply this to Hitchin representations in \autoref{sec:Hitchin}.

\item $w_0 = \begin{psmallmatrix} & & -1 \\ & 1 \\ 1 \end{psmallmatrix}$. Since $w_0^2 = \begin{psmallmatrix} -1 \\ & 1 \\ && -1 \end{psmallmatrix}$ is nontrivial, the minimal choice of $R$ is $\{1,w_0^2\}$ and we consider $R\backslash \widetilde{W}$. By \autoref{lem:order_projection} we get the Bruhat order on $R \backslash \widetilde W$ as the projection of \autoref{fig:SL3_Bruhat_order}. It is shown in \autoref{fig:SL3_mod_w02} alongside the action of $w_0$ on $R \backslash \widetilde W$.
  \begin{center}
    \begin{tikzpicture}
      \matrix (mtop) [matrix of math nodes,ampersand replacement=\&,row sep=1cm]{
        \begin{smatrix} && \!\!-1 \\ &  \!\!-1 \\    \!\!-1 \end{smatrix} \&
        \begin{smatrix} && 1      \\ &  1      \\    \!\!-1 \end{smatrix} \\
      };

      \matrix (m) [matrix of math nodes,ampersand replacement=\&,row sep=0.8cm,yshift=-2.5cm]{
        \begin{smatrix} && 1      \\    1      \\ &  1      \end{smatrix} \&
        \begin{smatrix} && \!\!-1 \\    \!\!-1 \\ &  1      \end{smatrix} \&
        \hspace{1cm} \&
        \begin{smatrix} &  1      \\ && \!\!-1 \\    \!\!-1 \end{smatrix} \&
        \begin{smatrix} &  1      \\ && 1      \\    1      \end{smatrix} \\
        \begin{smatrix} &  1      \\    1      \\ && \!\!-1 \end{smatrix} \&
        \begin{smatrix} &  1      \\    \!\!-1 \\ && 1      \end{smatrix} \&
        \hspace{1cm} \&
        \begin{smatrix}    1      \\ && \!\!-1 \\ &  1      \end{smatrix} \&
        \begin{smatrix}    \!\!-1 \\ && 1      \\ &  1      \end{smatrix} \\
      };

      \matrix (mbot) [matrix of math nodes,ampersand replacement=\&,row sep=1cm,yshift=-5cm]{
        \begin{smatrix}    1      \\ &  1      \\ && 1      \end{smatrix} \&
        \begin{smatrix}    1      \\ &  \!\!-1 \\ && \!\!-1 \end{smatrix} \\
      };

      \draw[red,thick] (mtop-1-1.south) -- (m-1-1.north);
      \draw[red,thick] (mtop-1-1.south) -- (m-1-2.north);
      \draw[red,thick] (mtop-1-2.south) -- (m-1-1.north);
      \draw[red,thick] (mtop-1-2.south) -- (m-1-2.north);

      \draw[ForestGreen,thick] (mtop-1-1.south) -- (m-1-4.north);
      \draw[ForestGreen,thick] (mtop-1-1.south) -- (m-1-5.north);
      \draw[ForestGreen,thick] (mtop-1-2.south) -- (m-1-4.north);
      \draw[ForestGreen,thick] (mtop-1-2.south) -- (m-1-5.north);

      \draw[red,thick] (m-1-1.south) -- (m-2-1.north);
      \draw[red,thick] (m-1-2.south) -- (m-2-2.north);

      \draw[red,thick] (m-1-4.south) -- (m-2-4.north);
      \draw[red,thick] (m-1-5.south) -- (m-2-5.north);

      \draw[ForestGreen,thick] (m-1-1.south east) -- (m-2-4.north west);
      \draw[ForestGreen,thick] (m-1-1.south east) -- (m-2-5.north west);
      \draw[ForestGreen,thick] (m-1-2.south east) -- (m-2-4.north west);
      \draw[ForestGreen,thick] (m-1-2.south east) -- (m-2-5.north west);

      \draw[blue,thick] (m-1-4.south west) -- (m-2-1.north east);
      \draw[blue,thick] (m-1-4.south west) -- (m-2-2.north east);
      \draw[blue,thick] (m-1-5.south west) -- (m-2-1.north east);
      \draw[blue,thick] (m-1-5.south west) -- (m-2-2.north east);

      \draw[red,thick] (m-2-1.south) -- (mbot-1-1.north);
      \draw[red,thick] (m-2-1.south) -- (mbot-1-2.north);
      \draw[red,thick] (m-2-2.south) -- (mbot-1-1.north);
      \draw[red,thick] (m-2-2.south) -- (mbot-1-2.north);

      \draw[ForestGreen,thick] (m-2-4.south) -- (mbot-1-1.north);
      \draw[ForestGreen,thick] (m-2-4.south) -- (mbot-1-2.north);
      \draw[ForestGreen,thick] (m-2-5.south) -- (mbot-1-1.north);
      \draw[ForestGreen,thick] (m-2-5.south) -- (mbot-1-2.north);
    \end{tikzpicture}
    \begin{tikzpicture}
      \matrix (mtop) [matrix of nodes,ampersand replacement=\&,column sep=0.5cm]{
        { } \& { } \\
      };

      \foreach \y in {1,2} \draw[fill=black] (mtop-1-\y) circle (0.1cm);

      \matrix (m) [matrix of nodes,ampersand replacement=\&,row sep=1cm,yshift=-2cm,column sep=0.5cm]{
        { } \& { } \& \hspace{1cm} \& { } \& { } \\
        { } \& { } \& \hspace{1cm} \& { } \& { } \\
      };

      \foreach \x in {1,2,4,5} {
        \foreach \y in {1,2} \draw[fill=black] (m-\y-\x) circle (0.1cm);
      }

      \matrix (mbot) [matrix of nodes,ampersand replacement=\&,column sep=0.5cm,yshift=-4cm]{
        { } \& { } \\
      };

      \foreach \y in {1,2} \draw[fill=black] (mbot-1-\y) circle (0.1cm);

      \draw[<->,red,thick] (mtop-1-1) -- (mbot-1-2);
      \draw[<->,red,thick] (mtop-1-2) -- (mbot-1-1);
      \draw[<->,red,thick] (m-1-1) -- (m-2-1);
      \draw[<->,red,thick] (m-1-2) -- (m-2-2);
      \draw[<->,red,thick] (m-1-4) -- (m-2-4);
      \draw[<->,red,thick] (m-1-5) -- (m-2-5);

      \path (0,1.1) --  (0, -5.1); 
    \end{tikzpicture}

    
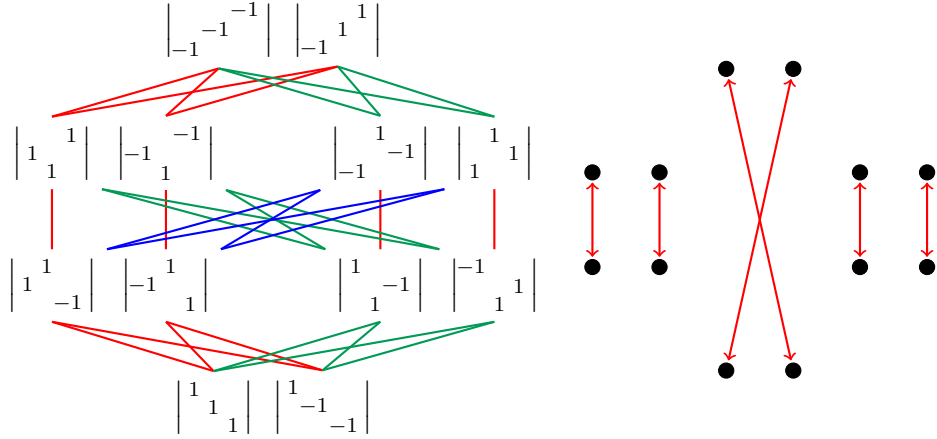
\captionof{figure}{The Bruhat order on $R \backslash \widetilde W$ and the involution given by $w_0$. Different colors are used purely for better visibility.}
    \label{fig:SL3_mod_w02}

  \end{center}

  In this case, the only balanced ideal is the lift of the unoriented one, and we do not obtain any new cocompact domains of discontinuity in oriented flag manifolds.
\end{enumerate}
\subsection{\texorpdfstring{Geometric interpretation of relative positions for $\SL(3,\bR)$}{Geometric interpretation of relative positions for SL(3,R)}}	\label{sec:geometric_interpretation}
In order to give a hands--on description of the various relative positions we saw in the previous subsection, we need notions of direct sums and intersections that take orientations into account. These notions appeared already in \cite{Guichard}, where they are used to describe curves of flags. First of all, let us fix some notation for oriented subspaces.

\begin{Def}
  Let $A,B \subset \bR^n$ be oriented subspaces. Then we denote by $-A$ the subspace $A$ with the opposite orientation. If $A$ and $B$ agree as oriented subspaces, we write $A \equalo B$.
\end{Def}
Now let $A,B \subset \bR^n$ be oriented subspaces. If they are (unoriented) transverse, taking a positive basis of $A$ and extending it by a positive basis of $B$ yields a basis of $A\oplus B$. Declaring this basis to be positive defines an orientation on $A \oplus B$.
The orientation on the direct sum depends on the order we write the two subspaces in,
\[ A\oplus B \equalo (-1)^{\dim(A)\dim(B)}B\oplus A. \]
The case of intersections is slightly more difficult. Assume that $A,B \subset \bR^n$ are oriented subspaces such that $A + B = \bR^n$, and fix a standard orientation on $\bR^n$. Let $A' \subset A$ be a subspace complementary to $A\cap B$ and analogously $B'\subset B$ a subspace complementary to $A\cap B$. We fix orientations on these two subspaces by requiring that
\[ A' \oplus B \equalo \bR^n \]
and
\[ A \oplus B' \equalo \bR^n.\]
Then there is a unique orientation on $A\cap B$ satisfying
\[ A' \oplus (A\cap B) \oplus B' \equalo \bR^n. \]
This is the induced orientation on the intersection. Since the set of subspaces of $A$ complementary to $A\cap B$ can be identified with $\Hom(A',A\cap B)$ and is therefore (simply) connected, the result does not depend on the choice of $A'$, and analogously does not depend on the choice of $B'$. Like the oriented sum, it depends on the order we write the two subspaces in,
\[ B \cap A \equalo (-1)^{\codim(A)\codim(B)}A \cap B. \]
With this terminology at hand, consider the oriented relative positions shown in \autoref{fig:SL3_Bruhat_order}. Let $f\in G/B_0$ be a reference complete oriented flag. We denote by $f^{(k)}$ the $k$--dimensional part of the flag $f$. Let $w = \begin{psmallmatrix} && -1\\ & -1\\ -1 \end{psmallmatrix} \in \widetilde{W}$. Then the refined Schubert stratum of flags at position $w$ with respect to $f$
\[ C_w(f) = \{ F \in G/B_0 \mid f^{(2)} \oplus F^{(1)} \equalo -\bR^3, \ f^{(1)} \oplus F^{(2)} \equalo -\bR^3 \}. \]
The other three positions of the highest level are characterized by the other choices of the two signs. Similarly, for the position $w' = \begin{psmallmatrix} && 1\\ 1\\ & 1 \end{psmallmatrix} \in \widetilde{W}$, we obtain
\[ C_{w'}(f) = \{ F \in G/B_0 \mid f^{(1)} \oplus F^{(1)} \equalo f^{(2)}, \ f^{(2)} \cap F^{(2)} \equalo F^{(1)} \}, \]
and the other three positions are characterized by the other choices of the two signs. Similar descriptions hold for the remaining oriented relative positions.
\subsection{A simple example in odd dimension: Halfspaces in spheres}	\label{sec:ideal_spheres}
As the previous subsection demonstrated, calculating the most general relative positions and the Bruhat order gets out of hand very quickly as one increases the dimension. For example, in $\SL(5,\bR)$, there are $120$ unoriented relative positions between complete flags and $1920$ oriented relative positions between complete oriented flags. For practical reasons, it thus makes sense to restrict to more special cases, i.e. to consider relative positions $\widetilde{W}_{R,S}$ for bigger $R,S$ than strictly necessary.

Let $G=\SL(2n+1,\bR)$ and $\theta, \eta \subset \Delta$ be the complements of $\ctheta = \{\alpha_{n},\alpha_{n+1}\}$, $\ceta = \{\alpha_1\}$, so that $\F_\theta$ is the space of partial flags consisting of the dimension $n$ and $n+1$ parts, and $\F_\eta$ is $\RP^{2n}$. Furthermore, let
\[E = \langle \Mbar_\theta, \I(\alpha_n)^2\I(\alpha_{n+1})^2 \rangle = \{m \in \Mbar \mid m_{n+1,n+1} = +1\},\]
$F=\Mbar_\eta$, $R = \langle \I(\theta),E \rangle$, and $S = \langle \I(\eta),F \rangle$. Then $\F_R$ is the space of oriented partial flags consisting of oriented $n$-- and $(n+1)$--dimensional subspaces up to changing both orientations simultaneously, and $\F_S$ is $S^{2n}$, the space of oriented lines on $\bR^{2n+1}$. Choose $w_0$ to be antidiagonal with $-1$ as the middle entry. The remaining entries are irrelevant for this example; if $2n+1 \in 4\bZ+3$, there should be an odd number of minus signs, if $2n+1 \in 4\bZ+1$, it has to be even.

\begin{figure}[ht]
  \centering
  \begin{tikzpicture}
    \matrix (m) [matrix of math nodes,ampersand replacement=\&,row sep=0.5cm]{
      \& \begin{smatrix} 0 \\ \svdots \\ 0 \\ 1 \end{smatrix} \\
      \begin{smatrix} 0 \\ \svdots  \\ 1  \\ \svdots \\ 0 \end{smatrix}
      \& \hspace{1cm} \& \begin{smatrix} 0 \\ \svdots  \\ -1 \\  \svdots \\ 0 \end{smatrix} \\
      \& \begin{smatrix} 1 \\ 0 \\ \svdots \\ 0 \end{smatrix} \\
    };

    \draw[thick] (m-1-2.south) -- (m-2-1.north);
    \draw[thick] (m-1-2.south) -- (m-2-3.north);
    \draw[thick] (m-2-1.south) -- (m-3-2.north);
    \draw[thick] (m-2-3.south) -- (m-3-2.north);

    \node[left=-2pt of m-2-1](mid1){\tiny{$n+1$}};

    \draw[<->,red,thick] ([yshift=-4pt]m-1-2.south) -- ([yshift=4pt]m-3-2.north);
    \draw[<->,red,thick] (m-2-1.east) -- (m-2-3.west);
  \end{tikzpicture}
  \caption{Oriented relative positions between $\F_R$ and $\F_S$}
  \label{fig:Bruhat_order_sphere}
\end{figure}
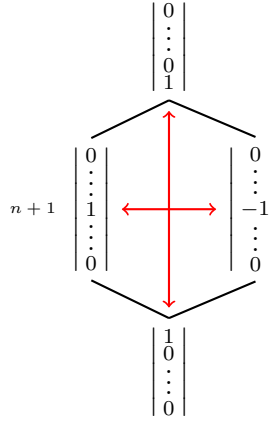

The Bruhat order on the space $\widetilde{W}_{R,S}$ of relative positions as well as the involution $w_0$ are shown in \autoref{fig:Bruhat_order_sphere}. We only need to keep track of the first column of the matrix representative since we right quotient by $S$. The left quotient by $R$ then reduces the possible relative positions further.
There are thus two balanced ideals, determined by choosing one of the two middle positions. This is in contrast to the unoriented case, where the two middle positions coincide and are a fixed point of the involution. In the case $n = 3$, the lifts of these two balanced ideals to $\widetilde W$ are among the 21 $w_0$--balanced ideals described in \autoref{sec:SL3R}. In the notation of case (i)b there, these are the ideals $(1,4)$ and $(2,3)$. In particular, the balanced ideal indicated in red in \autoref{sec:SL3R} is the lift of the balanced ideal we obtain here in $\widetilde W_{R,S}$ by choosing the first of the middle positions.

The geometric description of the relative positions is as follows. Let $f\in\F_R$ be a reference flag. Then $\dc{w} = \begin{bsmallmatrix} 0 \\ \svdots \\ 0 \\ 1 \end{bsmallmatrix}$ and $\dc{w'} = \begin{bsmallmatrix} 1 \\ 0 \\ \svdots \\ 0 \end{bsmallmatrix}$ correspond to
\[ C_{\dc{w}}(f) = \{ F\in S^{2n} \mid F \not\in f^{(n+1)} \}, \qquad C_{\dc{w'}}(f) = \{ F\in S^{2n} \mid F \in f^{(n)} \},\]
and $\dc{w''} = \begin{bsmallmatrix} 0 \\ \svdots  \\ 1  \\ \svdots \\ 0 \end{bsmallmatrix}$ corresponds to
\[ C_{\dc{w''}}(f) = \{ F\in S^{2n} \mid f^{(n)} \oplus F \equalo f^{(n+1)} \}. \]
This can be rephrased slightly: The codimension 1 subspace $f^{(n)} \subset f^{(n+1)}$ decomposes $f^{(n+1)}$ into two half--spaces, and we say that an oriented line $l \subset f^{(n+1)}$ is in the positive half--space if $f^{(n)} \oplus l \equalo f^{(n+1)}$. This is invariant under simultaneously changing the orientations of both $f^{(n)}$ and $f^{(n+1)}$ and therefore well--defined. Then $C_{\dc{w'}}(f)$ is the spherical projectivization of the positive half of $f^{(n+1)}$.

The half great circles in \autoref{fig:sphere} in the introduction are an example of this construction, and the associated cocompact domain of discontinuity in $S^2$ for a convex cocompact representation $\rho \colon F_k \to \SO_0(2,1)$ corresponds to this balanced ideal. The higher--dimensional cases apply to Hitchin representations and generalized Schottky representations into $\PSL(4n+3,\bR)$, yielding cocompact domains of discontinuity in $S^{4n+2}$ (see \autoref{sec:Hitchin} and \autoref{sec:example_gen_Schottky}). The latter case includes in particular the motivating example from the introduction.
\section{Examples of representations admitting cocompact domains of discontinuity} \label{sec:examples_of_representations}
\subsection{Hitchin representations}	\label{sec:Hitchin}
Let $\Gamma = \pi_1(S)$ be the fundamental group of a closed surface $S$ of genus at least 2. Particularly simple examples of representations of $\Gamma$ into $\PSL(n,\bR)$ are the Fuchsian ones: These are of the form $\iota \circ \rho_0$, where $\rho_0 \colon \Gamma \to \PSL(2,\bR)$ is injective with discrete image and $\iota \colon \PSL(2,\bR) \to \PSL(n,\bR)$ is the irreducible representation. Connected components of $\Hom(\Gamma,\PSL(n,\bR))$ which contain Fuchsian representations are called Hitchin components and their elements Hitchin representations.

If $n$ is odd, then the space $\Hom(\Gamma, \PSL(n,\bR))$ has 3 connected components, one of which is Hitchin. If $n$ is even, there are 6 components in total, and 2 of them are Hitchin components \cite{Hitchin}. Every Hitchin representation is $B$--Anosov \cite{Labourie}.

To find out if the limit map of a Hitchin representation lifts to an oriented flag manifold, let us first take a closer look at the irreducible representation. The standard Euclidean scalar product on $\bR^2$ induces a scalar product on the symmetric product $\Sym^{n-1}\bR^2$ by restricting the induced scalar product on the tensor power to symmetric tensors. Let $X = \binom{1}{0}$ and $Y = \binom{0}{1}$ be the standard orthonormal basis of $\bR^2$. Then $e_i = \sqrt{\binom{n-1}{i-1}}X^{n-i}Y^{i-1}$ for $1 \leq i \leq n$ is an orthonormal basis of $\Sym^{n-1} \bR^2$ and provides an identification $\bR^n \cong \Sym^{n-1}\bR^2$. For $A \in \SL(2,\bR)$ let $\iota(A) \in \SL(n,\bR)$ be the induced action on $\Sym^{n-1}\bR^2$ in this basis. The homomorphism
\[\iota \colon \SL(2,\bR) \to \SL(n,\bR),\]
defined this way is the (up to conjugation) unique irreducible representation. It maps $-1$ to $(-1)^{n-1}$ and is therefore also well--defined as a map $\iota \colon \PSL(2,\bR) \to \PSL(n,\bR)$. The induced action on $\Sym^{n-1} \bR^2$ preserves the scalar product described above, so $\iota(\PSO(2)) \subset \PSO(n)$.

It is easy to see that $\iota$ maps diagonal matrices to diagonal matrices. It also maps upper triangular matrices into $B_0 \subset \PSL(n,\bR)$ (that is, upper triangular matrices with the diagonal entries either all positive or all negative). Therefore, $\iota$ induces a smooth equivariant map
\[ \varphi \colon \RP^1 \to \F_{\{1\}}\]
between the complete oriented flag manifolds of $\PSL(2,\bR)$ and $\PSL(n,\bR)$.

\begin{Prop}	\label{prop:Hitchin_Anosov}
  Let $\rho \colon \Gamma \to \PSL(n,\bR)$ be a Hitchin representation. Then its limit map $\xi \colon \bdry \to \F_{\Mbar} = G/B$ lifts to the fully oriented flag manifold $\F_{\{1\}} = G / B_0$ with transversality type
  \[w_0 = \begin{pmatrix} &&& \iddots \\ && 1 & \\ & -1 && \\ 1 &&& \end{pmatrix}. \]
  So all Hitchin representations are $B_0$--Anosov.
\end{Prop}
\begin{Prf}
  Since the $B_0$--Anosov representations are a union of connected components of $B$--Anosov representations by \autoref{prop:clopen}, we can assume that $\rho = \iota \circ \rho_0$ is Fuchsian.

  Let $\xi_0 \colon \bdry \to \RP^1$ be the limit map of $\rho_0$ and $\pi \colon \F_{\{1\}} \to \F_\Mbar$ the projection forgetting all orientations. Then the limit map $\xi \colon \bdry \to \F_\Mbar$ of $\rho$ is just $\pi \circ \varphi \circ \xi_0$ (this is the unique continuous and dynamics--preserving map, see \cite[Remark 2.32b]{GueritaudGuichardKasselWienhard}). So $\widehat \xi = \varphi \circ \xi_0$ is a continuous and $\rho$--equivariant lift to $\F_{\{1\}}$. To calculate the transversality type, let $x, y \in \bdry$ with $\xi_0(x) = [1]$ and $\xi_0(y) = [w] \in \RP^1$, where $w \in \PSL(2,\bR)$ is the anti--diagonal matrix with $\pm 1$ entries.
  Then, since $\iota(w) = w_0$,
  \[\pos(\widehat\xi(x),\widehat\xi(y)) = \pos(\varphi([1]),\varphi([w])) = \pos([1],[\iota(w)]) = w_0 \in \widetilde W. \qedhere\]
\end{Prf}
\begin{Rem}
  Note that Hitchin representations map into $\PSL(n,\bR)$ and not $\SL(n,\bR)$. If $n$ is even, the fully oriented flag manifold $\F_{\{1\}}$ in $\PSL(n,\bR)$ is the space of flags $f^{(1)} \subset \dots \subset f^{(n-1)}$ with a choice of orientation on every part, but up to simultaneously reversing the orientation in every odd dimension (the action of $-1$). While we could lift $\rho$ to $\SL(n,\bR)$, its limit map would still only lift to $\F^{\SL(n,\bR)}_{\{\pm 1\}} = \F^{\PSL(n,\bR)}_{\{1\}}$ and not give us any extra information.
\end{Rem}
Now that we know that Hitchin representations are $B_0$--Anosov, we can apply \autoref{thm:prop_disc} and \autoref{thm:cocpt}. For every $w_0$--balanced ideal in $\widetilde W$ we get a cocompact domain of discontinuity in the oriented flag manifold $\F_{\{1\}}$ of $\PSL(n,\bR)$. These include lifts of the domains in unoriented flag manifolds constructed in \cite{KapovichLeebPortiFlagManifolds}, but also some new examples.

There are 21 different such $w_0$--balanced ideals if $n=3$ (see \autoref{sec:SL3R}) and already 4732 of them if $n=4$, which makes it infeasible to list all of them here. However, it is not difficult to state in which oriented Grassmannians these domains exist.
\begin{Prop}\label{prop:grassmannians}
  Let $n \geq 3$ and $\rho \colon \Gamma \to \PSL(n,\bR)$ be a Hitchin representation. Assume that either
  \begin{enumerate}
  \item $n$ is even and $k$ is odd, or
  \item $n$ is odd and $k(n+k+2)/2$ is odd.
  \end{enumerate}
  Then there exists a nonempty, open $\Gamma$--invariant subset $\Omega \subset \Gro(k, n)$ of the Grassmannian of oriented $k$--subspaces of $\bR^n$, such that the action of $\Gamma$ on $\Omega$ is properly discontinuous and cocompact.
\end{Prop}
\begin{Rems}\ \reallynopagebreak
  \begin{enumerate}
  \item The domain $\Omega$ is not unique, unless $n$ is even and $k \in \{1, n-1\}$.
  \item In case (i) of \autoref{prop:grassmannians}, there is a cocompact domain of discontinuity also in the unoriented Grassmannian, and $\Omega$ is just the lift of one of these. The domains in case (ii) are new (see \cite{Stecker}).
  \end{enumerate}
\end{Rems}
\begin{Prf}
  In the light of \autoref{thm:prop_disc_cocpt} it suffices to show that there is a $w_0$--balanced ideal in the set $\widetilde W / S$ where $S = \langle \I(\Delta \setminus \{\alpha_k\}) \rangle$. A $w_0$--balanced ideal exists if and only if the action of $w_0$ on $\widetilde W / S$ has no fixed points (see \autoref{lem:involution_order_reversing} and \autoref{lem:ideal_existence}).

  To see that $w_0$ has no fixed points on $\widetilde W / S$, observe that every equivalence class in $\widetilde W / S$ has a representative whose first $k$ columns are either the standard basis vectors $e_{i_1}, \dots e_{i_k}$ or $-e_{i_1}, e_{i_2}, \dots, e_{i_k}$, with $1 \leq i_1 < \dots < i_k \leq n$. So we can identify $\widetilde W / S$ with the set
  \[ \{\pm 1\} \times \{\text{$k$--element subsets of $\{1,\dots,n\}$}\}.\]
  The action of $w_0$ on this is given by
  \[(\varepsilon, \{i_1, \dots, i_k\}) \quad \mapsto \quad ((-1)^{k(k-1)/2 + \sum_j (i_j + 1)}\varepsilon, \{n+1-i_k, \dots, n+1-i_1\}).\]
  Only looking at the second factor, this can have no fixed points if $n$ is even and $k$ is odd, showing case (i). Otherwise, to get a fixed point it is necessary that $i_j + i_{k+1-j} = n+1$ for all $j \leq k$. But then
  \[\frac{k(k-1)}{2} + \sum_{j=1}^k (i_j + 1) =  \frac{k(k-1)}{2} + \frac{k}{2}(n+3) = \frac{k(n+k+2)}{2},\]
  so $w_0$ fixes these elements if and only if $k(n+k+2)/2$ is even. Note that this number is always even if $n$ and $k$ are both even, which is why assuming $n$ odd in case (ii) does not weaken the statement.

  It remains to show that every $\Omega \in \Gro(k,n)$ constructed from a balanced ideal $I \subset \widetilde W / S$ is nonempty. Consider the lifts  $\Omega' \subset \F_{\{1\}}$ of $\Omega$ and  $I' \subset \widetilde W$ of $I$. Then $\Omega'$ is the domain in $\F_{\{1\}}$ given by $I'$. We will show that $\F_{\{1\}} \setminus \Omega'$ has covering codimension\footnote{The (Lebesgue) covering dimension of a topological space $X$ is the smallest number $n$ such that every open cover of $X$ admits a refinement with the property that each point of $X$ is contained in at most $n+1$ of its elements.} at least $1$, so $\Omega'$ must be nonempty. Similar arguments were used in \cite{GuichardWienhardDomains} to prove the nonemptiness of certain domains. See \cite{Nagata} for more background on dimension theory. In this case, we could also use the dimension of $\F_{\{1\}} \setminus \Omega'$ as a CW-complex, but the present approach has the benefit of generalizing to word hyperbolic groups with more complicated boundaries.

  By \autoref{lem:properties_of_K} and the proof of \autoref{thm:cocpt}, $\F_{\{1\}} \setminus \Omega'$ is homeomorphic to a fiber bundle over $\bdry \cong S^1$ with fiber $\bigcup_{\dc{w}\in I} C_{\dc{w}}([1])$. The covering dimension is invariant under homeomorphisms and has the following locality property: If a metric space is decomposed into open sets of dimension (at most) $k$, then the whole space is (at most) $k$--dimensional\footnote{This follows from the equivalence of covering dimension and strong inductive dimension (\cite[Theorem II.7]{Nagata}) and locality of the strong inductive dimension.}. Therefore, the dimension of this fiber bundle equals the dimension of a local trivialization, that is, the dimension of the product $\bR \times \bigcup_{\dc{w}\in I} C_{\dc{w}}([1])$. By \cite[Theorem 2]{Morita}, the dimension of a product is the sum of the dimensions of the factors whenever one of the factors is a CW complex\footnote{In that paper, Kat\v etov--Smirnov covering dimension is used, which coincides with (Lebesgue) covering dimension for normal spaces.}. Thus
  \[\dim (\F_{\{1\}} \setminus \Omega') = 1 + \max_{w\in I'} \dim C_w([1]) = 1 + \max_{w \in I'} \ell(w).\]
  If we know that $\ell(w) \leq \ell(w_0) - 2$ for every $w \in I'$, then, since $\dim \F_{\{1\}} = \ell(w_0)$, the codimension of $\F_{\{1\}} \setminus \Omega'$ is at least 1, so $\Omega \neq \varnothing$.

  For $k < n$, if we write $w_k = \I(\alpha_1) \I(\alpha_2) \cdots \I(\alpha_k)$, and $\widetilde{w_k} = \I(\alpha_2) \cdots \I(\alpha_k)$, then by direct calculation, one verifies that
  \[ w_0 = w_{n-1} \cdots w_1 \]
  and
  \[ w_0\I(\alpha_k) = w_{n-1} \cdots w_{k+1} \widetilde {w_k} w_{k-1} \cdots w_1 \]
  (see \autoref{sec:SLnR} for an explicit description of $\widetilde W$). These are reduced expressions in the $\I(\alpha_i)$. So if $n \geq 3$, then $\I(\alpha_k) \leq w_0 \I(\alpha_k)$ by \autoref{prop:Bruhat_order_quotient}. Therefore, the balanced ideal $I' \subset \widetilde W$ cannot contain $w_0 \I(\alpha_k)$ and thus no element of length $\ell(w_0)-1$.
\end{Prf}
A special case of such cocompact domains of discontinuity for Hitchin representations $\rho \colon \Gamma \to \PSL(4n+3,\bR)$ is described by the balanced ideals in \autoref{sec:ideal_spheres}: Let $\widehat \xi \colon \bdry \to \F_{\{1\}}$ be the boundary map of $\rho$, with image in complete oriented flags in $\bR^{4n+3}$. Then the domain in $S^{4n+2}$ is obtained by removing the spherical projectivizations of the positive halves of $\widehat \xi (x)^{(2n+2)}, \ x\in\bdry$. Note that in the case of $\PSL(3,\bR)$, the result is not very interesting: Consider the base case $\bdry \xrightarrow{\rho_0} \PSL(2,\bR) \xrightarrow{\iota} \PSL(3,\bR)$, where $\rho_0$ is Fuchsian and $\iota$ is the irreducible representation. Since the limit set of $\rho_0$ is the full circle, the domain simply consists of two disjoint disks, and the quotient is two disjoint copies of the surface $S$ (compare \autoref{fig:sphere}). In higher dimension however, the domain is always connected and dense in $S^{4n+2}$.
\subsection{Generalized Schottky representations}	\label{sec:example_gen_Schottky}
In \cite{BurelleTreibFlags}, generalized Schottky groups in $G=\PSL(n,\bR)$ are introduced. The construction relies on the existence of a \emph{partial cyclic order} on the space $G/B_0 = \F_{\{1\}}$, which is an oriented version of Fock--Goncharov triple positivity \cite{FockGoncharov} and Labourie's 3--hyperconvexity \cite{Labourie}. We will give a brief overview of generalized Schottky groups and their properties before showing how they fit into our framework. More details and proofs can be found in \cite{BurelleTreibSchottky} and \cite{BurelleTreibFlags}.

For all odd dimensions, the partial cyclic order on $\F_{\{1\}}$ is given as follows. Recall that we defined oriented direct sums in \autoref{sec:geometric_interpretation}.
\begin{Def}
  Let $G=\PSL(2n+1,\bR)$ and $F_1,F_2,F_3 \in \F_{\{1\}}$ be complete oriented flags. Then the triple $(F_1,F_2,F_3)$ is \emph{positive} or \emph{increasing} if
  \[ F_1^{(i_1)} \oplus F_2^{(i_2)} \oplus F_3^{(i_3)} \equalo \bR^{2n+1} \]
  for every choice of $i_1,i_2,i_3 \geq 0$ such that $i_1 + i_2 + i_3 = 2n+1$.
\end{Def}
Note that positivity of a triple includes in particular the condition $F_i^{(i_1)} \oplus F_j^{(i_2)} \equalo \bR^{2n+1}$ for $i<j$ and $i_1+i_2=2n+1$, which we call \emph{oriented transversality} of $F_i$ and $F_j$. In terms of relative positions, this means that $\pos(F_i,F_j)=w_0$, where $w_0$ is the transversality type of Hitchin representations (see \autoref{prop:Hitchin_Anosov}). Having this partial cyclic order allows to define intervals in $\F_{\{1\}}$:
\begin{Def}
 Let $(F_1,F_3)$ be an oriented transverse pair of complete oriented flags. Then the \emph{interval} between $F_1$ and $F_3$ is given by
 \[ \ival{F_1}{F_3} = \{ F_2 \mid (F_1,F_2,F_3) \ \text{is increasing} \}. \]
\end{Def}
Consider a \emph{cycle} $(F_1,\ldots,F_{4k})$, that is, a tuple such that $(F_i,F_j,F_k)$ is increasing for any $i<j<k$. This defines the $2k$ intervals
\[ I_i = \ival{F_{2i-1}}{F_{2i}}, \ 1 \leq i \leq 2k. \]
 We say that a transformation $g \in \PSL(2n+1,\bR)$ \emph{pairs} two intervals $I_i$ and $I_j$, $i \neq j$, if
 \[ g(F_{2i-1}) = F_{2j} \  \text{and} \ g(F_{2i}) = F_{2j-1}. \]
Now pick $k$ generators $g_1,\ldots,g_k \in \PSL(2n+1,\bR)$ pairing all of the intervals in some way.
\begin{Def}
A \emph{purely hyperbolic generalized Schottky group in $\PSL(2n+1,\bR)$} is the group $\langle g_1,\ldots,g_k \rangle \subset \PSL(2n+1,\bR)$, where the generators $g_i$ are constructed as above. The associated representation $\rho \colon F_k \to \PSL(2n+1,\bR)$ is called a \emph{purely hyperbolic generalized Schottky representation}.
\end{Def}
In this definition, ``purely hyperbolic'' refers to the fact that all of the endpoints of intervals are distinct. This assumption is important to obtain contraction properties and provides a link to Anosov representations. It is an easy consequence of the Ping--Pong--Lemma that all generalized Schottky representations are faithful.

A map between partially cyclically ordered sets is called \emph{increasing} if it maps every increasing triple to an increasing triple. Moreover, the abstract group $\Gamma = F_k$ is identified with a subgroup of $\PSL(2,\bR)$ by picking a model Schottky group in $\PSL(2,\bR)$ admitting the same combinatorial setup of intervals and generators as in $\PSL(2n+1,\bR)$. This yields a homeomorphism between $\bdry$ and a Cantor set in $\partial \bH^2$. \\
The following result allows us to apply our theory of domains of discontinuity to purely hyperbolic generalized Schottky representations.
\begin{Thm}[\cite{BurelleTreibFlags}]
Let $\rho \colon \Gamma = F_k \to \PSL(2n+1,\bR)$ be a purely hyperbolic generalized Schottky representation. Then $\rho$ is $B_0$--Anosov. Moreover, the boundary map $\widehat \xi \colon \bdry \to \F_{\{1\}}$ is increasing.
\end{Thm}
As observed above, the definition of the partial cyclic order on $\F_{\{1\}}$ implies that the transversality type $w_0$ of $\widehat \xi$ is the same as for Hitchin representations. Consequently, the same balanced ideals can be used to obtain domains of discontinuity. In particular, the balanced ideals from \autoref{sec:ideal_spheres} yield cocompact domains in $S^{4n+2}$ for purely hyperbolic generalized Schottky representations $\rho \colon \Gamma \to \PSL(4n+3,\bR)$. The motivating example from the introduction is a special case of this: Every convex cocompact representation $\rho_0 \colon \Gamma \to \PSL(2,\bR)$ admits a Schottky presentation, and the composition $\Gamma \xrightarrow{\rho_0} \PSL(2,\bR) \xrightarrow{\iota} \PSL(3,\bR)$ with the irreducible representation $\iota$ is a purely hyperbolic generalized Schottky representation acting properly discontinuously and cocompactly on the domain shown in \autoref{fig:sphere}.
%
%
%
%
%
%
%
%

\section{\texorpdfstring{Block embeddings of irreducible representations of $\SL(2,\bR)$}{Block embeddings of irreducible representations of SL(2,R)}}	\label{sec:block_embeddings}

Let $n$ be odd and $k \leq n$. Let $\iota_k \colon \SL(2,\bR) \to \SL(k,\bR)$ be the irreducible representation (see \autoref{sec:Hitchin}). Then we define
\[b_k \colon \SL(2,\bR) \to \SL(n,\bR), \quad A \mapsto \begin{pmatrix}\iota_k(A) & 0 \\ 0 & \iota_{n-k}(A)\end{pmatrix}.\]
Let $\Gamma$ be the fundamental group of a closed surface of genus at least $2$ and $\rho \colon \Gamma \to \PSL(2,\bR)$ a Fuchsian (i.e. discrete and faithful) representation. Let $\overline\rho \colon \Gamma \to \SL(2,\bR)$ be a lift of $\rho$. We get every other lift of $\rho$ as $\overline\rho^\varepsilon$, where $\varepsilon \colon \Gamma \to \bZ/2\bZ$ is a group homomorphism and $\overline\rho^\varepsilon(\gamma) = (-1)^{\varepsilon(\gamma)} \overline\rho(\gamma)$.

In this section, we will consider representations $\rho_k^\varepsilon = b_k \circ \overline\rho^\varepsilon$ obtained by composing a Fuchsian representation with $b_k$. Our main result is the following proposition and its corollary: For different choices of $k$ and $\varepsilon$, the representations $\rho_k^\varepsilon$ land in different connected components of Anosov representations.

\begin{Prop}	\label{prop:block_embedding_Anosov}
  The representation $\rho_k^\varepsilon$ is $B$--Anosov and there exists $w_k \in T$ such that its limit map lifts to a continuous, equivariant map into $\F_{\{1,w_k^2\}}$ of transversality type $\dc{w_k}$. Thus $\rho_k^\varepsilon$ is $P_{\{1,w_k^2\}}$--Anosov. The choice $R = \{1,w_k^2\}$ is minimal in the sense of \autoref{prop:minimal_type}.

  Futhermore, $w_k$ is up to conjugation by elements of $\Mbar$ given by
  \begin{align*}
    w_k = \begin{pmatrix} & & & J \\ & & \delta \\ & K \\ L \end{pmatrix}
  \end{align*}
  with
  \[ \delta = \begin{cases} (-1)^{(k-1)/2}, & \text{$k$ odd} \\ (-1)^{(n-k-1)/2}, & \text{$k$ even} \end{cases} \]
  and $J \in \GL(\frac{n-1}{2},\bR)$, $K \in \GL(q-1,\bR)$, and $L \in \GL(\frac{Q-q+1}{2},\bR)$ denoting blocks of the form
\[
    J = \begin{pmatrix} && 1 \\ & \reflectbox{$\ddots$} \\ 1 \end{pmatrix} \quad
    K = \begin{pmatrix} &&& -1 \\ && 1 \\ & -1 \\ \reflectbox{$\ddots$} \end{pmatrix} \quad
    L = (-1)^{Q-1} \begin{pmatrix} && 1 \\ & \reflectbox{$\ddots$} \\ 1 \end{pmatrix}
  \]
  where $q=\min(k,n-k)$ and $Q=\max(k,n-k)$.
\end{Prop}
%
%

\begin{Prf}
  To simplify the description of the limit map, we will first modify the block embedding. For any $\lambda > 1$ the map $b_k$ maps $\begin{pmatrix} \lambda & 0 \\ 0 & \lambda^{-1} \end{pmatrix}$ to
  \[ g_\lambda = \begin{pmatrix} \lambda^{k-1} \\ & \lambda^{k-3} \\ && \ddots \\ &&& \lambda^{1-k} \\ &&&& \lambda^{n-k-1} \\ &&&&& \ddots \\ &&&&&& \lambda^{k+1-n} \end{pmatrix}. \]
  Let $z \in \SO(n)$ be the permutation matrix (or its negative) such that the entries of $z g_\lambda z^{-1}$ are in decreasing order, and consider
  \[ \rho' = z\rho_k^\varepsilon z^{-1} = \iota \circ \overline\rho^\varepsilon, \]
  where $\iota$ is the composition of $b_k$ and conjugation by $z$. The representation $\rho'$ is $B$--Anosov if and only $\rho_k^\varepsilon$ is, and, since $\SO(n) \subset \SL(n,\bR)$ is connected, $\rho'$ then lies in the same component of $\Hom_{\text{$B$--Anosov}}(\Gamma,\SL(n,\bR))$ as $\rho_k^\varepsilon$. So we can consider $\rho'$ instead of $\rho_k^\varepsilon$.

  We first show that $\rho'$ is $B$--Anosov. By \cite[Theorem 8.4]{BochiPotrieSambarino}, it suffices to show that there exist positive constants $c,d$ such that for every $\alpha \in \Delta$ and every element $\gamma \in \Gamma$, we have
  \begin{equation}
    \label{eq:lower_cli}
    \alpha(\mu_n(\rho'(\gamma))) \geq c|\gamma| - d,
  \end{equation}

  where $\mu_n$ is the Cartan projection in $\SL(n,\bR)$ and $|\cdot |$ denotes the word length in $\Gamma$. It follows from the description in \autoref{sec:Hitchin} that $\iota$ maps $\SO(2)$ into $\SO(n)$, and it maps $\begin{pmatrix} \lambda \\ & \lambda^{-1} \end{pmatrix}$ to $zg_\lambda z^{-1}$. Let $\alpha_0$ denote the (unique) simple root for $\SL(2,\bR)$ and $\alpha_i$ the $i$--th simple root for $\SL(n,\bR)$. Then by the above, for $h \in \SL(2,\bR)$,
  \[\alpha_i(\mu_n(\iota(h))) = \begin{cases}\frac{1}{2}\alpha_0(\mu_2(h)) & \text{if $\frac{n+1}{2} - q \leq i \leq \frac{n-1}{2} + q$,} \\ \alpha_0(\mu_2(h)) & \text{otherwise}\end{cases}\]
  Since $\overline\rho^\varepsilon$ is Fuchsian and therefore Anosov, there are positive constants $c_0,d_0$ such that
  \[ \alpha_0(\mu_2(\overline\rho^\varepsilon(\gamma))) \geq c_0|\gamma| - d_0 \qquad \forall \gamma \in \Gamma. \]
  This implies \eqref{eq:lower_cli} with $c = c_0/2$ and $d = d_0/2$, so $\rho'$ is $B$--Anosov.

  The map $\iota$ maps $B_0^{\SL(2,\bR)}$ into $B_0^{\SL(n,\bR)}$ and $-1$ to some diagonal matrix $m = \iota(-1) \in \Mbar$ with $\pm 1$ entries. So $\iota(B) \subset B_0 \cup mB_0 = P_{\{1,m\}}$ and $\iota$ therefore induces smooth maps
  \[\varphi \colon \RP^1 \to \F_{\{1,m\}} = G/P_{\{1,m\}}, \qquad \psi \colon S^1 \to \F_{\{1\}} = G/B_0\]
  which are $\iota$--equivariant. Let $\pi \colon \F_{\{1,m\}} \to \F_\Mbar$ be the projection which forgets orientations, and let $\overline\xi \colon \bdry \to \RP^1$ be the limit map of $\overline\rho^\varepsilon$, a homeomorphism (which does not depend on $\varepsilon$). Then the curve $\xi = \pi \circ \varphi \circ \overline\xi \colon \bdry \to \F_\Mbar$ is $\rho'$--equivariant and continuous. The definition of $z$ ensures that it is also dynamics--preserving. So by \cite[Remark 2.32b]{GueritaudGuichardKasselWienhard} $\xi$ is the limit map of $\rho'$, and $\widehat \xi = \varphi \circ \overline\xi$ is a continuous and equivariant lift to $\F_{\{1,m\}}$.

  We now show that $\xi$ does not lift to $\F_{\{1\}}^{\SL(n,\bR)}$. Write $\pi' \colon \F_{\{1\}} \to \F_\Mbar$ and $p \colon S^1 \to \RP^1$ for the projections. Then $\pi \circ \varphi \circ p = \pi' \circ \psi$. Now assume that $\xi$ lifts to $\F_{\{1\}}$. Then the curve $\xi \circ \overline\xi^{-1} = \pi \circ \varphi \colon \RP^1 \to \F_\Mbar$ also lifts to some curve $\widehat\varphi \colon \RP^1 \to \F_{\{1\}}$, i.e. $\pi' \circ \widehat\varphi = \pi \circ \varphi$. So
  \[\pi' \circ \widehat\varphi \circ p = \pi \circ \varphi \circ p = \pi' \circ \psi.\]
  By right--multiplication with an element of $\Mbar$ we can assume that $\widehat\varphi([1]) = [1]$. So $\widehat\varphi(p([1])) = [1] = \psi([1])$, and uniqueness of lifts implies that $\widehat\varphi \circ p = \psi$. But $p([1]) = p([-1])$, so then $[1] = \psi([1]) = \psi([-1]) = [m] \in \F_{\{1\}}$, which is false since either $k$ or $n-k$ has to be even and therefore $m \in \Mbar \setminus \{1\}$.

  To calculate the transversality type, let $w = \begin{psmallmatrix} 0 & 1 \\ -1 & 0 \end{psmallmatrix} \in \SL(2,\bR)$. Then $\iota(w) \in \widetilde W$ and we can easily compute the relative position of $\widehat \xi$ at the points $x = \overline\xi^{-1}([1])$ and $y = \overline\xi^{-1}([w])$. It is
  \[\pos(\widehat\xi(x), \widehat\xi(y)) = \pos(\varphi([1]),\varphi([w])) = \pos([1], [\iota(w)]) = \dc{\iota(w)}.\]
  Now $\iota_k(w)$ and $\iota_{n-k}(w)$ are antidiagonal, with alternating $\pm 1$ entries and starting with $+1$ in the upper right corner. Conjugation by $z$ interlaces the two blocks in the following way: The resulting matrix is antidiagonal, the middle entry is assigned to the odd--sized block and going towards the corners from there, entries are assigned alternatingly to the two blocks for as long as possible. Combined with the remarks on conjugation by $\Mbar$ at the beginning of \autoref{sec:SLnR} and careful bookkeeping, this proves the claim about the transversality type $\dc{w_k}$. Since $\iota(w)^2 = \iota(-1) = m$, we have $w_k^2 = m$ (recall from \autoref{rem:boundary_maps} \itemnr{2} that $w_k$ is well--defined up to conjugation with $\Mbar$, which does not change the square since $\Mbar$ is abelian).
\end{Prf}
\begin{Cor}\label{cor:number_of_components}
  Let $n$ be odd, $0 \leq k_1 \leq k_2 \leq \frac{n-1}{2}$ and $\rho_{k_1}^{\varepsilon_1}, \rho_{k_2}^{\varepsilon_2}$ be as in the previous proposition. If $\rho_{k_1}^{\varepsilon_1}$ and $\rho_{k_2}^{\varepsilon_2}$ are in the same connected component of $\Hom_{\text{\upshape $B$--Anosov}}(\Gamma,\SL(n,\bR))$, then $k_1 = k_2$ and either $k_1 = k_2 = 0$ or $\varepsilon_1 = \varepsilon_2$.

  As a consequence, $\Hom_{\text{\upshape $B$--Anosov}}(\Gamma,\SL(n,\bR))$ has at least $2^{2g-1}(n-1)+1$ components.
\end{Cor}
\begin{Prf}
  We saw before that $\rho_k^\varepsilon$ is $P_{\{1,w_k^2 \}}$--Anosov, with a limit map of transversality type $\dc{w_k}$, and that this is the minimal oriented parabolic for which $\rho_k^\varepsilon$ is Anosov. By \autoref{cor:connected_components}, if $\rho_{k_1}^{\varepsilon_1}$ and $\rho_{k_2}^{\varepsilon_2}$ were in the same connected component, then $\dc{w_{k_1}}$ and $\dc{w_{k_2}}$ would be conjugate by $\Mbar$, which only occurs when if $k_1 = k_2$ by the discussion at the beginning of \autoref{sec:SLnR}.

  Now assume that $k_1 = k_2 = k \neq 0$ but $\varepsilon_1(\gamma) \neq \varepsilon_2(\gamma)$ for some $\gamma \in \Gamma$. Then $\overline\rho^{\varepsilon_1}(\gamma) = - \overline\rho^{\varepsilon_2}(\gamma)$, so one of them, say $\overline\rho^{\varepsilon_1}(\gamma)$, has two negative eigenvalues while the eigenvalues of $\overline\rho^{\varepsilon_2}(\gamma)$ are both positive. Then $\rho_{k}^{\varepsilon_1}(\gamma)$ has $k$ (if $k$ is even) or $n-k$ (if $k$ is odd) negative eigenvalues, while $\rho_{k}^{\varepsilon_2}(\gamma)$ has only positive eigenvalues. But since $\rho(\gamma)$ has only real non--zero eigenvalues for every $B$--Anosov representation $\rho$, there can be no continuous path from $\rho_{k}^{\varepsilon_1}(\gamma)$ to  $\rho_{k}^{\varepsilon_2}(\gamma)$ in this case.

  In summary, we have $\frac{n-1}{2}$ different possible non--zero values for $k$ and $2^{2g}$ different choices for $\varepsilon$ (its values on the generators of $\Gamma$), giving $2^{2g-1}(n-1)$ connected components, plus the Hitchin component, $k = 0$.
\end{Prf}
\appendix

\section{The \texorpdfstring{$B_0 \times B_0$--action}{B0xB0-action} on \texorpdfstring{$G$}{G}}\label{sec:B0B0action}
\subsection{Orbits}\label{sec:doublecosets}
Let $B \subset G$ be the minimal parabolic subgroup, and $B_0$ its identity component, i.e. the connected subgroup of $\fg$ with Lie algebra $\fb = \bigoplus_{\alpha \in \Sigma^+ \cup \{0\}} \fg_\alpha$. We consider the action of $B_0 \times B_0$ on $G$ by left-- and right--multiplication. In this section, we give a proof of the refinend Bruhat decomposition (\autoref{prop:bruhat}), which shows that the double quotient $B_0 \backslash G / B_0$ is described by the group $\widetilde W = N_K(\fa)/Z_K(\fa)_0$. It will be an important ingredient for our description of relative positions of oriented flags. The proof requires some rather technical preparations.

As restricted root systems are not necessarily reduced, we will work with the set $\Sigma_0$ of indivisible roots, i.e. the roots $\alpha \in \Sigma$ such that $\alpha/2 \not \in \Sigma$. For any $\alpha \in \Sigma_0^+ = \Sigma_0 \cap \Sigma^+$ let $\fu_\alpha = \fg_\alpha \oplus \fg_{2\alpha}$. Then $\fu_\alpha$ is a subalgebra of $\fg$. Let $U_\alpha \subset G$ be the connected subgroup with Lie algebra $\fu_\alpha$.

For $\alpha, \beta \in \Sigma_0^+$ let $(\alpha, \beta) \subset \Sigma_0^+$ be the set of all indivisible roots which can be obtained as positive linear combinations of $\alpha$ and $\beta$. Then $[\fu_\alpha, \fu_\beta] \subset \bigoplus_{\gamma \in (\alpha, \beta)} \fu_\gamma$. For every $w \in W$ the set $\Psi_w = \Sigma_0^+ \cap w\Sigma_0^-$ has the property that $(\alpha, \beta) \subset \Psi_w$ for all $\alpha, \beta \in \Psi_w$. Let $U_w$ be the connected subgroup of $G$ with Lie algebra $\fu_w = \bigoplus_{\alpha \in \Psi_w} \fu_\alpha$.
\begin{Lem}\label{lem:nilpotent_product}
  Let $\Psi' \subset \Psi \subset \Sigma_0^+$ such that $(\alpha, \beta) \subset \Psi$ and $(\alpha, \gamma) \subset \Psi'$ for all $\alpha, \beta \in \Psi$ and $\gamma \in \Psi'$. Let $\fu = \bigoplus_{\alpha \in \Psi} \fu_\alpha$ and $\fu' = \bigoplus_{\alpha \in \Psi'} \fu_\alpha$ and let $U,U' \subset G$ be the corresponding connected subgroups. Let $\Psi \setminus \Psi' = \{\alpha_1, \dots, \alpha_n\}$, in arbitrary order. Then
  \[U = U' U_{\alpha_1} \cdots U_{\alpha_n}.\]
  In particular, for $\Psi = \Psi_w$ and $\Psi' = \varnothing$, we have $U_w = \prod_{\alpha \in \Psi_w} U_\alpha$, where the product can be written in arbitrary order.
\end{Lem}
\begin{Prf}
  First note that the conditions ensure that $\fu, \fu' \subset \fg$ are subalgebras and that $\fu'$ is an ideal of $\fu$. We proceed by induction on $n = |\Psi \setminus \Psi'|$. The case $n = 0$ is trivial and for $n = 1$ the statement is shown in \cite[Lemma 7.97]{Knapp}.

  If $n \geq 1$ then choose a longest root $\alpha_k$ among $\alpha_1, \dots, \alpha_n$ and let $\Psi'' = \Psi' \cup \{\alpha_k\}$. Then $(\alpha, \beta) \subset \Psi' \subset \Psi''$ for all $\alpha \in \Psi$ and $\beta \in \Psi''$, since every element of $(\alpha, \alpha_k)$ will be longer than $\alpha_k$ and therefore in $\Psi'$. So $\fu'' = \bigoplus_{\alpha \in \Psi''}\fu_\alpha$ is an ideal of $\fu$. Let $U'' \subset G$ be its connected subgroup. Since $\fu'$ and $\fu''$ are ideals of $\fu$, $U'$ and $U''$ are normal subgroups of $U$. Now let $g \in U$. By the induction hypothesis there are $g'' \in U''$, $g_- \in U_{\alpha_1}\cdots U_{\alpha_{k-1}}$ and $g_+ \in U_{\alpha_{k+1}}\cdots U_{\alpha_n}$ such that $g = g'' g_- g_+$. Since $U'' \subset U$ is normal, $g = g_- \overline g'' g_+$ for $\overline g'' = g_-^{-1} g'' g_- \in U''$. Since $\Psi' \subset \Psi''$ satisfy the assumptions of the Lemma for $n = 1$, we get $\overline g'' = g' g_0$ for some $g' \in U'$ and $g_0 \in U_{\alpha_k}$. Now set $\overline g' = g_- g' g_-^{-1} \in U'$, then $g = \overline g' g_- g_0 g_+ \in U' U_{\alpha_1}\cdots U_{\alpha_n}$, as required.
\end{Prf}
\begin{Lem}\label{lem:components}
  Let $w \in N_K(\fa)$. Then the map
  \begin{equation} \label{eq:BwB_parametrization}
    U_w \times B \to G, \quad (u,b) \mapsto uwb
  \end{equation}
  is a smooth embedding with image $BwB$. The restriction of \eqref{eq:BwB_parametrization} to $U_w \times B_0$ maps onto $U_w w B_0 = B_0 w B_0$.
\end{Lem}
\begin{Prf}
  We get \eqref{eq:BwB_parametrization} as a composition
  \[U_w \times B \xrightarrow{\mathrm{conj}_{w^{-1}} \times \id} w^{-1} U_w w \times B \hookrightarrow N^- \times B \to G \xrightarrow{L_w} G.\]
  The first and last map are diffeomorphisms, the inclusion is a smooth embedding and the multiplication map $N^- \times B \to G$ is a diffeomorphism onto an open subset of $G$ by \cite[Lemma 6.44, Proposition 7.83(e)]{Knapp}. So the composition is a smooth embedding. We only have to compute its image, i.e. that $U_wwB = BwB$.

  To prove this, use the Iwasawa decomposition $B = NAZ_K(\fa)$ to get $BwB = NwB$ and then write, using \autoref{lem:nilpotent_product},
  \[NwB = \Bigl( \prod_{\alpha \in \Psi_w} U_\alpha \Bigr)\Bigl( \prod_{\alpha \in \Sigma_0^+ \setminus \Psi_w} U_\alpha \Bigr) w B = U_w w \Bigl( \! \prod_{\alpha \in \Sigma_0^+ \setminus \Psi_w} \!\!\! w^{-1} U_\alpha w \Bigr) B.\]
  For all $\alpha \in \Sigma_0^+ \setminus \Psi_w$ we have $w^{-1}\alpha \in \Sigma_0^+$, so $\Ad_{w^{-1}}\fu_\alpha = \fu_{w^{-1}\alpha} \subset \fn$ and thus $w^{-1}U_\alpha w \subset N \subset B$, so $BwB = U_wwB$. If we restrict \eqref{eq:BwB_parametrization} to the connected component $U_w \times B_0$, its image is $U_w w B_0$. The Iwasawa decomposition shows that $B_0 = NAZ_K(\fa)_0$, so $B_0 w B_0 = N w B_0$ and this equals $U_w w B_0$ by the same argument as above.
\end{Prf}
\begin{Prop}[Refined Bruhat decomposition]	\label{prop:bruhat_appendix}
  $G$ decomposes into the disjoint union
  \[G = \bigsqcup_{w \in \widetilde W} B_0 w B_0.\]
\end{Prop}
\begin{Prf}
  Let $\pi \colon \widetilde W \to W$ be the projection to the Weyl group. By the Bruhat decomposition \cite[Theorem 7.40]{Knapp}, $G$ decomposes disjointly into $BwB$ for $w \in W$, so we only have to show that
  \begin{equation}
    \label{eq:bruhat_subcomponents}
    BwB = \bigsqcup_{w' \in \pi^{-1}(w)} B_0 w' B_0.
  \end{equation}
  \Autoref{lem:components} identifies $BwB$ with $U_w \times B$, the connected components of which are the sets $U_w \times mB_0$ for $m \in \Mbar$. These correspond via the map from \autoref{lem:components} to the subsets $U_w w m B_0 = B_0 w m B_0 \subset B w B$. Also $\pi^{-1}(w) = \{wm \mid m \in \Mbar\}$, proving \eqref{eq:bruhat_subcomponents}.
\end{Prf}
\subsection{Orbit closures}\label{sec:doublecosets2}
We now turn to analyzing the closures of $B_0 \times B_0$--orbits in $G$. We call such orbits \emph{refined Bruhat cells}. The main result is \autoref{prop:combinatorial_bruhat_order}, which gives a combinatorial description of closures of refined Bruhat cells. This part is similar to \cite[Section 3]{BorelTits}, where Borel and Tits describe the left and right action of the Borel subgroup for an algebraic group $G$. Some of their arguments also work in our setting.

Let $\widetilde W = N_K(\fa)/Z_K(\fa)_0$ as before and $\pi \colon \widetilde W \to W$ the projection to the Weyl group $W = N_K(\fa)/Z_K(\fa)$. As described in \autoref{sec:parabolics}, $\Delta$ is realized as a generating set of $W$ and we write $\ell$ for the word length with respect to $\Delta$. Recall that, in \autoref{def:generators}, we defined a lift $\I \colon \Delta \to \widetilde W$ of this generating set by $\I(\alpha) = \exp(\frac{\pi}{2}(E_\alpha + \Theta E_\alpha))$, where $E_\alpha \in \fg_\alpha$ were chosen such that $\|E_\alpha\|^2 = 2\|\alpha\|^{-2}$.

Recall that $\Psi_w = \Sigma_0^+ \cap w \Sigma_0^-$.
\begin{Lem}\label{lem:properties_of_Psi}
  Let $w_1, w_2 \in W$. Then
  \[\Psi_{w_1} \cap w_1(\Psi_{w_2}) = \varnothing, \quad \Psi_{w_1w_2} \subset \Psi_{w_1} \cup w_1(\Psi_{w_2}), \quad |\Psi_{w_1}| = \ell(w_1).\]
  Furthermore, $\Psi_w \subset \Span \theta$ for $\theta \subset \Delta$ if and only if $w \in \langle \theta \rangle \subset W$.
\end{Lem}
\begin{Prf}
  First observe that $\alpha \Sigma_0^+ = \Sigma_0^+ \setminus \{\alpha\} \cup \{-\alpha\}$ and therefore $\Psi_\alpha = \{\alpha\}$ for any $\alpha \in \Delta$. The first two identities follow easily from the definition of $\Psi_w$ and the inequality $|\Psi_w| \leq \ell(w)$ is a direct consequence of $\Psi_{w\alpha} \subset \Psi_w \cup w \Psi_\alpha$ for every $\alpha \in \Delta$.

  We want to show that $|\Psi_w| = r$ implies $\ell(w) = r$ by induction on $r \in \bN$. For $r = 0$ this follows from the fact that $W$ acts freely on positive systems of roots. If $|\Psi_w| = r > 0$, then $\Delta \not\subset w\Sigma_0^+$, as otherwise $\Sigma_0^+ \subset w\Sigma_0^+$ and thus $\Psi_w = \varnothing$. So choose $\alpha \in \Delta \cap w\Sigma_0^- \subset \Psi_w$. Then
  \[\alpha\Psi_{\alpha w} = \alpha \Sigma_0^+ \cap w\Sigma_0^- = (\Sigma_0^+ \setminus \{\alpha\} \cup \{-\alpha\}) \cap w\Sigma_0^- = \Psi_w \setminus \{\alpha\},\]
  so $|\Psi_{\alpha w}| = r - 1$ and thus $\ell(\alpha w) = r - 1$ by the induction hypothesis. So $\ell(w) = r$.

  To prove the remaining statement, note that the reflection along a root $\alpha$ maps every other root $\beta$ into $\Span(\alpha, \beta)$. So $\Span(\theta)$ is invariant by every $w \in \langle \theta \rangle$. Assume the equivalence of $\Psi_w \subset \Span \theta$ and $w \in \langle \theta \rangle$ was already proved for $w \in W$ and let $\ell(\alpha w) = \ell(w) + 1$. Then $\Psi_{\alpha w} = \Psi_\alpha \sqcup \alpha \Psi_w = \{\alpha\} \sqcup \alpha \Psi_w$ is contained in $\Span \theta$ if and only if $\alpha \in \theta$ and $w \in \langle \theta \rangle$, proving what we wanted.
\end{Prf}
The two cases $\ord(\I(\alpha)) = 4$ or $\ord(\I(\alpha)) = 2$ in \autoref{rem:choice_of_r} are related to whether there is an associated embedded $\SL(2,\bR)$ or $\PSL(2,\bR)$:
\begin{Lem}\label{lem:SL2_immersion}
  Let $\alpha \in \Delta$ and $E_\alpha$ as in \autoref{def:generators}. Then there is a Lie group homomorphism $\Phi \colon \SL(2,\bR) \to G$, which is an immersion and satisfies
  \begin{enumerate}
  \item $\dd_1\Phi \colon \fsl(2,\bR) \to \fg$ maps $\left(\begin{smallmatrix}0 & 1 \\ 0 & 0\end{smallmatrix}\right)$ to $E_\alpha$, $\left(\begin{smallmatrix}0 & 0 \\ 1 & 0\end{smallmatrix}\right)$ to $-\Theta E_\alpha$, and $\left(\begin{smallmatrix}1 & 0 \\ 0 & -1\end{smallmatrix}\right)$ to $2\|\alpha\|^{-2}H_\alpha$,
  \item $\Phi\left(\begin{smallmatrix}0 & 1 \\ -1 & 0\end{smallmatrix}\right) = \I(\alpha)$,
  \item $\Phi$ is an isomorphism if $\ord(\I(\alpha)) = 4$, and $\ker\Phi = \{\pm 1\}$ if $\ord(\I(\alpha)) = 2$.
  \end{enumerate}
\end{Lem}
\begin{Prf}
  $\dd_1\Phi$ as defined in \itemnr{1} is a monomorphism of Lie algebras \cite[Proposition 6.52]{Knapp}, so it integrates to an immersive Lie group homomorphism $\widetilde\Phi \colon \widetilde{\SL(2,\bR)} \to G$. Since $\ker \widetilde\Phi \subset \widetilde{\SL(2,\bR)}$ is a normal subgroup and discrete and $\widetilde{\SL(2,\bR)}$ is connected, conjugation actually fixes $\ker \widetilde\Phi$ pointwise, i.e.
  \[\ker \widetilde\Phi \subset Z(\widetilde{\SL(2,\bR)}) = \exp(\bZ X), \quad X = \begin{psmallmatrix}0 & \pi \\ -\pi & 0\end{psmallmatrix}.\]
  Now $\widetilde\Phi(\exp(k X)) = \exp(\dd_1\Phi(kX)) = \exp(k\pi(E_\alpha + \Theta E_\alpha)) = \I(\alpha)^{2k}$. If $\ord(\I(\alpha)) = 2$, then $\ker \widetilde\Phi = \exp(\bZ X)$ and if $\ord(\I(\alpha)) = 4$, then $\ker \widetilde\Phi = \exp(2\bZ X)$. By \autoref{rem:choice_of_r}\itemnr{1}, these are the only possibilities. In any case, $\widetilde\Phi$ descends to a homomorphism $\Phi$ on $\SL(2,\bR) = \widetilde{\SL(2,\bR)}/\exp(2\bZ X)$, having the desired properties.
\end{Prf}
We can now already understand the closures of double cosets containing one of the generators, and show how to decompose parabolic subgroups determined by one simple restricted root according to \autoref{prop:bruhat_appendix}.
\begin{Lem}\label{lem:doublecoset_generator}
  Let $\alpha \in \Delta$ and let $s = \I(\alpha) \in \widetilde W$. Then
  \begin{align}
    (P_\alpha)_0 &= B_0 \cup B_0 s B_0 \cup B_0 s^2 B_0 \cup B_0 s^3 B_0, \label{eq:Ps0} \\
    \overline {B_0 s B_0} &= B_0 \cup B_0 s B_0 \cup B_0 s^2 B_0, \label{eq:B0sB0clo} \\
    B_0 s B_0 s^{-1} B_0 &= B_0 \cup B_0 s B_0 \cup B_0 s^{-1} B_0. \label{eq:B0sB0sB0}
  \end{align}
  Note that these unions are disjoint unless $s$ has order 2.
\end{Lem}
\begin{Prf}
  First note that $P_\alpha = B \cup B s B$. This follows from the Bruhat decomposition and the following argument: An element $w \in W$ is contained in $P_\alpha = N_G(\fp_\alpha)$ if and only if $\Ad_w \fp_\alpha \subset \fp_\alpha$. This holds if and only if $w$ preserves $\Sigma_0^+ \cup \Span (\alpha)$, or equivalently $\Psi_w \subset \{\alpha\}$. By \autoref{lem:properties_of_Psi} this is true if and only if $w \in \{1, \alpha \}$.

  We now distinguish two cases, depending on the dimension of $\fu_\alpha$. First assume that $\dim \fu_\alpha > 1$. In this case, $(P_\alpha)_0 \cap B \subset (P_\alpha)_0$ is a closed subgroup of codimension at least 2. Therefore, its complement $(P_\alpha)_0 \cap BsB$ in $(P_\alpha)_0$ is connected, thus equal to $B_0 s B_0$, which is a connected component of $BsB$ by \autoref{lem:components}. This implies $(P_\alpha)_0 \cap Z_K(\fa) = Z_K(\fa)_0$, as otherwise there would be $m \in \Mbar \setminus \{1\}$ with $B_0 s B_0 m = B_0 sm B_0 \subset (P_\alpha)_0$, but this is disjoint from $B_0 s B_0$ by \autoref{prop:bruhat_appendix}. So $(P_\alpha)_0 = B_0 \cup B_0 sB_0$. Since $s \in \widetilde W$ must have order $2$ in this case by \autoref{prop:bruhat_appendix}, this is \eqref{eq:Ps0} as we wanted.

  To see \eqref{eq:B0sB0clo} and \eqref{eq:B0sB0sB0} in this case, we only have to prove that the inclusions $B_0 s B_0 \subset \overline{B_0 s B_0}$ and $B_0 \subset B_0 s B_0 s^{-1} B_0$ are strict: Using $B_0$--invariance from both sides, this will imply $\overline{B_0 s B_0} = B_0 s B_0 s^{-1} B_0 = (P_\Delta \setminus \alpha)_0 = B_0 \cup B_0 s B_0$. And indeed, $B_0 s B_0 = \overline{B_0 s B_0}$ would imply that $B_0 s B_0$ and $B_0$ are closed, so $B_0 s B_0 \cup B_0$ would not be connected. If $B_0 = B_0 s B_0 s^{-1} B_0$, then $sB_0s^{-1}\subset B_0$, so $\fg_{-\alpha} = \Ad_s \fg_\alpha \subset \Ad_s\fb \subset \fb$, a contradiction.

  Now we consider the case $\dim\fu_\alpha = 1$. Then $P_\alpha / B_0$ is a compact 1--dimensional manifold, i.e. a disjoint union of circles. Denote by $\pi$ the projection from $P_\alpha$ to the quotient. Let $e = \begin{psmallmatrix} 0 & 1 \\ 0 & 0\end{psmallmatrix}$ and let $\Phi$ be the map from \autoref{lem:SL2_immersion}. The map $\gamma \colon \bR \to P_\alpha / B_0$ defined by $\gamma(t) = \pi(\Phi(\exp(te))s)$ is an injective smooth curve in this 1--manifold. This is because the map $\bR \to U_\alpha, \ t \mapsto \Phi(\exp(te)) = \exp(tE_\alpha)$ is injective, and the map $U_\alpha \to P_\alpha/B_0, \ u \mapsto \pi(us)$ is injective as a consequence of \autoref{lem:components}. Therefore, its limits for $t \to \pm\infty$ exist and $\overline{\gamma(\bR)} = \gamma(\bR) \cup \{\gamma(\pm\infty)\}$. Now by \autoref{lem:components}
  \[B_0 s B_0 = U_\alpha s B_0 = \exp(\fu_\alpha) s B_0 = \pi^{-1}(\pi(\exp(\fu_\alpha) s)) = \pi^{-1}(\gamma(\bR)),\]
  so
  \[\overline{B_0 s B_0} = \pi^{-1}(\overline{\gamma(\bR)}) = B_0 s B_0 \cup \pi^{-1}(\gamma(\infty)) \cup \pi^{-1}(\gamma(-\infty)).\]
  To compute the limits, note that
  \[\Phi\begin{pmatrix}|t|^{-1} & \sgn(t) \\ 0 & |t|\end{pmatrix} = \exp\!\left[\dd_1\Phi\begin{pmatrix}- \log |t| & 0 \\ 0 & \log |t|\end{pmatrix}\right]\exp\!\left[\dd_1\Phi\begin{pmatrix}0 & t \\ 0 & 0\end{pmatrix}\right] \in B_0,\]
  so
  \begin{align*}
    \lim_{t \to \pm\infty} \gamma(t) &= \lim_{t \to \pm\infty} \pi(\Phi(\exp(te)s)) = \lim_{t \to \pm\infty} \pi \left[ \Phi \left[\begin{pmatrix}1 & t \\ 0 & 1\end{pmatrix}\begin{pmatrix}0 & 1 \\ -1 & 0\end{pmatrix}\begin{pmatrix}|t|^{-1} & \sgn(t) \\ 0 & |t|\end{pmatrix}\right]\right] \\
                                     &= \pi \left[ \Phi \left[\lim_{t\to\pm\infty}\begin{pmatrix}-\sgn(t) & 0 \\ -|t|^{-1} & -\sgn(t)\end{pmatrix}\right] \right] = \pi(\Phi(\mp 1)) = \pi(s^{1 \pm 1}).
  \end{align*}
  So $\overline{B_0 s B_0} = B_0 s B_0 \cup B_0 \cup B_0 s^2 B_0$, which is \eqref{eq:B0sB0clo}.

  Since $B_0 \cup B_0 s B_0 \cup B_0 s^2 B_0 \cup B_0 s^3 B_0 = \overline{B_0 s B_0} \cup \overline{B_0 s^3 B_0} \subset P_\alpha$, $P_\alpha$ decomposes into the disjoint union

  \[P_\alpha = \bigsqcup_{m \in \Mbar} (B_0 \cup B_0 s B_0) \, m = \!\!\! \bigsqcup_{m \in \langle s^2 \rangle \backslash \Mbar} \!\!\! (\overline{B_0 s B_0} \cup \overline{B_0 s^3 B_0}) \, m.\]
  Therefore, $\overline {B_0 s B_0} \cup \overline {B_0 s^3 B_0}$ is closed and open in $P_\alpha$, hence equal to $(P_\alpha)_0$.

  Finally, to prove \eqref{eq:B0sB0sB0}, we claim that, for $t, \tau \in \bR$,
  \begin{equation}
    \label{eq:exp_s_exp_s}
    \pi(\exp(tE_\alpha) s \exp(\tau E_\alpha) s^{-1}) = \begin{cases}\pi(1) & \text{if $\tau = 0$,} \\ \pi(\exp((t-\tau^{-1})E_\alpha)s^{\sgn(\tau)}) & \text{if $\tau \neq 0$.} \end{cases}
  \end{equation}
  This then shows $\pi(U_\alpha s U_\alpha s^{-1}) = \pi(1) \cup \pi(U_\alpha s) \cup \pi(U_\alpha s^{-1})$ and therefore
  \[B_0 s B_0 s^{-1} B_0 = U_\alpha s U_\alpha s^{-1} B_0 = B_0 \cup U_\alpha s B_0 \cup U_\alpha s^{-1} B_0 = B_0 \cup B_0 s B_0 \cup B_0 s^{-1} B_0.\]
  The claim \eqref{eq:exp_s_exp_s} is clear if $\tau = 0$, since $\exp(tE_\alpha) \in B_0$. So let $\tau \neq 0$. Then
  \begin{align*}
    \pi(e^{tE_\alpha}s e^{\tau E_\alpha} s^{-1}) &= \pi \left[\Phi \left[\begin{pmatrix}1 & t \\ 0 & 1\end{pmatrix}\begin{pmatrix}0 & 1 \\ -1 & 0\end{pmatrix}\begin{pmatrix}1 & \tau \\ 0 & 1\end{pmatrix}\begin{pmatrix}0 & -1 \\ 1 & 0\end{pmatrix}\right]\right] \\
                                                 &= \pi \left[ \Phi \left[\begin{pmatrix}1 - t\tau & t \\ -\tau & 1\end{pmatrix}\begin{pmatrix}|\tau|^{-1} & \sgn(\tau) \\ 0 & |\tau|\end{pmatrix}\right]\right] \\
                                                 &= \pi \left[ \Phi \left[\begin{pmatrix}|\tau|^{-1}-t\sgn(\tau) & \sgn(\tau) \\ - \sgn(\tau) & 0\end{pmatrix}\begin{pmatrix} 0 & -\sgn(\tau) \\ \sgn(\tau) & 0\end{pmatrix}\right]s^{\sgn(\tau)}\right] \\
                                                 &= \pi \left[ \Phi \left[\begin{pmatrix} 1 & t - \tau^{-1} \\ 0 & 1\end{pmatrix}\right]s^{\sgn(\tau)}\right] = \pi(\exp((\tau^{-1}-t)E_\alpha)s^{\sgn(\tau)}).\qedhere
  \end{align*}
\end{Prf}
In preparation for the general case, the next three lemmas show how products of refined Bruhat cells behave.
\begin{Lem}\label{lem:doublecoset_additivity1}
  For any $w_1, w_2 \in \widetilde W$ with $\ell(w_1 w_2) = \ell(w_1) + \ell(w_2)$ we have $B_0 w_1 w_2 B_0 = B_0 w_1 B_0 w_2 B_0$.
\end{Lem}
\begin{Prf}
  We get $B_0 w_1 B_0 w_2 B_0 = U_{w_1} w_1 U_{w_2} w_2 B_0$ and $B_0 w_1 w_2 B_0 = U_{w_1w_2} w_1 w_2 B_0$ by \autoref{lem:components}. We want to show that $U_{w_1w_2} = U_{w_1}w_1U_{w_2}w_1^{-1}$. By \autoref{lem:nilpotent_product} both sides can be written as products of $U_\alpha$ for some set of $\alpha$. For the left hand side, the product is taken over all $\alpha \in \Psi_{w_1w_2}$ while for the right hand side we need all $\alpha \in \Psi_{w_1} \cup w_1 \Psi_{w_2}$. But it follows from \autoref{lem:properties_of_Psi} that $\Psi_{w_1w_2} = \Psi_{w_1} \cup w_1 \Psi_{w_2}$ if $\ell(w_1w_2) = \ell(w_1) + \ell(w_2)$.
\end{Prf}
\begin{Lem}\label{lem:doublecoset_additivity2}
  Let $w \in \widetilde W$ and $s = \I(\alpha) \in \widetilde W$ for some $\alpha \in \Delta$. Then $\ell(sw) = \ell(w) \pm 1$ and
  \[ B_0 s B_0 w B_0 =
    \begin{cases}
      B_0 s w B_0 & \text{if $\ell(sw) = \ell(w) + 1$}, \\
      B_0 w B_0 \cup B_0 sw B_0 \cup B_0 s^2w B_0 & \text{if $\ell(sw) = \ell(w) - 1$}.
    \end{cases}\]
\end{Lem}
\begin{Prf}
  Clearly $|\ell(sw) - \ell(w)| \leq 1$ since $\pi(s) \in W$ is in the generating system, but also $\ell(sw) \neq \ell(w)$ by the property of Coxeter groups that only words with an even number of letters can represent the identity. If $\ell(sw) = \ell(w) + 1$, then the statement follows from \autoref{lem:doublecoset_additivity1}. Assume $\ell(sw) = \ell(w) - 1$. Then $\ell(s^{-1}w) = \ell(s^{-2}sw) = \ell(w) - 1$ since $\ell(s^{-2}) = 0$. So $B_0 s B_0 s^{-1}w B_0 = B_0 w B_0$ by the first part and also $B_0 s B_0 s B_0 = B_0 s B_0 \cup B_0 s^2 B_0 \cup B_0 s^3 B_0$ by \autoref{lem:doublecoset_generator}. We thus get
  \begin{align*}
    B_0 s B_0 w B_0 &= B_0 s B_0 s B_0 s^{-1}w B_0 = B_0 s B_0 s^{-1}w B_0 \cup B_0 s^2 B_0 s^{-1}w B_0 \cup B_0 s^3 B_0 s^{-1}w B_0 \\
                    &= B_0 w B_0 \cup B_0 sw B_0 \cup B_0 s^2w B_0,
  \end{align*}
  again using the first part of the lemma for the last equality.
\end{Prf}
\begin{Lem}\label{lem:doublecoset_productclosure}
  Let $w \in \widetilde W$ and $s = \I(\alpha) \in \widetilde W$ for some $\alpha \in \Delta$. Then $\overline {B_0 s B_0 w B_0} = \overline{B_0 s B_0}\;\overline{B_0 w B_0}$ and $\overline{B_0 s B_0} B_0 w B_0 = B_0 w B_0 \cup B_0 sw B_0 \cup B_0 s^2w B_0$.
\end{Lem}
\begin{Prf}
  For the first part, note that $\overline {B_0 s B_0}\;\overline{B_0 w B_0} \subset \overline{B_0 s B_0 w B_0}$. We want to prove that $\overline{B_0 s B_0}\;\overline{B_0 w B_0}$ is closed. Consider the map
  \[f \colon G/B_0 \to \C(G/B_0), \quad gB_0 \mapsto g\overline{B_0 w B_0} \]
  where $\C(G/B_0)$ is the set of closed subsets of $G/B_0$. Since $G/B_0$ is compact and $f$ is $G$--equivariant, the space $\C(G/B_0)$ is compact with the Hausdorff metric, $f$ is continuous and for any closed subset $A \subset G/B_0$ the union $\bigcup_{x \in A} f(x)$ is closed (see e.g. \autoref{prop:hausdorff_compact}, \autoref{lem:properties_of_K}\itemnr{1}, and \autoref{lem:properties_of_K}\itemnr{2}).

  In particular, the union of all elements of $f(\overline{B_0 s B_0}/B_0)$ is a closed subset of $G/B_0$, and so is its preimage in $G$. But this is just $\overline{B_0 s B_0}\;\overline{B_0 w B_0}$, which is therefore a closed set containing $B_0 s B_0 w B_0$, hence equal to $\overline {B_0 s B_0 w B_0}$.

  For the second part, \autoref{lem:doublecoset_generator} implies that
  \[\overline{B_0 s B_0} B_0 w B_0 = B_0 w B_0 \cup B_0 s B_0 w B_0 \cup B_0 s^2 B_0 w B_0\]
  and in both cases of \autoref{lem:doublecoset_additivity2} this equals what we want.
\end{Prf}
We now arrive at the following combinatorial description of closures of refined Bruhat cells.
\begin{Prop}\label{prop:combinatorial_bruhat_order}
  Let $w \in \widetilde W$ and $\pi(w) = \alpha_1 \dots \alpha_k$ a reduced expression by simple root reflections for the projection $\pi(w) \in W$ of $w$ to the Weyl group. Then $w = \I(\alpha_1) \dots \I(\alpha_k) m$ for some $m \in \Mbar$. Let
  \[A_w = \{\I(\alpha_1)^{i_1} \dots \I(\alpha_k)^{i_k}\,m \mid i_1, \dots, i_k \in \{0, 1, 2\}\} \subset \widetilde W.\]
  be the set of words that can be obtained by deleting or squaring some of the letters. Then
  \[\overline{B_0 w B_0} = \bigcup_{w' \in A_w} B_0 w' B_0.\]
  In particular, $A_w$ does not depend on the choice of reduced word for $\pi(w)$.
\end{Prop}
\begin{Prf}
  First of all, since $\pi(\I(\alpha_1)\ldots \I(\alpha_k)) = \alpha_1 \dots \alpha_k = \pi(w)$, there exists $m \in \Mbar$ such that $w = \I(\alpha_1) \dots \I(\alpha_k) m$.\\
  We now prove the second statement by induction on $\ell(w)$. If $\ell(w) = 0$, then $w \in \Mbar$, so $B_0 w B_0 = w B_0$ is already closed, and $A_w = \{w\}$. Now let $\ell(w) > 0$ and assume the statement is already proven for all $\widetilde w \in \widetilde W$ with $\ell(\widetilde w) < \ell(w)$. Assume that $w = \I(\alpha_1)\ldots \I(\alpha_k)m$ as above. Then we can write $w = s \widetilde w$ with $s=\I(\alpha_1)$ and $\ell(\widetilde w) = \ell(w) - 1$. Using \autoref{lem:doublecoset_additivity2} and \autoref{lem:doublecoset_productclosure} we get
  \begin{align*}
    \overline{B_0 w B_0} &= \overline{B_0 s \widetilde w B_0} = \overline{B_0 s B_0 \widetilde w B_0} = \overline{B_0 s B_0} \; \overline{B_0 \widetilde w B_0} = \bigcup_{w' \in A_{\widetilde w}} \overline{B_0 s B_0}\,B_0 w' B_0 \\
                         &= \bigcup_{w' \in A_{\widetilde w}} B_0 w' B_0 \cup B_0 s w' B_0 \cup B_0 s^2 w' B_0 = \bigcup_{w' \in A_w} B_0 w' B_0. \qedhere
  \end{align*}
\end{Prf}
\section{Group actions on compact homogeneous spaces}\label{sec:compact_homogeneous}
The goal of this section is to give a detailed proof of \autoref{prop:expansion_implies_cocompact_appendix}. This is the equivalent of \cite[Proposition 5.30]{KapovichLeebPortiFlagManifolds} in our more general setting. All key arguments of that paper still work.

For a compact metric space $Z$ let $\C(Z)$ be the space of closed subsets, equipped with the Hausdorff metric. The following fact will be useful to us later on (see for example \cite[Lemma 5.31]{BridsonHaefliger} for a proof).\index{C_(Z)@$\C(Z)$ space of closed subsets of $Z$}
\begin{Prop}\label{prop:hausdorff_compact}
  The space $\C(Z)$, equipped with the Hausdorff metric, is a compact metric space.
\end{Prop}
\begin{Def}
  Let $Z$ be a metric space, $g$ a homeomorphism of $Z$ and $\Gamma$ a group acting on $Z$ by homeomorphisms.
  \begin{enumerate}
  \item $g$ is \emph{expanding at $z \in Z$} if there exists an open neighbourhood $z \in U \subset Z$ and a constant $c > 1$ (the expansion factor) such that
    \[d(gx, gy) \geq c\,d(x,y)\]
    for all $x, y \in U$.
  \item Let $A \subset Z$ be a subset. The action of $\Gamma$ on $Z$ is \emph{expanding at $A$} if for every $z \in A$ there is a $\gamma \in \Gamma$ which is expanding at $z$.
  \end{enumerate}
\end{Def}
In \cite{KapovichLeebPortiFlagManifolds}, expansion was used together with a new and more general definition of transverse expansion to prove cocompactness.
\begin{Def}[{\cite[Definition 5.28]{KapovichLeebPortiFlagManifolds}}]
  Let $Z$ be a compact metric space, $g$ a homeomorphism of $Z$ and $\Q \colon \Lambda \to \C(Z)$ a map from any set $\Lambda$.\\
  Then $g$ is \emph{expanding at $z \in Z$ transversely to $\Q$} if there is an open neighbourhood $z \in U$ and a constant $c > 1$ such that
  \[d(gx, g\Q(\lambda)) \geq c\,d(x,\Q(\lambda))\]
  for all $x \in U$ and all $\lambda \in \Lambda$ with $\Q(\lambda) \cap U \neq \varnothing$.
\end{Def}
\begin{Lem}[{\cite[Remark 5.22]{KapovichLeebPortiFlagManifolds}}]	\label{lem:expanding_arbitrarily}
  If the action of $\Gamma$ on $Z$ is expanding at a closed $\Gamma$--invariant subset $A \subset Z$ then it is arbitrarily strongly expanding, i.e. for every $z \in A$ and $c > 1$ there is a $\gamma \in \Gamma$ which is expanding at $z$ with expansion factor $c$.
\end{Lem}
\begin{Prf}
  If the action is expanding at $z \in A$ with some expansion factor, then it is expanding by the same factor in a neighbourhood of $z$. By covering $A$ with finitely many such neighborhoods, we can assume that the action is expanding with a uniform expansion factor $C > 1$. Now let $z \in A$ and let $\gamma_1 \in \Gamma$ be expanding at $z$ by the factor $C$. Let $\gamma_2 \in \Gamma$ be expanding at $\gamma_1z \in A$ by $C$. Then $\gamma_2\gamma_1$ expands at $z$ by $C^2$. Iterating this, we get an element $\gamma_n \cdots \gamma_1 \in \Gamma$ which is expanding with expansion factor $C^n \geq c$ at $z$.
\end{Prf}
Let $G$ be a Lie group and $X$, $Y$ be compact $G$--homogeneous spaces. Fix Riemannian metrics on $X$, $Y$ and a left--invariant Riemannian metric on $G$. Recall that smooth maps between manifolds are locally Lipschitz with respect to any Riemannian distances.
\begin{Lem}\label{lipschitz}
  There exists a compact subset $S \subset G$ such that for every pair $(x,y) \in X^2$, there exists $s_{xy} \in S$ satisfying $s_{xy}x=y$.
\end{Lem}
\begin{Prf}
  Fix a basepoint $x_0 \in X$ and let $V$ be a precompact open neighborhood of the identity in $G$. Since $G \to X, \ g \mapsto gx_0$ is a submersion, $Vx_0$ is a neighborhood of $x_0$. Then by compactness there are finitely many $g_1, \dots, g_n \in G$ such that the sets $g_iVx_0$ cover $X$. So $S = \overline{g_1V} \cup \cdots \cup \overline{g_nV}$ is a compact subset of $G$ which maps $x_0$ to any point in $X$. The set $SS^{-1}$ is compact and satisfies the desired transitivity property.
\end{Prf}
\begin{Lem}	\label{lem:comparable_g}
  There exists a constant $C>0$ such that the following holds: For any two points $x,y \in X$, there exists $g\in G$ satisfying $gx = y$ and $d(1,g) \leq C d(x,y)$.
\end{Lem}
\begin{Prf}
  Assume by contradiction that there are sequences $(x_n),(y_n) \in X^\bN$ such that every $g_n$ sending $x_n$ to $y_n$ must satisfy $d(1,g_n) > n d(x_n,y_n)$. After taking subsequences, we have $x_n \to x, \ y_n \to y$. If $x\neq y$, we obtain in particular that $d(1,g_n) \to \infty$ for every choice of $g_n$ sending $x_n$ to $y_n$. But by \autoref{lipschitz}, a compact subset of $G$ already acts transitively on $X$, so $g_n$ can be chosen such that $d(1,g_n)$ remains bounded. We are thus left with the case $x=y$. Since the map $G \to X, \ g \mapsto gx$ is a smooth submersion, there exists a local section at $x$: There is a neighborhood $x\in U$ and a smooth map $s: U \to G$ satisfying $s(x) = 1$ and $s(z)x = z$ for every $z \in U$. After shrinking $U$ if necessary, $s$ is $C'$--Lipschitz for some $C'>0$. For large $n$, $x_n$ and $y_n$ are inside $U$, and we have $s(y_n)s(x_n)^{-1}x_n = y_n$. Since inversion in $G$ is a smooth map and therefore $C''$--Lipschitz close to the identity, it follows that (after possibly shrinking $U$ some more)
  \[ d\!\left(1,s(y_n)s(x_n)^{-1}\right) = d\!\left(s(y_n)^{-1},s(x_n)^{-1}\right) \leq C'' d(s(y_n),s(x_n)) \leq C'C'' d(y_n,x_n), \]
  a contradiction.
\end{Prf}
\begin{Lem}	\label{lem:lipschitz_on_compact}
  Let $A \subset G$ be a compact set. Then there exists a constant $C>0$ such that:
  \begin{itemize}
  \item The map $A \to X, \ g \mapsto gx$ is $C$--Lipschitz.
  \item For every $g \in A$, the diffeomorphism $X \to X, \ x \mapsto gx$ is $C$--Lipschitz.
  \end{itemize}
\end{Lem}
\begin{Prf}
  The map $G \times X \to X, \ (g,x) \mapsto gx$ is smooth and thus locally Lipschitz. Its restriction to the compact set $A \times X$ is therefore Lipschitz. This implies both parts of the claim.
\end{Prf}
The following auxiliary lemma is a combination of the corresponding statements in Lemmas 7.1, 7.3 and 7.4 in \cite{KapovichLeebPortiFlagManifolds}, transferred to our setting.
\begin{Lem}\label{lem:properties_of_K}
  Let
  \[\Q \colon X \to \mathcal C(Y)\]
  be a $G$--equivariant map. Then there are constants $L, D > 0$ such that
  \begin{enumerate}
  \item $\Q$ is $L$--Lipschitz.
  \item If $A \subset X$ is compact, then $\bigcup_{x \in A} \Q(x)$ is compact.
  \item For all $x \in X$ and $y \in Y$ there exists $x' \in X$ such that $y \in \Q(x')$ and
    \[d(x',x) \leq D\,d(y, \Q(x)).\]
  \item If $\Lambda \subset X$ is compact with $\Q(\lambda) \cap \Q(\lambda') = \varnothing$ for all distinct $\lambda,\lambda' \in \Lambda$, then the map
    \[\pi \colon \bigcup_{\lambda \in \Lambda} \Q(\lambda) \to \Lambda\]
    mapping every point of $\Q(\lambda)$ to $\lambda$ is a uniformly continuous fiber bundle (in the subspace topologies).
  \end{enumerate}
\end{Lem}
\begin{Prf}\ \reallynopagebreak
  \begin{enumerate}
  \item Let $x,y \in X$ be arbitrary. Using \autoref{lem:comparable_g}, we choose $g\in G$ such that $gx = y$ and $d(1,g) \leq Cd(x,y)$. By equivariance, $g\Q(x) = \Q(y)$ holds. As $\mathrm{diam}(X)$ is finite, \autoref{lem:lipschitz_on_compact} implies that $d(gz,hz) \leq C'd(g,h)$ for any $g,h \in \overline{B_{C\cdot \mathrm{diam}(X)}(1)}$ and $z\in Y$. Both constants $C,C'$ do not depend on the choice of $x$ and $y$. Therefore,
    \begin{align*}
      d_H(\Q(x),\Q(y)) & = \max \left\{ \max\limits_{a\in \Q(x)} d(a,g\Q(x)), \max\limits_{gb\in g\Q(x)} d(\Q(x),gb) \right\} \\
                     & \leq \max \left\{ \max\limits_{a\in \Q(x)} \! d(a,ga), \!\!\max\limits_{gb\in g\Q(x)}\!\! d(b,gb) \right\} \leq C' d(1,g) \leq CC' d(x,y).
    \end{align*}
  \item Let $(y_n)$ be a sequence in $\bigcup_{x \in A}\Q(x)$ and $x_n \in A$ such that $y_n \in \Q(x_n)$. Passing to a subsequence we can assume that $y_n \to y \in Y$ and $x_n \to x \in A$. But
    \[d(y_n, \Q(x)) \leq d_H(\Q(x_n),\Q(x)) \to 0\]
    by \itemnr{1}, so $d(y, \Q(x)) = 0$, which means $y \in \Q(x)$ since $\Q(x)$ is closed.
  \item Let $a\in \Q(x)$ be such that $d(y,\Q(x)) = d(y,a)$. By \autoref{lem:comparable_g}, there is an element $g$ with $ga = y$ and $d(1,g) \leq C d(y,a)$. Moreover, since $\mathrm{diam}(Y)$ is finite, \autoref{lem:lipschitz_on_compact} implies that $d(x,gx) \leq C' d(1,g)$. Therefore, $gx$ is the point $x'$ we were looking for.
  \item We start by showing continuity of $\pi$. Assume that $y_n \in \Q(x_n)$, $y_n \to y \in \bigcup_{\lambda \in \Lambda} \Q(\lambda)$ and $\pi(y_n) = x_n \to x \in \Lambda$. We need to show that $\pi(y) = x$. Since $\Q$ is continuous, we have $\Q(x_n) \to \Q(x)$ in $\C(Y)$. Therefore, $d(y,\Q(x)) = 0$, so $y\in \Q(x)$ and $\pi(y)=x$ as $\Q(x)$ is closed. By compactness of $\bigcup_{\lambda \in \Lambda} \Q(\lambda)$ (according to \itemnr{2}), $\pi$ is uniformly continuous.

    Now we construct a local trivialization. Let $x \in\Lambda$ be a point, $U$ a neighborhood of $x$ in $X$, and $s:U \to G$ a smooth local section of the submersion $G \to X, \ g \mapsto gx$. Then the map
    \begin{align*}
      (\Lambda \cap U) \times \Q(x) & \to \bigcup_{\lambda \in \Lambda} \Q(\lambda) \\
      (\lambda,y) & \mapsto s(\lambda)y
    \end{align*}
    is a homeomorphism onto its image, since its inverse is given by
    \begin{align*}
      y &\mapsto (\pi(y),s(\pi(y))^{-1}y). \qedhere
    \end{align*}
  \end{enumerate}
\end{Prf}
The following key lemma shows how expansion in $X$ leads to expansion transverse to the map $\Q \colon X \to \mathcal C(Y)$ in $Y$ (compare \cite[Lemma 7.5]{KapovichLeebPortiFlagManifolds}).
\begin{Lem}	\label{lem:expanding_implies_transverse_expanding}
  Let $\Q \colon X \to \mathcal C(Y)$ be $G$--equivariant and $\Lambda \subset X$ compact with $\Q(\lambda) \cap \Q(\lambda') = \varnothing$ for all distinct $\lambda, \lambda' \in \Lambda$. Let $g \in G$ be expanding at $\lambda \in \Lambda$ with expansion factor $c > LD$, where $L$ and $D$ are the constants from \autoref{lem:properties_of_K}. Then $g$ is expanding at every $y \in \Q(\lambda)$ transversely to $\Q|_\Lambda$.
\end{Lem}
\begin{Prf}
  We give a short outline of the proof before delving into the details. Let $V$ be a neighborhood of the point $y \in \Q(\lambda)$, $y'\in V$ and $\lambda' \in \Lambda$ such that $\Q(\lambda') \cap V \neq \varnothing$. We want to choose $x\in X$ with $y' \in \Q(x)$ such that the following string of inequalities holds:
  \[ d(y',\Q(\lambda')) \leq d_H(\Q(x),\Q(\lambda')) \leq L\,d(x,\lambda') \stackrel{(*)}{\leq} c^{-1}\!L\,d(gx,g\lambda') \stackrel{(**)}{\leq} c^{-1}\!LD\,d(gy',\Q(g\lambda')) \]
  The first two inequalities are true for any choice of $x$. For $(*)$, $x$ and $\lambda'$ need to be close to $\lambda$ so $g$ is expanding. For $(**)$, we need $gx$ to be a ``good'' choice in the sense of \autoref{lem:properties_of_K}\itemnr{3}. Our task is to make sure that the choices of $V$ and $x$ can be made accordingly.

  Since $g$ is expanding, there is an open neighbourhood $U \subset X$ of $\lambda$ with $d(gz,gz') \geq c\,d(z,z')$ for all $z,z' \in U$. We can assume that $U = B_\varepsilon(\lambda)$ for some $\varepsilon > 0$. Let $\pi$ be the bundle map from \autoref{lem:properties_of_K}\itemnr{4}. Since $\pi$ is uniformly continuous there is $\delta > 0$ with
  \begin{align}	\label{eq:projection_continuous}
    d(\pi(y),\pi(y')) < \frac{\varepsilon}{2} \quad \text{whenever} \quad d(y,y') < \delta.
  \end{align}
  We can assume that $\delta \leq \frac{\varepsilon}{2}\alpha\beta D$ where $\alpha$ and $\beta$ are Lipschitz constants for the action of $g^{-1}$ on $X$ and the action of $g$ on $Y$.

  Let $V = B_\delta(y)$ and let $\lambda' \in \Lambda$ with $\Q(\lambda') \cap V \neq \varnothing$ and $y' \in V$. Then by \autoref{lem:properties_of_K}\itemnr{3} (applied to $g\lambda'$ and $gy'$) there is $gx$ such that $gy' \in \Q(gx)$ and thus $y' \in \Q(x)$ and
  \[d(gx, g\lambda') \stackrel{(**)}{\leq} D\,d(gy', \Q(g\lambda')).\]
  Next we want to show that $x,\lambda' \in U = B_\varepsilon(\lambda)$ by bounding $d(\lambda,\lambda')$ and $d(x,\lambda')$. First, since $\Q(\lambda')$ intersects $V$, there is a point $p \in \Q(\lambda')$ with $d(p,y) < \delta$. So $d(\lambda,\lambda') = d(\pi(y),\pi(p)) < \varepsilon/2$ by \eqref{eq:projection_continuous}. Second,
  \[d(x,\lambda') \leq \alpha\,d(gx,g\lambda') \leq \alpha D\,d(gy',g\Q(\lambda')) \leq \alpha\beta D\,d(y',\Q(\lambda')) \leq \alpha\beta D\delta \leq \varepsilon/2,\]
  so $x \in U$ and also $\lambda' \in U$. This implies $(*)$. Therefore,
  \[d(y', \Q(\lambda')) \leq c^{-1}LD\,d(gy', \Q(g\lambda')),\]
  and $c^{-1}LD < 1$, so $g$ is transversely expanding.
\end{Prf}
\begin{Lem}	\label{lem:expansion'_implies_cocompact}
  Let $Y$ be a compact metric space acted upon by a group $\Gamma$ and $\Xi\subset Y$ a compact, $\Gamma$--invariant subset. Assume that for every $y\in\Xi$, there exists a neighborhood $U_y$, an element $\gamma\in\Gamma$ and a constant $c>1$ such that
  \[ d(\gamma y', \Xi) \geq c \, d(y',\Xi) \quad \forall y' \in U. \]
  Then $\Gamma$ acts cocompactly on $Y \setminus \Xi$.
\end{Lem}
\begin{Prf}
  By compactness of $\Xi$, we can find finitely many points $y_1,\ldots,y_n \in Y$ such that their associated neighborhoods $U_{y_i}$ cover $\Xi$. Moreover, there exists a $\delta>0$ such that their union $\bigcup_i U_{y_i}$ contains the $\delta$--neighborhood $N_\delta(\Xi)$. Let $c>1$ be the minimal expansion factor of the corresponding elements $\gamma_i$. We will show that every orbit $\Gamma y, \ y\in Y\setminus \Xi$ has a representative in $Y \setminus N_\delta(\Xi)$. This will prove the lemma since $Y \setminus N_\delta(\Xi)$ is compact.\\
  Indeed, if $y\in N_\delta(\Xi) \setminus \Xi$, there exists $i_0 \in \{1,\ldots,n\}$ such that $y\in U_{y_{i_0}}$. Therefore, $d(\gamma_{i_0}y,\Xi) \geq cd(y,\Xi)$. If $d(\gamma_{i_0}y,\Xi) \geq \delta$, we are done. Else, we repeat the procedure until we obtain a point in the orbit which does not lie in $N_\delta(\Xi)$.
\end{Prf}
After these preparations, we now turn to our goal for this section: A criterion for group actions to be cocompact. Let $\rho \colon \Gamma \to G$ be a representation of a discrete group. This defines an action of $\Gamma$ on $X$ and $Y$. Connecting expansion, transverse expansion and \autoref{lem:expansion'_implies_cocompact} yields the following useful result (compare \cite[Proposition 5.30]{KapovichLeebPortiFlagManifolds}):
\begin{Prop}	\label{prop:expansion_implies_cocompact_appendix}
  Let $\Q \colon X \to \mathcal C(Y)$ be $G$--equivariant and $\Lambda \subset X$ compact and $\Gamma$--invariant with $\Q(\lambda) \cap \Q(\lambda') = \varnothing$ for all distinct $\lambda, \lambda' \in \Lambda$. Also assume that the action of $\Gamma$ on $X$ is expanding at $\Lambda$. Then $\Gamma$ acts cocompactly on $\Omega = Y \setminus \bigcup_{\lambda \in \Lambda} \Q(\lambda)$.
\end{Prop}
\begin{Prf}
  Since $\Lambda$ is closed and $\Gamma$--invariant, we know by \autoref{lem:expanding_arbitrarily} that the action of $\Gamma$ is expanding arbitrarily strongly at every point $\lambda \in \Lambda$. \autoref{lem:expanding_implies_transverse_expanding} therefore shows that the action of $\Gamma$ is expanding at every point of $\bigcup_{\lambda\in\Lambda} \Q(\lambda)$ transversely to $\Q$. We will show that this implies the prerequisites of \autoref{lem:expansion'_implies_cocompact} and thus the action on $\Omega$ is cocompact.\\
  First of all, we observe that for any point $z\in Y$, we have
  \begin{equation}	\label{eq:closest_fiber}
    d\Big(z,\bigcup_{\lambda\in\Lambda} \Q(\lambda) \Big) = d(z,\Q(\lambda_z))
  \end{equation}
  for some $\lambda_z \in \Lambda$ (this follows from compactness of $\bigcup_{\lambda\in\Lambda} \Q(\lambda)$, \autoref{lem:properties_of_K}\itemnr{2}). Now let $y \in \bigcup_{\lambda\in\Lambda} \Q(\lambda)$ and let a neighborhood $U\ni y$ and $\gamma \in \Gamma$ be chosen such that $\gamma$ is $c$--expanding on $U$ transversely to $\Q$. There exists $\epsilon > 0$ satisfying
  \[ B_\epsilon(\gamma y) \subset \gamma U. \]
  Let $\delta > 0$ be sufficiently small such that
  \[ \gamma B_\delta(y) \subset B_{\epsilon/2}(\gamma y). \]
  For any point $y' \in B_\delta(y)$, let $\lambda_{\gamma y'}$ be chosen as in \eqref{eq:closest_fiber}. Since $d(\gamma y', \gamma y) < \epsilon/2$ and $\gamma y \in \bigcup_{\lambda\in\Lambda} \Q(\lambda)$, we necessarily have $\Q(\lambda_{\gamma y'}) \cap B_\epsilon(\gamma y) \neq \varnothing$. Therefore,
  \[ \gamma^{-1}\Q(\lambda_{\gamma y'}) \cap U \neq \varnothing. \]
  Since $y' \in B_\delta(y)\subset U$, transverse expansion now implies
  \[ d\Big(\gamma y',\bigcup_{\lambda\in\Lambda} \Q(\lambda)\Big) = d(\gamma y', \Q(\lambda_{\gamma y'})) \geq c d(y',\gamma^{-1} \Q(\lambda_{\gamma y'})) \geq c d\Big(y', \bigcup_{\lambda\in\Lambda} \Q(\lambda)\Big). \qedhere \]
\end{Prf}

\printindex

\printbibliography

Florian Stecker \\
Heidelberg Institute for Theoretical Studies \\
Mathematics Department, Heidelberg University \\
\emph{E--mail address:} \texttt{fstecker@mathi.uni-heidelberg.de}

Nicolaus Treib \\
Heidelberg Institute for Theoretical Studies \\
Mathematics Department, Heidelberg University \\
\emph{E--mail address:} \texttt{ntreib@mathi.uni-heidelberg.de}

\end{document}